\documentclass[12pt,a4paper]{article}

\usepackage{amsmath}    % need for subequations
\usepackage{amssymb}
\usepackage{graphicx}   % need for figures
\usepackage{verbatim}   % useful for program listings
\usepackage{color}      % use if colour is used in text

\usepackage{hyperref}   % use for hypertext links, including those to external documents and URLs

\usepackage{enumerate}
\usepackage{mathrsfs}
\usepackage{empheq}
\usepackage{bbold}

\usepackage[final]{showlabels}

\usepackage{amsthm}
\newtheorem{definition}{Definition}[section]
\newtheorem{lemma}{Lemma}[section]
\newtheorem{corollary}{Corollary}[section]
\newtheorem{theorem}{Theorem}[section]

\usepackage[makeroom]{cancel}

\usepackage{array}
\newlength\mylen
\newcolumntype{C}{>{\hfil$}p{\mylen}<{$\hfil}} % centered, in math mode, fixed width

\newcommand{\del}{\partial}
\renewcommand{\theta}{\vartheta}
\renewcommand{\phi}{\varphi}
\newcommand{\vecc}[2]{\left ( \begin{array}{c}#1\\#2\\ \end{array}\right )}
\newcommand{\veccc}[3]{\left ( \begin{array}{c}#1\\#2\\#3\\ \end{array}\right )}

\newcommand{\dd}{\mathrm{d}}

\newcommand{\id}{\mathbb{1}}

\renewcommand{\div}{\mathrm{div\,}}
\renewcommand{\vec}{\mathbf}

\newcommand{\ii}{\mathbb{i}}

\renewcommand{\title}{A node-conservative vorticity-preserving Finite Volume method for linear acoustics on unstructured grids}

\newcommand{\authorOne}{Wasilij Barsukow\footnote{Bordeaux Institute of Mathematics, Bordeaux University and CNRS/UMR5251, Talence, 33405 France, wasilij.barsukow@math.u-bordeaux.fr}}
\newcommand{\authorTwo}{Rapha\"el Loubère\footnote{Bordeaux Institute of Mathematics, Bordeaux University and CNRS/UMR5251, Talence, 33405 France, raphael.loubere@u-bordeaux.fr}}
\newcommand{\authorThree}{Pierre-Henri Maire\footnote{CEA/Cesta, 33116 Le Barp, France, pierre-henri.maire@cea.fr}}
\begin{document}

\begin{center} \Large
\title

\vspace{1cm}

\date{}
\normalsize

\authorOne, \authorTwo, \authorThree
\end{center}

\begin{abstract}

% #############################################################################################################################
Instead of ensuring that fluxes across edges add up to zero, we split the edge in two halves and also associate different fluxes to each of its sides. This is possible due to non-standard Riemann solvers with free parameters. We then enforce conservation by making sure that the fluxes around a node sum up to zero, which fixes the value of the free parameter. We demonstrate that for linear acoustics one of the non-standard Riemann solvers leads to a vorticity preserving method on unstructured meshes.
% #############################################################################################################################



\end{abstract}

\section{Introduction}

The aim of this paper is to explore a strategy for deriving structure preserving methods from first principles. We focus on the involutional constraint of stationary vorticity for linear acoustics.

First, the reasons for studying structure preserving methods will be recalled. A standard requirement for numerical methods is that they are convergent, i.e. that the numerical solution converges in some sense to the solution of the PDE as the scales of discretization $\Delta t, \Delta x, \Delta y, \ldots$ in time and space go to zero. In view of limited computational resources (in particular for multi-dimensional simulations) convergence alone is not enough, and the focus must lie on the quality of simulations for finite $\Delta x, \Delta t$. Conventional strategies include improving the speed of convergence (i.e. employing a high-order method), or to perform grid refinement selectively in regions of particular interest (adaptive mesh refinement). Both strategies require additional memory and additional computational time. As is discussed below, however, there exist examples of numerical methods (\textbf{structure preserving methods}) that show drastically better results on coarse grids and which have been obtained by merely changing a few terms, i.e. without any additional memory requirement or additional computational time. These methods are derived following the paradigm that the numerical methods should mimic some key properties of the PDE.

Some examples of the kind of properties that are considered will be given next. 
In the context of the $\nabla \cdot \vec B =0 $ involutional constraint\footnote{An involutional constraint is a constraint that needs to be fulfilled only at initial time. In other words, the evolution equations already imply its stationarity.} of magnetohydrodynamics its violation in a numerical simulation has been linked to numerical instability. Another example of an involutional constraint is found for the equations of linear acoustics
\begin{align}
 \del_t \vec v + \nabla p &= 0\\
 \del_t p + \nabla \cdot \vec v &= 0.
\end{align}
The curl $\nabla \times \vec v$ remains stationary, although $\vec v$ might be evolving. This paper focuses on preserving discrete involutions. Even if the violation of an involutional constraint does not necessarily lead to instability or make the simulation crash, it has been found in practice that preserving a discrete version of the constraint (\textbf{involution preserving methods}) dramatically improves the quality of the results on coarse grids, as discussed next. 

Involution preserving methods are also \textbf{stationarity preserving} (\cite{barsukow17a}), a property that has a large practical impact. If the initial data are a projection of a stationary state of the PDE onto the discrete degrees of freedom, then they in general are not kept stationary by the numerical method. This is true for both stationarity preserving methods and those which do not have this property. As long as the numerical method is von Neumann stable, all instationary modes are decaying exponentially quickly in time. The time evolution of any initial data will be such that it evolves exponentially quickly towards one of the stationary states of the numerical method; more slowly if the grid is finer. This, again, is true for both stationarity preserving methods and those which are not. The difference between them is in the amount of stationary states of the numerical method. 

Numerical methods that are not stationarity preserving have a very poor set of stationary states (e.g. only constant states). For linear acoustics, all divergencefree velocity fields remain stationary (as long as $p$ is constant), i.e. the set of stationary states is much richer than just constants. Initial data close to a stationary state of the PDE will be diffused away and will become a constant in the limit of long time, i.e. after some time will no longer be recognizable. Stationarity preserving numerical methods discretize all stationary states, i.e. for each projection of a stationary state of the PDE there exists a stationary state of the numerical method that is in some sense close. This means that upon time evolution, a projection of a stationary state of the PDE will exponentially quickly settle onto some stationary state of the numerical method, which will be similar; more similar if the grid is finer. After this transition the numerical solution will no longer change in time; a crucial qualitative difference to the exponentially decaying stationary states of a numerical method that fails to be stationarity preserving. The precise definitions and statements involve an analysis of the discrete Fourier modes, and are presented in \cite{barsukow17a} and exemplified in \cite{barsukow18hypproceeding} in great detail.  

Whenever the involutions involve differential operators, the question of their numerical preservation encounters the fundamental difficulty that only \emph{approximations} of the differential operators can be preserved. As can be derived from the above heuristic explanation, the qualitative ability to preserve stationary states is independent of the details of the discrete divergence operator characterizing the discrete stationary states, as long as it is consistent. The same is true for the preservation of involutions. It thus seems to be of great importance that \emph{some} discretization of the involution is preserved at all, and the question which particular one it is decides over more subtle aspects, such as the order of accuracy of the involution approximation.

Next, the reason why standard finite volume methods are not involution preserving will be reviewed in a heuristic way. For stability, finite volume methods add numerical diffusion. In general, also the evolution of discretizations of the involutional constraint will diffuse away in time. If there exists some discretization which is exempt from the effect of diffusion, then such a method is called involution-preserving. An efficient way of finding such discretizations, or proving that no such discretizations exist has been described in \cite{barsukow17a} in the case of linear equations and involutions and employs the discrete Fourier transform.

In the tradeoff between maintaining stability and enlarging the kernel of the operator of numerical diffusion several strategies have been applied. Having no numerical diffusion, i.e. when employing central discretizations, stabilization becomes more difficult. Forward-Euler time integration is no longer possible, but some high-order Runge-Kutta schemes or leap-frog-type time integration can be used. The latter were successfully employed in the derivation of explicit involution-preserving methods in \cite{barsukow21yee}. A popular way is also to resort to implicit, or partially implicit (IMEX) time integrators. 

In this work, explicit, in particular forward-Euler time integration of the method will remain possible, which means that approaches based only on central derivatives cannot be applied here and presence of upwinding/numerical diffusion is absolutely necessary. 

It is sometimes possible to remove only certain terms from the (matrix) of numerical diffusion to achieve involution preservation while maintaining stability. It can be shown that \textit{ad-hoc} modifications commonly applied to finite volume methods in the regime of low Mach number flow, once applied to the equations of linear acoustics, lead to vorticity preservation and allow explicit time integration (\cite{barsukow17a}). For illustration, consider the modified equation for the dimensionally split upwind method for linear acoustics in two spatial dimensions
\begin{align}
 \del_t u + \del_x p &= 0\label{eq:linacu2d}\\
 \del_t v + \del_y p &= 0\label{eq:linacv2d}\\
 \del_t p + \del_x u + \del_y v &= 0 \label{eq:linacp2d}
\end{align}
which reads
\begin{align}
 \del_t u + \del_x p &= \Delta x \del_x^2 u \\
 \del_t v + \del_y p &= \Delta y \del_y^2 v \\
 \del_t p + \del_x u + \del_y v &= \Delta x \del_x^2 p + \Delta y \del_y^2 p.
\end{align}
The problematic terms are those on the right hand side of the $u$ and $v$ equations. On the one hand, they diffuse away divergence-free ($\del_x u + \del_y v = 0$) initial data instead of keeping them stationary. On the other hand, $\vecc{\Delta x \del_x^2 u}{\Delta y \del_y^2 v}$ is not a gradient, such that the curl $\del_y u - \del_x v$ is no longer an involution: $\del_t (\del_y u - \del_x v) \neq 0$. This can be rectified by removing these terms. Based on the modified equation alone one cannot judge about the stability of the resulting method, but there exist numerical methods that are stable upon centering the $u$, $v$ equations (\cite{dellacherie10,barsukow16}).

The alternative path, and the one pursued in this work, is to complete these terms such that they become a gradient of the divergence, see \eqref{eq:multidmodifiedu}--\eqref{eq:multidmodifiedv}:
\begin{align}
 \del_t u + \del_x p &= \Delta x \del_x(\del_x u + \del_y v) \label{eq:multidmodifiedu}\\
 \del_t v + \del_y p &= \Delta y \del_y(\del_x u + \del_y v)\label{eq:multidmodifiedv}\\
 \del_t p + \del_x u + \del_y v &= \Delta x \del_x^2 p + \Delta y \del_y^2 p. \label{eq:multidmodifiedp}
\end{align}
Equation \eqref{eq:multidmodifiedp} has been left unmodified. Observe that the new terms $\del_x \del_y v$ and $\del_x \del_y u$ are cross-derivatives which are not available in dimensionally split methods: the method whose modified equations are \eqref{eq:multidmodifiedu}--\eqref{eq:multidmodifiedp} must be truly multi-dimensional. It is important to mention that one cannot just use \emph{any} set of discretizations of the derivatives appearing in \eqref{eq:multidmodifiedu}--\eqref{eq:multidmodifiedp} . It turns out (see \cite{barsukow17a} for more details) that one needs to find discrete second derivatives $D_{(2,0)} \simeq \del_x^2$, $D_{(1,1)} \simeq \del_x \del_y$, $D_{(0,2)} \simeq \del_y^2$ and discrete first derivatives $$D_{(1,0)} \simeq \del_x, D_{(0,1)} \simeq \del_y$$ such that for any grid functions $u, v$ the statement 
\begin{align}
 D_{(1,0)} u + D_{(0,1)} v = 0 \quad \text{everywhere on the grid}
\end{align}
implies
\begin{align}
 D_{(2,0)} u + D_{(1,1)} v &= 0\\
 D_{(1, 1)} u + D_{(0,2)} v &= 0.
\end{align}
To this end, the different discretizations have to ``match'', and it is \textit{a priori} not clear whether such tuples of discrete derivatives exist at all. On a Cartesian grid one can show that choosing central derivatives 
\begin{align}
 (D_{(1,0)} u)_{ij} &= \frac{u_{i+1,j} - u_{i-1,j}}{2\Delta x},  & (D_{(0,1)} v)_{ij} &= \frac{v_{i,j+1} - v_{i,j-1}}{2 \Delta y},
\end{align}
does not allow to find suitable $D_{(2,0)}$, $D_{(1,1)}, D_{(0, 2)}$, for example. In fact, if attention is restricted to $3\times 3$ stencils on Cartesian grids and if one insists on all the discrete derivatives being (anti)symmetric, then there is a unique choice of discrete derivatives (see \cite{barsukow17a} for a proof):

\settowidth\mylen{-11}

\begin{align}
 8\Delta x D_{(1,0)} &= \begin{array}{C|C|C} -1 & 0 & 1 \\\hline -2 & 0 & 2 \\\hline -1 & 0 & 1  \end{array} ,&
 8\Delta y D_{(0,1)} &= \begin{array}{C|C|C} 1 & 2 & 1 \\\hline 0 & 0 & 0 \\\hline -1 & -2 & -1  \end{array}, \label{eq:cartesainstencils1}\\
 4\Delta x^2 D_{(2,0)} &= \begin{array}{C|C|C} 1 & -2 & 1 \\\hline 2 & -4 & 2 \\\hline 1 & -2 & 1  \end{array}, &
 4\Delta y^2 D_{(0,2)} &= \begin{array}{C|C|C} 1 & 2 & 1 \\\hline -2 & -4 & -2 \\\hline 1 & 2 & 1  \end{array}, \\
 4\Delta x \Delta y D_{(1,1)} &= \begin{array}{C|C|C} -1 & 0 & 1 \\\hline 0 & 0 & 0 \\\hline 1 & 0 & -1  \end{array} .\label{eq:cartesainstencils2}
\end{align}
This uniqueness result explains why many vorticity preserving numerical methods that have been suggested in the literature (\cite{morton01,sidilkover02,jeltsch06,mishra09preprint}) essentially use the same finite differences on Cartesian grids. One also observes that now, the stationary states of \eqref{eq:multidmodifiedu}--\eqref{eq:multidmodifiedp} are again given by constant $p$ and the vanishing divergence $\del_x u + \del_y v$. As shown in \cite{barsukow17a}, a vorticity preserving method is also stationarity preserving, i.e. its stationary states are a discretization of all the stationary states of the PDE (and not just of a subset of those).

Both strategies, removing troublesome terms or adding new terms, are not first-principle derivations, although the latter strategy, at least, comes at the advantage of not modifying the one-dimensional method. Now, Godunov's method is a first-principle derivation of a stable numerical method. However, the application of the one-dimensional Godunov method in a dimensionally-split way does not yield involution preservation for linear acoustics. One might object that, as the involution is a truly multi-dimensional object (the curl being trivial in 1D), one should rather consider the fully multi-dimensional Godunov method that takes into account all the multi-dimensional Riemann problems. Despite its complexity, such a truly multi-dimensional solver for linear acoustics has been derived and studied in \cite{barsukow17}, and, suprisingly, no improvement in terms of structure preservation has been found. As the three steps of Godunov's procedure, reconstruction--evolution--average, have been made without any approximations, a structure preserving method coming from a first-principle derivation cannot be a Godunov method in the usual sense.

The aim of the paper is to show how, based on previous work \cite{gallice22} and the new concept of conservation around a node, a truly multi-dimensional numerical method with modified equation \eqref{eq:multidmodifiedu}--\eqref{eq:multidmodifiedp} can be derived without \textit{a posteriori} modifications of terms. The new strategy yields a structure-preserving method on general grids. On Cartesian grids it reduces to a method that again uses the well-known discrete derivatives \eqref{eq:cartesainstencils1}--\eqref{eq:cartesainstencils2}. 

The paper is organized as follows. After introducing a new notion of conservation in Section \ref{sec:conservation}, an approximate Riemann solver is constructed in Section \ref{sec:riemann}, which contains a free parameter. In Section \ref{sec:method} the new method is constructed on unstructured grids and in Section \ref{sec:statioanritypreservation} it is shown that upon a suitable choice of the parameter in the Riemann solver the resulting method is involution preserving. Section \ref{sec:cartesian} is reviewing both the concepts and the new vorticity preserving method on the special case of Cartesian grids. The method is extended to second-order accuracy in Section \ref{sec:secondorder}. The time integration is discussed in Section \ref{sec:timeintegration}. Numerical results are presented in Section \ref{sec:numerical}.

We set variables with as many components as the dimensions of space in boldface. Indices never denote derivatives. $\vec a \times \vec b$ ($\vec a \cdot \vec b$) is the cross product (scalar product) of the two vectors $\vec a$ and $\vec b$ and $\vec a \otimes \vec b$ is the matrix with entries $a_i b_j$.

\section{Conservative methods} \label{sec:conservation}

\subsection{Conservation at PDE level}

\begin{figure}
 \centering
 \includegraphics[width=0.6\textwidth]{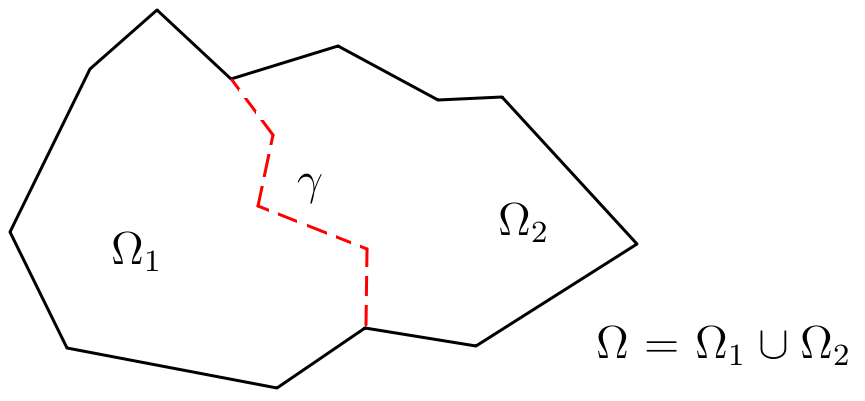}
 \caption{Illustration of conservation: The change in mass in $\Omega_1$ due to the flux through $\gamma$ (dashed) is the negative of the change of mass in $\Omega_2$.}
 \label{fig:internaledge}
\end{figure}

Consider the conservation law 
\begin{align}
\del_t q + \nabla \cdot \vec f(q) &= 0 & q &\colon \mathbb R^+_0 \times \mathbb R^2 \to \mathbb R^m \\
 \vec f = (f_1, f_2), \qquad f_i &\colon \mathbb R^m \to \mathbb R^m.
\end{align}
The update of the average $q_\Omega := \frac{1}{|\Omega|} \int_\Omega  q\, \dd \vec x $ of $q$ over a compact set $\Omega \subset \mathbb R^2$ with piecewise differentiable boundary $\del \Omega$ reads
\begin{align}
 \del_t q_\Omega = - \frac{1}{|\Omega|} \oint_{\del \Omega} \vec f \cdot \vec n  \, \dd s =  - \frac{1}{|\Omega|} \sum_{e \subset \del \Omega} |e| \hat f_{e,\Omega}. \label{eq:finvol}
\end{align}
where, in the case of a polygonal domain $\Omega$, $\hat f_{e,\Omega} = \frac{1}{|e|} \int_e \vec f \cdot \vec n \,\dd s$ is the average of $\vec f \cdot \vec n$ with $\vec n$ the outward unit vector associated to the straight line segment $e$ of $\del \Omega$. Assume that $\Omega = \Omega_1 \cup \Omega_2$ such that $\gamma := \overline \Omega_1 \cap \overline \Omega_2$ is a piecewise straight curve (see Figure \ref{fig:internaledge}), and assume that the only non-vanishing fluxes are those through $\gamma$. Then the total mass does not change: $\del_t (|\Omega_1| q_{\Omega_1} + |\Omega_2|q_{\Omega_2}) = 0$, while according to Equation \eqref{eq:finvol}

\begin{align}
 \del_t q_{\Omega_1} &= - \frac{1}{|\Omega_1|} \sum_{e \subset \gamma} |e| \hat f_{e,\Omega_1}, &
 \del_t q_{\Omega_2} &= - \frac{1}{|\Omega_2|} \sum_{e \subset \gamma} |e| \hat f_{e,\Omega_2} .
\end{align}

Therefore conservation means 
\begin{align}
 -  \sum_{e \subset \gamma} |e| \hat f_{e,\Omega_1} = \sum_{e \subset \gamma} |e| \hat f_{e,\Omega_2}. \label{eq:conservationconditionacrossedgegeneral}
\end{align}
The mass flowing from cell $\Omega_1$ to $\Omega_2$ through $\gamma$ is the negative of the mass flowing from $\Omega_2$ to $\Omega_1$. To this end, it is sufficient to have
\begin{align}
 f_{e,\Omega_1}  = -f_{e,\Omega_2}  \qquad \forall e \subset \gamma. \label{eq:conservationconditionacrossedgegeneral2}
\end{align}
The expression on the left side is the flux through $e$ as seen from $\Omega_1$ (with the normal pointing out of $\Omega_1$), while the one on the right is the flux through $e$ as seen from $\Omega_2$ (with the normal thus pointing in the opposite direction).

Global conservation means that the total mass inside some volume (e.g. $\Omega$) is constant up to fluxes through the boundaries. Local conservation is the same statement for smaller subvolumes (e.g. $\Omega_1$ and $\Omega_2$).

\subsection{The classical and the new concepts of conservation}

Consider now numerical methods on a computational grid consisting of polygonal cells $c$. We do not make a notational distinction between the cell $c$ as a subset of $\mathbb R^d$ (which would allow to write $\del c$ for its boundary) and index $c$ of a cell (which allows to write $c \in \{ c_1, c_2, c_3 \}$ for any of the cells indexed by $c_1, c_2, c_3$). We treat similarly other objects, such as nodes and edges.

\begin{definition}
 Denote
 \begin{enumerate}[(i)]
  \item by $|e|$ the length of an edge and by $|c|$ the area of a cell;
  \item by $\mathcal C$ the set of all cells, by $\mathcal E$ the set of all edges and by $\mathcal N$ the set of all nodes of the grid;
  \item by $\mathcal C(n)$ the set of all cells around node $n$, and by $\mathcal N(c)$ the set of all nodes contained in the boundary of cell $c$;
  \item by $\mathcal E(c)$ the set of all edges that form the boundary of cell $c$ and by $\mathcal E(n)$ the set of all edges containing node $n$;
  \item by $\mathcal N(e)$ the set of (two) nodes at the boundaries of an edge $e$ and by $\mathcal C(e)$ the set of (two) cells adjacent to an edge $e$;
 \end{enumerate}
\end{definition}

In the context of a numerical method, the cells are the smallest volumes under consideration. Denote by $q_c$ the average of $q$ over a cell $c$. By applying \eqref{eq:finvol} with $\Omega = c$ one is led to define a classical finite volume method as
\begin{align}
 \del_t q_c &= - \frac{1}{|c|} \sum_{e \in \mathcal E(c)} |e| \hat f_{e,c} \qquad \forall c \in \mathcal C
\end{align}
where $\hat f_{e,c}$ is the numerical flux associated to an edge $e$ as seen by cell $c$. Global conservation means that
\begin{align}
 0 = \sum_{c \in \mathcal C} |c| \del_t q_c = - \sum_{c \in \mathcal C} \sum_{e \in \mathcal E(c)} |e| \hat f_{e,c} &= \sum_{e  \in \mathcal E} |e| \sum_{c \in \mathcal C(e)} \hat f_{e, c}. \label{eq:globalconservationtonodal}
\end{align}

Observe how all the fluxes on the grid have been grouped by edge in \eqref{eq:globalconservationtonodal}, thus allowing for a passage from global to local conservation. Usually, global conservation is ensured by imposing 
\begin{align} \sum_{c \in \mathcal C(e)} \hat f_{e, c} = 0 \end{align}
for all edges $e \in \mathcal E$, which is the same as \eqref{eq:conservationconditionacrossedgegeneral2}. This classical concept will be referred to as \emph{edge-conservation}. It is also possible to regroup the sum differently, and here the focus will be on the nodes.

\begin{figure}
 \centering
 \includegraphics[width=0.6\textwidth]{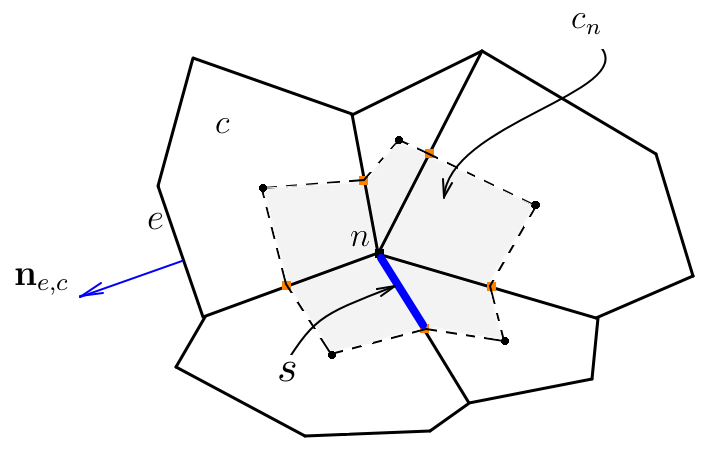} 
 \caption{The dual cell $c_n$ associated to a node $n$ (bounded by the dashed line that joins cell centroids and edge midpoints).}
 \label{fig:celldual}
\end{figure}

\begin{definition}
 
 \begin{enumerate}[(i)]
  \item Denote by $c_n$ the dual cell associated to node $n$: the polygonal region whose vertices are the centroids of all cells $c \in \mathcal C(n)$ and the edge midpoints for all the edges $e \in \mathcal E(n)$ (see Figure \ref{fig:celldual}); 
  \item refer to the halved edges $s$ that appear in the construction of $c_n$ as \emph{subedges}, i.e. $s = c_n \cap e$ for some $e \in \mathcal E(n)$; 
  \item denote by $\mathcal S \mathcal E(n)$ the set of subedges around a node $n$ and by $\mathcal S \mathcal E(c)$ the set of subedges forming the boundary of a cell $c$;
  \item denote finally by $\mathcal S \mathcal E(n, c) = \mathcal S \mathcal E(n) \cap \mathcal S \mathcal E(c)$ the set of (two) subedges of a cell $c$ that are adjacent to one of its nodes $n \in \mathcal N(c)$.
  \end{enumerate}
\end{definition}

The last part of the Definition implies that the summation $\sum_{s \in \mathcal S \mathcal E(c)}$ can equivalently be replaced by $\sum_{n \in \mathcal N(c)} \sum_{s \in \mathcal {SE}(n, c)}$.

The natural definition of a finite volume method now is
\begin{align}
 \del_t q_c &= - \frac{1}{|c|} \sum_{s \in \mathcal S \mathcal E(c)} |s| \hat f_{s,c} = - \frac{1}{|c|} \sum_{n \in \mathcal N(c)} \sum_{s \in \mathcal {SE}(n, c)} |s| \hat f_{s,c}  \label{eq:finvolsubedges}
\end{align}
where $\hat f_{s,c}$ is the numerical flux associated to a subedge $s$ as seen by cell $c$. Global conservation can be ensured by grouping the fluxes by \emph{node}
\begin{align}
 0 = \sum_{c \in \mathcal C} |c| \del_t q_c = - \sum_{c \in \mathcal C} \sum_{s \in \mathcal S \mathcal E(c)} |s| \hat f_{s,c} &= \sum_{n \in \mathcal N} \sum_{c \in \mathcal C(n)} \sum_{s \in \mathcal {SE}(n,c)}  |s| \hat f_{s, c}
\end{align}
and imposing
\begin{align}
 \sum_{c \in \mathcal C(n)} \sum_{s \in \mathcal S \mathcal E(n,c)}  |s| \hat f_{s, c} = 0 \qquad \forall n \label{eq:nodalconservationgeneral}
\end{align}
Apart from having doubled the fluxes, this concept of \emph{nodal conservation} merely amounts to a regrouping of terms in a sum. Similar concepts already appear in the context of residual distribution methods (\cite{abgrall06}), as well as in a cell-centered Lagrangian Finite Volume method (\cite{despres05,maire07}). This method is also the first occurrence of the so-called nodal solver for approximating the sub-face fluxes.

\begin{figure}
 \centering
 \includegraphics[width=0.9\textwidth]{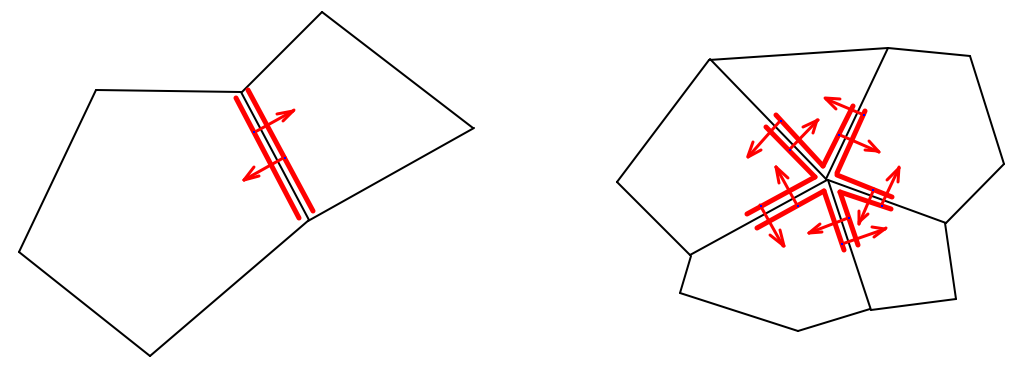}
 \caption{Different conservation concepts. \emph{Left}: Local conservation focusing on an edge. \emph{Right}: Local conservation focusing on a node.}
 \label{fig:conservationconcepts}
\end{figure}

Note that nodal conservation is a truly multi-dimensional concept, not existing in one spatial dimension.
Edge conservation implies nodal conservation, because if the two fluxes across any subedge cancel out already (see Figure \ref{fig:conservationconcepts}), the sum over subedges around a node contains only zeros. In the following we will be deriving a method where this does not happen, and when computing fluxes using a Riemann solver we thus are forced to use a Riemann solver conceptually different from the usual ones.  

Observe that by considering fluxes on subedges the number of numerical fluxes has been increased by two. While edge conservation provides $\# \mathcal E$ equations, nodal conservation provides $\#\mathcal N$ equations. One thus is also forced to make a few identifications, i.e. to add a few equations in order to be able to close the system and obtain one value per subedge-flux. Before coming to that, we will next derive an unusual Riemann solver for linear acoustics that will be used to compute fluxes through subedges. The solver will contain a free parameter, that ultimately will be used to enforce nodal conservation.

\section{An approximate Riemann solver for linear acoustics} \label{sec:riemann}

\subsection{General linear systems in multi-d}\label{sec:generallinear}

Consider a linear $m \times m$ system
\begin{align}
 \del_t q + (\vec J \cdot \nabla) q &= 0 & q &\colon \mathbb R^+_0 \times \mathbb R^d \to \mathbb R^m
\end{align}
where $\vec J$ is a set of $d$ square $m \times m$ Jacobian matrices $J_1, J_2, \ldots, J_d$ and
\begin{align}
 \vec J \cdot \nabla = J_1 \del_{x_1} + J_2 \del_{x_2} + \ldots + J_d \del_{x_d}.
\end{align}
The system is called hyperbolic if $\vec J \cdot \vec k$ is real diagonalizable for any $\vec k \in \mathbb R^d \backslash \{0\}$.

\subsection{Linear acoustics}

Linear acoustics is the following $m = d+1$ system
\begin{align}
 \del_t \vec v + \nabla p &= 0 & \vec v &\colon \mathbb R^+_0 \times \mathbb R^d \to \mathbb R^d \label{eq:linacv}\\
 \del_t p + \nabla \cdot \vec v &= 0 & p &\colon\mathbb R^+_0 \times \mathbb R^d \to \mathbb R. \label{eq:linacp}
\end{align}
The eigenvalues of $\vec J \cdot \vec k$ are $0$ and $\pm|\vec k|$. In the following, we will restrict ourselves to the case $d=2$ and will write $\vec v \equiv (u, v)$, which is system \eqref{eq:linacu2d}--\eqref{eq:linacp2d}.

The form \eqref{eq:linacv}--\eqref{eq:linacp} allows to immediately conclude that linear acoustics is covariant under rotations $\mathcal R \in \mathrm{SO}(2)$ of the $x$-$y$-plane around its origin. This is because the divergence and the gradient transform under rotations as, respectively, scalars and vectors. It thus suffices to consider a Riemann problem in $x$-direction to cover all Riemann problems where the initial discontinuity is a straight line. Note that this is due to the special form of linear acoustics and not generally true for all linear systems.

Along the $x$-direction the system then reads
\begin{align}
 \del_t \veccc{u}{v}{p} + \left( \begin{array}{ccc} 0 & 0 & 1 \\ 0 & 0 & 0 \\ 1 & 0 &0 \end{array} \right ) \del_x \veccc{u}{v}{p} &= 0.
\end{align}
The Jacobian has eigenvalues $0, \pm1$. The exact solution of the Riemann problem is
\begin{align}
 q(t, x) &= \begin{cases} q_\text{L} & \frac{x}{t} < -1 \\ q_1^* & -1 \leq \frac{x}{t} < 0 \\ q_2^* & 0 \leq \frac{x}{t} < 1 \\ q_\text{R} & 1 \leq \frac{x}{t} \end{cases}
\end{align}
with
\begin{align}
 q_1^* &= \veccc{u^*}{v_\text{L}}{p^*} &q_2^* &= \veccc{u^*}{v_\text{R}}{p^*}
\end{align}
\begin{align}
 u^* &= \frac{u_\text{R} + u_\text{L}}{2} - \frac{p_\text{R} - p_\text{L}}{2} &
 p^* &= \frac{p_\text{R} + p_\text{L}}{2} - \frac{u_\text{R} - u_\text{L}}{2} .\label{eq:riemannproblemacousticsexact}
\end{align}

A Riemann problem between states $q_\text{L}$ and $q_\text{R}$ separated by a discontinuity with unit normal\footnote{The notation L/R being chosen such that the normal points from L to R.} $\vec n = (n_x, n_y)^\text{T}$ can be rotated by
\begin{align}
 \mathcal R = \left( \begin{array}{cc} n_x & n_y \\ - n_y & n_x \end{array} \right)
\end{align}
to become a Riemann problem in $x$-direction as considered above. The normal vector is thus obviously mapped to $(1, 0)^\text{T}$, while the velocity $(u,v)^\text{T}$ is mapped to 
\begin{align}
\vecc{n_x u  + n_y v }{ -n_y u + n_x v } =: \vecc{\vec v \cdot \vec n}{ \vec v \cdot \vec n^\perp}.
\end{align}
The exact, or approximate solution $q(t, x) =: (u(t, x), v(t, x), p(t, x))$ of the Riemann problem is then to be rotated back by applying 
\begin{align}
\mathcal R^{-1} = \left( \begin{array}{cc} n_x & -n_y \\ n_y & n_x \end{array} \right)
\end{align}
to $(u(t,x), v(t,x))^\text{T}$ and keeping $p(t,x)$ as it is. This is also true for normal fluxes $\vec f \cdot \vec n$, because they transform just as the respective variables.

\subsection{Modified solution of the Riemann problem}

Following \cite{harten83}, consider an approximate solution $\mathscr R(\frac{x}{t}, q_\text{L}, q_\text{R})$ to the Riemann problem with the quickest/slowest waves moving at speeds $\pm\lambda$, $\lambda > 0$. Integration of the conservation law over the space-time volume $[0, \Delta t] \times [0, \lambda \Delta t]$ gives
\begin{align}
 \Delta t (f(q_\text{R}) - f^*_\text{R}) + \int_0^{\lambda \Delta t} \mathscr R\left(\frac{x}{\Delta t}, q_\text{L}, q_\text{R}\right) \, \dd x - q_\text{R} \lambda \Delta t  = 0\label{eq:approxriemannsolverR}
\end{align}
while integration over $[0, \Delta t] \times [-\lambda \Delta t, 0]$ gives
\begin{align}
 \Delta t (f^*_\text{L} - f(q_\text{L})) + \int_{-\lambda \Delta t}^{0} \mathscr R\left(\frac{x}{\Delta t}, q_\text{L}, q_\text{R}\right) \, \dd x - q_\text{L} \lambda \Delta t = 0.\label{eq:approxriemannsolverL}
\end{align}
Here, $f^*_\text{L/R}$ are the fluxes along the edge, and \cite{harten83} remark that enforcing $f^*_\text{R} = f^*_\text{L}$ (edge-conservation) results in the well-known consistency condition for approximate Riemann solvers. We will eventually not do this here.

If one chooses the approximate Riemann solution to consist of two constant states $q^*_\text{L/R}$, then \eqref{eq:approxriemannsolverR}--\eqref{eq:approxriemannsolverL} would read
\begin{align}
 f(q_\text{R}) - f^*_\text{R} &= \lambda(q_\text{R} -  q^*_\text{R}) \\
 f^*_\text{L} - f(q_\text{L}) &= -\lambda (q^*_\text{L} - q_\text{L} )
\end{align}
which are the Rankine-Hugoniot conditions across the left and right waves.

\subsubsection{Classical solver}

For linear acoustics, denote the fluxes of $u$ by $\bar p_\text{L/R}$ and those of $p$ by $\bar u_\text{L/R}$. The Rankine-Hugoniot conditions over all three waves, whose speeds we choose as $-1, 0, 1$ read
\begin{align}
 -(u^*_\text{L} - u_\text{L}) &= \bar p_\text{L} - p_\text{L},  & \bar p_\text{L} &= \bar p_\text{R} ,&  u_\text{R} - u^*_\text{R} &= p_\text{R} - \bar p_\text{R}, \label{eq:approxriemanngeneral1}\\
 -(v^*_\text{L} - v_\text{L}) &= 0,  &  0&=0, &  v_\text{R} - v^*_\text{R} &= 0,\\
 -(p^*_\text{L} - p_\text{L}) &= \bar u_\text{L} - u_\text{L},  & \bar u_\text{L} &= \bar u_\text{R}, &  p_\text{R} - p^*_\text{R} &= u_\text{R} - \bar u_\text{R}. \label{eq:approxriemanngeneral3}
\end{align}
Recall that there is no $x$-flux of $v$. The middle row implies
\begin{align}
 v^*_\text{L} &= v_\text{L} &  v_\text{R} &= v^*_\text{R}
\end{align}
The remaining two rows of system \eqref{eq:approxriemanngeneral1}--\eqref{eq:approxriemanngeneral3} form 6 equations for 10 variables: $u^*_\text{L/R}, v^*_\text{L/R}, p^*_\text{L/R}$ as well as $\bar u_\text{L/R}, \bar p_\text{L/R}$, an underdetermined system.

As the wave speeds are exact, one obtains the exact solution \eqref{eq:riemannproblemacousticsexact} of the Riemann problem upon \emph{defining} the fluxes as $f^*_{\text{L/R}} := J_x q^*_\text{L/R}$, i.e.
\begin{align}
 \bar u_\text{L} &= u^*_\text{L}  & \bar p_\text{L} &= p^*_\text{L}\\ 
 \bar u_\text{R} &= u^*_\text{R}  & \bar p_\text{R} &= p^*_\text{R}.
\end{align}

Alternatively, one can also \emph{assume} $u^*_\text{L} = u^*_\text{R} =: u^*$ and $p^*_\text{L} = p^*_\text{R} =: p^*$, which leads to the same result.

\subsubsection{Solver with a free velocity variable} \label{ssec:solvervelocity}

Ignore now the equation $\bar p_\text{L} = \bar p_\text{R}$, but \emph{assume} again $u^*_\text{L} = u^*_\text{R}$. Upon defining the shorthand notations $u^* := u^*_\text{L} = u^*_\text{R}$ and $\bar u := \bar u_\text{L} = \bar u_\text{R}$ (a consequence of the middle column of \eqref{eq:approxriemanngeneral1}--\eqref{eq:approxriemanngeneral3}), the equations read
\begin{align}
 -(u^* - u_\text{L}) &= \bar p_\text{L} - p_\text{L}   &  u_\text{R} - u^* &= p_\text{R} - \bar p_\text{R} \\
 -(p^*_\text{L} - p_\text{L}) &= \bar u - u_\text{L}  &   p_\text{R} - p^*_\text{R} &= u_\text{R} - \bar u .
\end{align}
Having removed one equation we expect to be left with one free variable, and we will keep $u^*$ for this purpose. Still, the above system of four equations contains 5 remaining variables ($\bar p_\text{L/R}, p^*_\text{L/R}, \bar u$). 

By assuming either $\bar u = u^*$, or $p_\text{L}^* = \bar p_\text{L}$, or $p_\text{R}^* = \bar p_\text{R}$, one is uniquely led to the solution
\begin{align}
 p_\text{L}^* = \bar p_\text{L} &= p_\text{L} - u^* + u_\text{L}  &
 p_\text{R}^* = \bar p_\text{R} &= p_\text{R} + u^* -u_\text{R} .\label{eq:ustarstatesp}
\end{align}
This is the solver introduced in \cite{gallice22}. Recall that $u^*$ is at this point a free variable to be determined later.

\subsubsection{Solver with a free pressure variable} \label{ssec:solverpressure}

Starting out from \eqref{eq:approxriemanngeneral1}--\eqref{eq:approxriemanngeneral3} and now ignoring the equation $\bar u_\text{L} = \bar u_\text{R}$, keeping the equation $\bar p_\text{L} = \bar p_\text{R} =: \bar p$, and assuming $p^*_\text{L} = p^*_\text{R} =: p^*$ as well as $p^* = \bar p$, yields a mirror image of the previous result:
\begin{align}
 u^*_\text{L} = \bar u_\text{L} &= u_\text{L}-p^* + p_\text{L} &  \bar u_\text{R} = u^*_\text{R} &= u_\text{R} + p^* - p_\text{R}. \label{eq:pstarstatesu}
\end{align}
Here, $p^*$ is a free unknown to be determined elsewhere, consistent with the framework in \cite{gallice22}.

Despite the symmetry between Equations \eqref{eq:ustarstatesp} and \eqref{eq:pstarstatesu}, the solvers that emerge from these two modifications of the Riemann problem have fundamentally different properties, as will be seen later. This is due to the fact that the values of the free parameters will be determined from nodal conservation, which will naturally see the multi-dimensionality of the problem (Section \ref{ssec:nodalconservation}). While $u$ and $p$ indeed can be simply exchanged in the one-dimensional Riemann solvers discussed above, the multi-dimensional conservation requirement will be aware of the fact that it is not possible to exchange $\vec v$ and $p$ in multiple dimensions.

\section{Combining the new Riemann solvers with nodal conservation} \label{sec:method}

\subsection{General idea}

Consider the solvers from sections \ref{ssec:solvervelocity} and \ref{ssec:solverpressure}. In order to ``embed'' those one-dimensional solvers into the multi-dimensional setting, the fluxes $f^*$ are understood as those across an edge $e$ of a cell $c$ with the outward normal vector $(1, 0)^\text{T}$. This edge is, of course, shared by two cells, and for the other cell ($c'$, say) the outward unit vector associated to $e$ is $(-1,0)^\text{T}$. Cell $c'$ therefore is updated with $-f^*$. For general orientations of the edges we make use of the rotational covariance of the acoustic system.

\begin{definition}
 Denote by $\vec n_{e,c}$ the \emph{outward unit normal} at edge $e \in \mathcal E$, pointing out of cell $c \in \mathcal C$ (see Figure \ref{fig:celldual}) and by $\vec n_e$ a normal associated to edge $e \in \mathcal E$ that is fixed once for all. Let the subedges inherit the normals and the direction from the edges they are part of, i.e. $\vec n_s = \vec n_e$, as well as $\vec n_{e, c} = \vec n_{s, c}$ if $s$ is a subedge of $e$.
\end{definition}

The one-dimensional solvers can thus be used with $u_\text{L/R} := \vec v_\text{L/R} \cdot \vec n_e$ where $\text{L/R}$ is defined such that $\vec n_e$ is pointing from cell L to cell R. As will be seen below, the particular choice of orientation of $\vec n_e$ does not matter, as long as some choice is made once for all.

The solvers from section \ref{ssec:solvervelocity} and \ref{ssec:solverpressure} have both been designed in such a way that they leave a free variable, $u^*$ or $p^*$, that will now be determined by imposing conservation around a node (Equation \eqref{eq:nodalconservationgeneral}). There is still a mismatch in the number of free variables and the number of equations, though, that will be dealt with by appropriately identifying the variables associated to half-edges around each node. This will be detailed next in the derivation of the methods.

\subsection{A method with velocities at nodes} \label{ssec:methodgallice}

For the solver from section \ref{ssec:solvervelocity}, the free variable is $u^*$, which in the multi-dimensional context has the meaning of the normal component $\vec v \cdot \vec n_e$. This is one scalar variable per subedge impinging on node $n$. The fluxes of the pressure cancel out on every subedge (edge-conservation), such that for them nodal conservation is trivially fulfilled and only the conservation of $\vec v$ remains to be ensured. These are two equations per node $n$. It is thus natural to identify all the free variables $u^*$ on the subedges $s \in \mathcal{SE}(n)$ impinging on $n$ as being the projections onto the respective $\vec n_s$ of some velocity $\vec v^*_n$ associated to the node $n$:
\begin{align}
 u^*\Big|_{s \in \mathcal S \mathcal E(n)} := \vec v^*_n \cdot \vec n_s \qquad \forall n \in \mathcal N.
\end{align}

At every node $n$, this results in two equations for the two components $\vec v^*_n$. This solver has been derived in \cite{gallice22}. As can be shown by considering it on Cartesian grids, it is not vorticity preserving (see Section \ref{app:nonvortpres}).

\subsection{The nodal pressure}

In the case of the free variable being $p^*$ (solver from section \ref{ssec:solverpressure}), there is one free variable per subedge adjacent to node $n$, and there is only one non-trivial equation: that of conservation of $p$, because the fluxes of $\vec v$ cancel out on every subedge (edge-conservation). It is natural to identify all the $p^*$ on the subedges impinging at node $n$ with just a single pressure $p^*_n$ associated to the node:
\begin{align}
 p^*\Big|_{s \in \mathcal S \mathcal E(n)} := p^*_n  \qquad \forall n \in \mathcal N.
\end{align}

At every node, one is thus left with one equation for one variable. To base the method on a nodal pressure is the new contribution of this present work and will be investigated in detail next.

\begin{figure}
 \centering
 \includegraphics[width=0.3\textwidth]{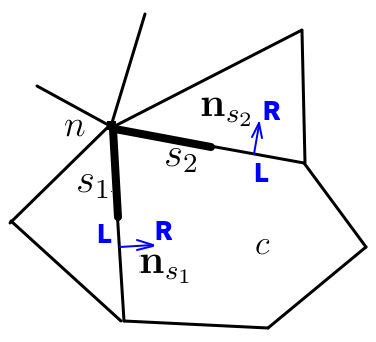} \hspace{2cm} \includegraphics[width=0.45\textwidth]{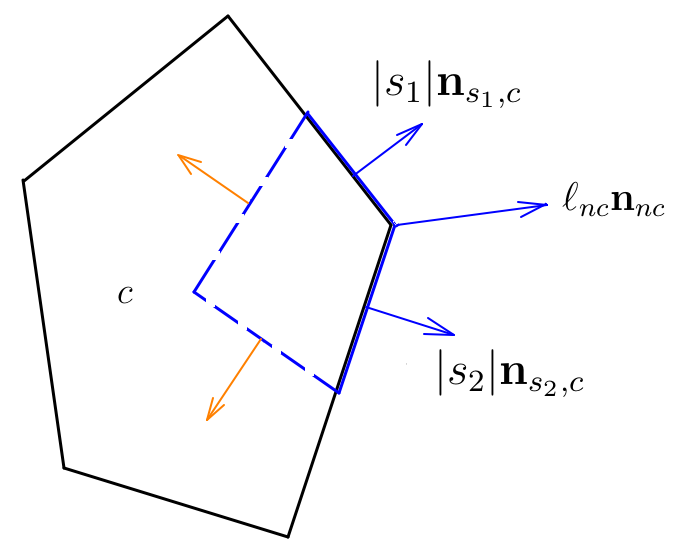}
 \caption{\emph{Left:} The definition of L/R depends on the choice of the normal vector $\vec n_e$ on edge $e$. If outward normals with respect to cell $c$ are used, then cell $c$ is always L. \emph{Right}: The node normal.}
 \label{fig:normals}
\end{figure}

Node-conservation \eqref{eq:nodalconservationgeneral}, i.e.
\begin{align}
 0 = \sum_{c \in \mathcal C(n)} \sum_{s \in \mathcal S \mathcal E(c, n)} |s| \hat f_{s,c} \label{eq:nodalconservationrepeat}
\end{align}
requires summing the fluxes over all subedges $s \in \mathcal S \mathcal E(n)$ adjacent to the node $n$. Recall that we associate a normal vector $\vec n_s$ to any subedge $s \in \mathcal S \mathcal E(n)$. If the whole setup is rotated such that this normal is $(1, 0)^\text{T}$, then this normal is the outward unit normal for cell L, and minus the outward unit normal for cell R. As $\hat f_{s, c}$ is defined with respect to the outward unit normal, we need to use the flux $\bar u_\text{L}$ as updating cell L and flux $-\bar u_\text{R}$ as updating cell R. Thus, \eqref{eq:nodalconservationrepeat} becomes

\begin{align}
 0 &=  \sum_{s \in \mathcal{SE}(n)} |s| (-\bar u_\text{R} + \bar u_\text{L})\\
  &\overset{\eqref{eq:pstarstatesu}}{=} -\sum_{s \in \mathcal{SE}(n)} |s| (\vec v_\text{R} \cdot \vec n_{s} - p_\text{R} +2 p^*_{n} - \vec v_\text{L} \cdot \vec n_{s} - p_\text{L}).
\end{align}
This gives
\begin{align}
  p^*_{n}  &= \frac{\displaystyle\sum_{s \in \mathcal{SE}(n)} |s| \left( \frac{p_\text{R} + p_\text{L}}{2} - \frac{(\vec v_\text{R}- \vec v_\text{L}) \cdot \vec n_{s}}{2} \right)}{\displaystyle\sum_{s \in \mathcal{SE}(n)} |s|}. \label{eq:nodalpressureunstructuredraw}
\end{align}
The choice of normal direction $\vec n_s$ does not influence the result, but only the definition of L/R, because the expression
\begin{align}
\sum_{s \in \mathcal{SE}(n)} |s| \left(\vec v_\text{R} - \vec v_\text{L}\right ) \cdot \vec n_{s} \label{eq:unstrdivlr}
\end{align}
can be rewritten without reference to L/R as follows (see also Figure \ref{fig:normals}, left). For any cell $c \in \mathcal C(n)$ there are two subedges $\mathcal{SE}(n,c) = \mathcal{SE}(n) \cap \mathcal{SE}(c)$ involved in the summation \eqref{eq:unstrdivlr}; let us call them $s_1$ and $s_2$. If $\vec n_{s_1}$ is pointing out of cell $c$, i.e. if $\vec n_{s_1} = \vec n_{s_1, c}$, then $c$ is the ``left'' cell and the contribution of cell $c$ to the sum is
\begin{align}
- |s_1| \vec v_c \cdot \vec n_{s_1,c}.
\end{align}
Otherwise, if $\vec n_{s_1} = -\vec n_{s_1, c}$, then $c$ is the ``right'' cell, and the contribution is the same. Therefore the total contribution of cell $c$ is
\begin{align}
-  (|s_1| \vec n_{s_1,c} + |s_2| \vec n_{s_2,c}) \cdot \vec v_c
\end{align}
which allows to rewrite Equation \eqref{eq:unstrdivlr} into
\begin{align}
\sum_{s \in \mathcal{SE}(n)} |s| \left(\vec v_\text{R} - \vec v_\text{L}\right ) \cdot \vec n_{s} = - \sum_{c \in \mathcal C(n)} \left( \sum_{s \in \mathcal{SE}(n, c)} |s| \vec n_{s,c} \right ) \cdot \vec v_c.
\end{align}

This motivates the following definition (see also \cite{loubere16}):

\begin{definition} \label{def:nodenormal}
 Denote by $\ell_{nc} \vec n_{nc} := \displaystyle \sum_{s \in \mathcal{SE}(n,c) } |s| \vec n_{s,c}$ the \emph{node-normal} ($|\vec n_{nc}| = 1$).
\end{definition}

Note that the node-normal is equal to the negative of the sum of weighted normals to the two lines connecting the centroid of the cell with the edge midpoints, see Figure \ref{fig:normals}, right.

\begin{theorem} \label{thm:normal} The following identities hold on any cell $c \in \mathcal C$:
 \begin{enumerate}[(i)]
  \item \label{it:thmnormal1} $\sum_{n \in \mathcal N(c)} \ell_{nc} \vec n_{nc} = 0$,
  \item \label{it:thmnormal2}$\sum_{n \in \mathcal N(c)} \ell_{nc} \vec n_{nc} \otimes \vec x_n = |c| \id$.
 \end{enumerate}

\end{theorem}
\begin{proof}
 \begin{enumerate}[(i)]
  \item The node-normal $\ell_{nc} \vec n_{nc}$ sums the normals of the two subedges adjacent to node $n$, weighted by their respective lengths. Summed over all nodes this simply gives a sum over all the edge normals of the polygon, weighted by the lengths of the edges:
  \begin{align}
  \sum_{n \in \mathcal N(c)} \ell_{nc} \vec n_{nc} = \sum_{e \in \mathcal E(c)} \vec n_{e, c} |e| = 0.
  \end{align}
  \item In the sum, consider the contribution of the two nodes $n_1, n_2$ joined by an edge $e$, and denote by $s_1 \in \mathcal{SE}(n)$ the subedge of $e$ adjacent to $n_1$, and analogously $s_2$. Their contributions concerning $e$ are $|s_1| \vec n_{s_1,c} \otimes \vec x_{n_1}$ and $|s_2| \vec n_{s_2,c} \otimes \vec x_{n_2}$. Using $ \vec n_{s_1,c}  =  \vec n_{s_2,c}  = \vec n_{e,c}$ and $|s_1| = |s_2| = \frac12 |e|$ yields
  \begin{align}
   \sum_{n \in \mathcal N(c)} \ell_{nc} \vec n_{nc} \otimes \vec x_n &= \sum_{e \in \mathcal E(c)} |e| \vec n_{e, c} \otimes \frac12 \sum_{n \in \mathcal N(e)} \vec x_n = \sum_{e \in \mathcal E(c)} |e|\vec n_{e, c} \otimes \int_0^1 \vec x \, \dd s.
 \end{align}
  The average coordinates of the two nodes are the coordinates of the midpoint of the edge (first equality), and also the average of the coordinates of points along the edge (second equality). Defining $\vec x_{n_1} + s (\vec x_{n_2} - \vec x_{n_1})$ as parametrization of these points one finally obtains  (see also \cite{loubere16})
    \begin{align}
   \sum_{n \in \mathcal N(c)} \ell_{nc} \vec n_{nc} \otimes \vec x_n   &= \oint_{\del c} \vec n_c \otimes \vec x \, \dd \ell  \overset{\text{Gauss}}{=} \int_c (\nabla \vec x) \dd \vec x = |c| \id.
  \end{align}
  
 \end{enumerate}
 
%\hfill $\Box$ 
\end{proof}

\begin{theorem}
 Consider a cell $c \in \mathcal C$.
 \begin{enumerate}[(i)]
  \item The discrete cell-centered operator 
  \begin{align}(\vec G \phi)_c := \frac{1}{|c|} \sum_{n \in \mathcal N(c)} \ell_{nc} \vec n_{nc} \phi_n , \quad \vec G \colon \mathcal N \to \mathcal C^d \label{eq:discretegradunstructured} \end{align}
  is a consistent discretization of the gradient $\nabla \phi$ at the centroid of $c$, in the sense that it is exact on affine functions $\phi$. 
  \item The discrete node-centered operator 
  \begin{align}
  (D \vec v)_n := -\frac{1}{|c_n|} \sum_{c \in \mathcal C(n)} \ell_{nc} \vec n_{nc} \cdot \vec v_c, \quad D \colon \mathcal C^d \to \mathcal N \label{eq:discretedivunstructured}
  \end{align}
  is a weakly consistent discretization of the divergence $\nabla \cdot \vec v$ at the location of the node $n$, in the sense that it is the $\ell^2$ dual of $\vec G$.
 \end{enumerate}
\end{theorem}
\emph{Note}: As much as we denote by $\phi$ a function and by $\phi_n$ its value at a node, we denote by $\vec G \phi$ the function obtained by applying the discrete gradient $\vec G$ to $\phi$ everywhere, and by $(\vec G \phi)_c$ the value of this function at cell $c$. The operator itself will be denoted by $\vec G$, or $\vec G_c$ if the discrete gradient operator at a specific cell $c$ is meant. This means that $\vec G_c \phi \equiv (\vec G \phi)_c$. Analogous notation is used for $D$.

\begin{proof}
  \begin{enumerate}[(i)]
  \item Consider an affine field $\phi$, such that $\phi_n = a + \vec b \cdot \vec x_n$. Then
  \begin{align}
   \vec G_c \phi = \frac{1}{|c|} a \underbrace{\sum_{n \in \mathcal N(c)} \ell_{nc} \vec n_{nc}}_{= 0} + \frac{1}{|c|} \underbrace{ \sum_{n \in \mathcal N(c)} \ell_{nc} \vec n_{nc} \otimes \vec x_n }_{= |c| \id}\vec b \overset{Thm. \ref{thm:normal}}{=} \vec b = \nabla \phi
  \end{align}
  \item The aim is to mimic $\int (\nabla \cdot \vec v) \phi \,\dd \vec x= - \int \vec v \cdot \nabla \phi \, \dd \vec x + \mathrm{b.t.}$ Consider the $\ell^2$ scalar product $\langle a, b \rangle := \sum_{n \in \mathcal N} a_n b_n |c_n| $ with a scalar function $\phi$, defined on nodes:
  \begin{align}
   \sum_{n \in \mathcal N} D_n \vec v \phi_n |c_n| &= -\sum_{n \in \mathcal N} \sum_{c \in \mathcal C(n)} \ell_{nc} \vec n_{nc} \cdot \vec v_c \phi_n \\
   &= -\sum_{c \in \mathcal C} \vec v_c \cdot \sum_{n \in \mathcal N(c)}  \ell_{nc} \vec n_{nc} \phi_n = -\sum_{c \in \mathcal C} \vec v_c \cdot \vec G_c \phi |c|
  \end{align}
  \end{enumerate}

 % \hfill $\Box$

\end{proof}

Note that the operator $D_n$ is generally not exact on affine functions. Moreover, $\vec G$ and $D$ are not the operators $\widetilde{\mathcal{GRAD}}$ and $ \mathcal{DIV}$ from \cite{lipnikov14}, Section 6, but would be their duals.

Having established these operators, one can rewrite Equation \eqref{eq:nodalpressureunstructuredraw} as
\begin{align}
  p^*_{n}  &= \frac{\displaystyle\sum_{s \in \mathcal{SE}(n)} |s|  \frac{p_\text{R} + p_\text{L}}{2}  -  \frac12 |c_n| D_n \vec v}{\displaystyle\sum_{s \in \mathcal{SE}(n)} |s|} \label{eq:nodalpressureunstructuredrawdiv}
\end{align}
This expression is worked out on Cartesian grids in Equation \eqref{eq:pstardefinition}.

\subsection{The new method with pressure at nodes} \label{ssec:nodalconservation}

Finally the update equations will be derived, starting out from Equation \eqref{eq:finvolsubedges} and using the nodal pressure \eqref{eq:nodalpressureunstructuredrawdiv}. Note that the normal vectors $\vec n_{s,c}$ are chosen to point outwards of cell $c$, such that cell $c$ is always the left (L) cell and every subedge contributes with its own fluxes: $\bar u_\text{L}$ as flux of $p$, and $\vec n_{s, c} p^*_n$ as flux of $\vec v$. The semi-discrete update of cell $c$ thus reads as follows:
\begin{align}
\frac{\dd}{\dd t} \vec v_c &= - \frac{1}{|c|} \sum_{n \in \mathcal N(c)} \sum_{s \in \mathcal{SE}(n, c)} |s| \vec n_{s, c} p^*_n \overset{\text{Def. \ref{def:nodenormal}}}{=} - \frac{1}{|c|} \sum_{n \in \mathcal N(c)} \ell_{nc} \vec n_{nc} p^*_n \overset{\eqref{eq:discretegradunstructured}}{=} - \vec G_c p^*, \label{eq:pstarsolverevoutionofugeneral}\\
 \frac{\dd}{\dd t} p_c &=  - \frac{1}{|c|} \sum_{n \in \mathcal N(c)} \sum_{s \in \mathcal{SE}(n, c)} |s| \bar u_\text{L} \\
 &= - \frac{1}{|c|} \sum_{n \in \mathcal N(c)} \sum_{s \in \mathcal{SE}(n, c)} |s| (\vec v_c \cdot \vec n_{s,c} - p^*_n + p_c) \\
 &= - \frac{1}{|c|} \sum_{n \in \mathcal N(c)} \sum_{s \in \mathcal{SE}(n, c)} |s| (- p^*_n + p_c).  \label{eq:pstarsolverevoutionofpgeneral}
\end{align}
 Theorem \ref{thm:normal}(\ref{it:thmnormal1}) has been used to obtain \eqref{eq:pstarsolverevoutionofpgeneral}. Observe the appearance of the natural discrete gradient $\vec G$ in \eqref{eq:pstarsolverevoutionofugeneral}. While in \cite{morton01}, the authors merely suggest to use this discretization of the gradient, here it is a strict consequence of the chosen Riemann solver and of nodal conservation.

 Using Equation \eqref{eq:nodalpressureunstructuredrawdiv}, the pressure equation \eqref{eq:pstarsolverevoutionofpgeneral} can be rewritten as
 \begin{align}
  \frac{\dd}{\dd t} p_c &=- \frac{1}{|c|} \sum_{n \in \mathcal N(c)} \left(\sum_{s \in \mathcal{SE}(n, c)} |s|\right ) \left(p_c- \frac{\displaystyle\sum_{s \in \mathcal{SE}(n)} |s|  \frac{p_\text{R} + p_\text{L}}{2}  -  \frac12 |c_n| D_n \vec v}{\displaystyle\sum_{s \in \mathcal{SE}(n)} |s|} \right) 
 \end{align}
 which makes explicit the cell-centered approximation of the divergence 
 \begin{align}
  \tilde D_c \vec v := \frac{1}{|c|} \sum_{n \in \mathcal N(c)} \frac{ \displaystyle \sum_{s \in \mathcal{SE}(n, c)} |s|}{\displaystyle\sum_{s \in \mathcal{SE}(n)} |s|}  \frac12 |c_n| D_n \vec v
 \end{align}
 which is a linear combination of the node-centered divergence $D_n \vec v$. It is \emph{not} an average, i.e. the coefficients
 \begin{align}
  \alpha_{nc} := \frac{1}{|c|} \frac{ \displaystyle \sum_{s \in \mathcal{SE}(n, c)} |s|}{\displaystyle\sum_{s \in \mathcal{SE}(n)} |s|}  \frac12 |c_n|
 \end{align}
 do not generally sum up to 1. This is clear because $|c_n|$ can be varied while keeping all the lengths involved in the summation over $\mathcal{SE}(n)$ constant. But then again, $D_n \vec v$ is not a strong divergence either; and one does indeed confirm a weak form of the averaging property:
 
 \begin{theorem}
 For any scalar function $\phi$ defined on the nodes
 \begin{align}
  \sum_{c \in \mathcal C} |c| \sum_{n \in \mathcal N(c)} \alpha_{nc} \phi_n &= \sum_{n \in \mathcal N} |c_n| \phi_n
  \end{align}
  \end{theorem}
  \begin{proof}
  \begin{align}  
  \sum_{c \in \mathcal C} |c| \sum_{n \in \mathcal N(c)} \alpha_{nc} \phi_n  &= \sum_{c \in \mathcal C} \sum_{n \in \mathcal N(c)} \frac{ \displaystyle \sum_{s \in \mathcal{SE}(n, c)} |s|}{\displaystyle\sum_{s \in \mathcal{SE}(n)} |s|}  \frac12 |c_n| \phi_n 
  \end{align}
  A given node $n$ appears several times in the summation $\sum_{c \in \mathcal C} \sum_{n \in \mathcal N(c)}$ that one would like to replace by a summation $\sum_{n \in \mathcal N} \sum_{c \in \mathcal C(n)}$. Most terms depend just on $n$, and only $\sum_{s \in \mathcal{SE}(n, c)} |s|$ requires attention. Here, focusing on some particular node $n$, every subedge appears twice because it belongs to two cells. Upon reversing the summations one therefore picks up a factor of 2: $$\sum_{c \in \mathcal C(n)} \sum_{s \in \mathcal{SE}(n, c)} |s| = 2 \sum_{s \in \mathcal{SE}(n)} |s|. $$
  Thus,
  \begin{align}
 \sum_{c \in \mathcal C} |c| \sum_{n \in \mathcal N(c)} \alpha_{nc} \phi_n  &= \sum_{n \in \mathcal N}  \frac{ 1}{\displaystyle\sum_{s \in \mathcal{SE}(n)} |s|}  \frac12 |c_n| \phi_n \cdot  \sum_{c \in \mathcal C(n)}   \sum_{s \in \mathcal{SE}(n,c)} |s| \\
  &= \sum_{n \in \mathcal N}    |c_n| \phi_n   
 \end{align}
 
 \end{proof}

 Cartesian versions of Equations \eqref{eq:pstarsolverevoutionofugeneral}--\eqref{eq:pstarsolverevoutionofpgeneral} are given as Equations \eqref{eq:pstarsolverevoutionofp}, \eqref{eq:pstarsolverevoutionofu}. In the following, it will be shown that the semi-discretization \eqref{eq:pstarsolverevoutionofugeneral}--\eqref{eq:pstarsolverevoutionofpgeneral} is vorticity preserving.

\section{Vorticity and stationarity preservation} \label{sec:statioanritypreservation}

In \eqref{eq:pstarsolverevoutionofugeneral}--\eqref{eq:pstarsolverevoutionofpgeneral}, the discrete node-centered divergence $D_n$ contained in $p_n^*$ appears both in the ``physical'' divergence of the evolution equation of $p$, and in the diffusion in the evolution of $\vec v$. This allows to prove the following:

\begin{theorem}
If the initial data are such that $p$ is uniformly constant and $D_n \vec v = 0$ for all nodes $n$, then they remain stationary.
\end{theorem}
\begin{proof}
It follows from the assumptions that $p^* = p$ and then $\frac{\dd}{\dd t} p_c = 0$ as well as 
\begin{align}
 \frac{\dd}{\dd t} \vec v_c &= - \frac{1}{|c|} \sum_{n \in \mathcal N(c)} \ell_{nc} \vec n_{nc} p= 0 
\end{align}

%\hfill $\Box$
\end{proof}

Whether this implies stationarity preservation depends on how many solutions the system of equations
\begin{align}
  D_n \vec v = 0 \qquad \forall n \label{eq:discretedivzero}
\end{align}
has.

\begin{definition}
 A linear numerical method of the form \eqref{eq:pstarsolverevoutionofugeneral}--\eqref{eq:pstarsolverevoutionofpgeneral} is called \emph{stationarity preserving} if it has as many linearly independent stationary states as there are nodes on the grid (up to modifications related to boundaries).
\end{definition}

This definition aims at capturing discretely the richness of non-trivial stationary states of linear acoustics. Next, it will be investigated whether this is the case for the presented method.

Equations \eqref{eq:discretedivzero} are as many equations as there are nodes in the grid ($\#\mathcal N$), while the number of variables is twice the number of cells $\#\mathcal C$, as $\vec v$ has two components. Thus, a priori the system has $2\#\mathcal C - \#\mathcal N$ solutions, and more if it does not have full rank.

For methods on Cartesian grids, involution preservation is equivalent to stationarity preservation (\cite{barsukow17a}). On unstructured grids these two properties remain tightly linked. In order to show involution preservation of method \eqref{eq:pstarsolverevoutionofugeneral}--\eqref{eq:pstarsolverevoutionofpgeneral} one is led to investigate those discrete operators that annihilate the discrete gradient $\vec G_c$ appearing in \eqref{eq:pstarsolverevoutionofugeneral}. Written as a matrix, it has $2\#\mathcal C$ rows and $\#\mathcal N$ columns. Its left kernel (the kernel of its transpose) has at least dimension $2\#\mathcal C  - \#\mathcal N$, because its elements are constrained by $\#\mathcal N$ equations. These are as many as there are stationary states solving \eqref{eq:discretedivzero}. Next, this number will be estimated on different grids.
%a linear system, these would be $\#\mathcal N$ equations for $2 \#\mathcal C$ variables, with at least $2 \#\mathcal C - \#\mathcal N$ solutions. 

For doubly periodic boundaries, the domain is a torus with $\#\mathcal C - \#\mathcal E + \#\mathcal N = 0$. As will be seen below, the operators annihilating the cell-centered gradient will indeed be node-centered discrete curl operators, which we are able to characterize for the case of triangular-quadrangular grids, i.e. grids with $\# \mathcal E(c) \leq 4$ for all cells $c$. This means that $\#\mathcal E \leq 2\#\mathcal C$ and thus
\begin{align}
2 \#\mathcal C - \#\mathcal N &\geq \#\mathcal N,
\end{align}
i.e. we are sure to find as many solutions as there are nodes. In general, if $\# \mathcal E(c) \leq \alpha$, then 
\begin{align}
2 \#\mathcal C - \#\mathcal N &\geq \frac{6-\alpha}{\alpha-2}  \#\mathcal N.
\end{align}
For hexagonal grids one thus faces the interesting issue that $\vec G_c$ amounts to a square matrix, which might have only a trivial kernel if it is full rank.

\begin{corollary}
 The method from Section \ref{ssec:solverpressure} with a nodal pressure is stationarity preserving on triangular-quadrangular grids ($\# \mathcal E(c) \leq 4$).
\end{corollary}

\begin{figure}
  \centering
  \includegraphics[width=\textwidth]{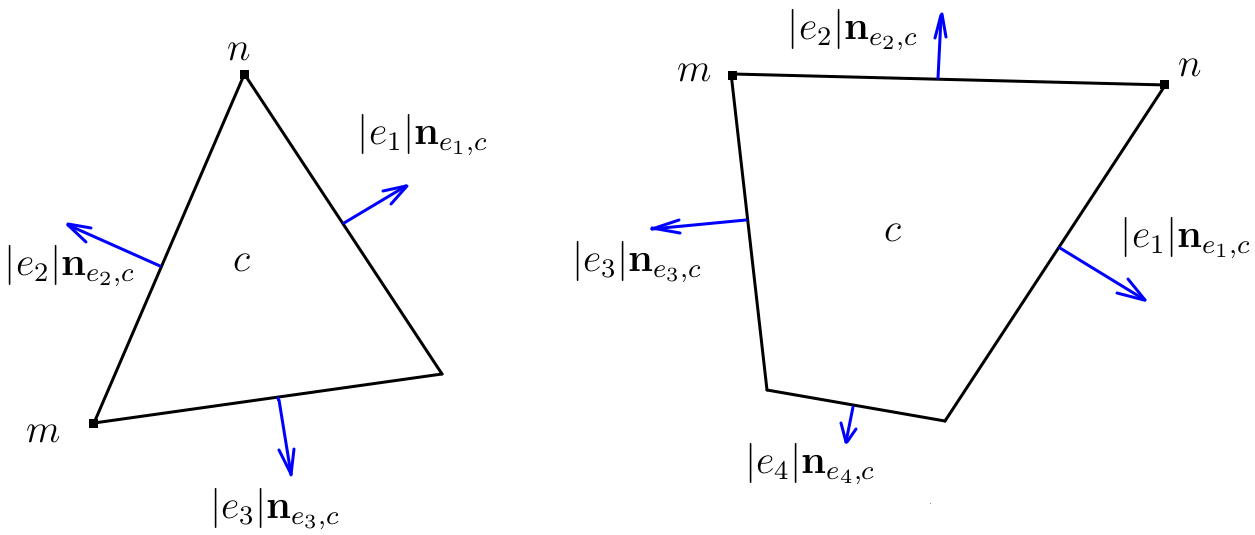}  \includegraphics[width=0.3\textwidth]{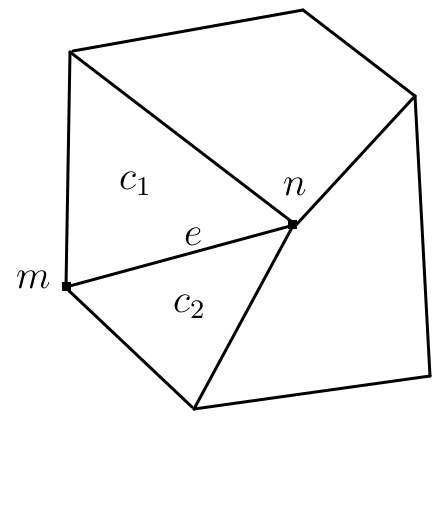}
  \caption{\emph{Top}: Illustration and motivation for the proof of Lemma \ref{lemma:area}. \emph{Bottom}: Notation for the proof of Theorem \ref{thm:curl} (the case of $n \neq m$, but $n$ and $m$ sharing an edge $e$).}
  \label{fig:area}
 \end{figure}

\begin{lemma} \label{lemma:area}
 For a cell $c$ with $\# \mathcal E(c) \leq 4$ we have
 \begin{align}
  |c| = 2 \ell_{nc} \vec n_{nc} \times \ell_{mc} \vec n_{mc}
 \end{align}
 if nodes $n$ and $m$ share an edge and are ordered counterclockwise along the boundary of $c$.
\end{lemma}
\begin{proof}
 Number the edges counterclockwise as $1, 2, 3, (4)$ (see Figure \ref{fig:area}). Then the area of a 
 triangular cell is given by
 \begin{align}
  |c| &= \frac12 |e_1| \vec n_{e_1,c} \times |e_2| \vec n_{e_2,c} \\
  &= \frac12 |e_1| \vec n_{e_1,c} \times \Big ( |e_1| \vec n_{e_1,c} + |e_2| \vec n_{e_2,c}  \Big ) \\
  &= \frac12 \Big (|e_2| \vec n_{e_2,c} + |e_1| \vec n_{e_1,c} - |e_2| \vec n_{e_2,c}\Big) \times \Big ( |e_1| \vec n_{e_1,c} + |e_2| \vec n_{e_2,c}  \Big ) \\
  &= -\frac12 \Big (|e_3| \vec n_{e_3,c} + |e_2| \vec n_{e_2,c}\Big) \times \Big ( |e_1| \vec n_{e_1,c} + |e_2| \vec n_{e_2,c}  \Big ) \\
  &= 2 \ell_{nc} \vec n_{nc} \times \ell_{mc} \vec n_{mc}
 \end{align}

 by calling $m$ the common node of $e_2$ and $e_3$ and $n$ the common node of $e_1$ and $e_2$.
 Similarly, by triangulation the area of a quadrangular cell is 
 \begin{align}
  |c| &= \frac12 |e_1| \vec n_{e_1,c} \times |e_2| \vec n_{e_2,c}   +    \frac12 |e_3| \vec n_{e_3,c} \times |e_4| \vec n_{e_4,c} \\
  &= \frac12 |e_1| \vec n_{e_1,c} \times |e_2| \vec n_{e_2,c}   -  \frac12 |e_3| \vec n_{e_3,c} \times \Big(|e_1| \vec n_{e_1,c} + |e_2| \vec n_{e_2,c} + |e_3| \vec n_{e_3,c} \Big) \\
  &= \frac12 |e_1| \vec n_{e_1,c} \times |e_2| \vec n_{e_2,c}   -  \frac12 |e_3| \vec n_{e_3,c} \times \Big( |e_1| \vec n_{e_1,c} + |e_2| \vec n_{e_2,c} \Big)\\
  &=   -  \frac12 \Big(|e_2| \vec n_{e_2,c} + |e_3| \vec n_{e_3,c} \Big) \times \Big( |e_1| \vec n_{e_1,c} + |e_2| \vec n_{e_2,c} \Big) \\
  &=   2 \ell_{nc} \vec n_{nc}  \times \ell_{mc} \vec n_{mc}.
 \end{align}
 
% \hfill $\Box$

\end{proof}

\begin{theorem} \label{thm:curl}
 Define the discrete node-centered operator 
 \begin{align} 
 \vec C_n \vec v := -\frac{1}{|c_n|} \sum_{c \in \mathcal C(n)} \ell_{nc} \vec n_{nc} \times \vec v_c
 \label{eq:unstructuredcurl}
 \end{align}
 \begin{enumerate}[(i)]
  \item $\vec C_n \vec v$ is a weakly consistent discretization of the curl $\nabla \times \vec v$ at the location of the node $n$.
  \item On grids with $\# \mathcal E(c) \leq 4$ for all cells $c$ we have:
 \begin{align}
  \mathrm{im}\,\vec G \subset \ker \vec C.
 \end{align}
 \end{enumerate}

\end{theorem}
\begin{proof}
 \begin{enumerate}[(i)]
  \item 
    We aim at a discrete version of $\int (\nabla \times \vec v)  \phi \,\dd \vec x = - \int (\nabla \psi) \times \vec v \, \dd \vec x + \text{b.t.}$ To this end, consider a scalar function $\phi$ defined on the nodes and 
    \begin{align}
      \sum_{n \in \mathcal N} \vec C_n \vec v \phi_n |c_n| &= -\sum_{n \in \mathcal N} \sum_{c \in \mathcal C(n)} \ell_{nc} \vec n_{nc} \times \vec v_c \phi_n 
      = \sum_{n \in \mathcal N} \sum_{c \in \mathcal C(n)} \ell_{nc} \vec v_c \times \vec n_{nc} \phi_n \\
      &= \sum_{c \in \mathcal C}  \vec v_c \times\sum_{n \in \mathcal N(c)} \ell_{nc} \vec n_{nc} \phi_n 
      %= \sum_{c \in \mathcal C}  \vec v_c \times \vec G_c \phi |c| \\
      = -\sum_{c \in \mathcal C}  \vec G_c \phi \times \vec v_c |c| 
    \end{align}
  \item
 For any function $\phi$ with node values $\phi_n$ one has
 \begin{align}
  \vec C_n (\vec G \phi) &= -\frac{1}{|c_n|} \sum_{c \in \mathcal C(n)} \ell_{nc} \vec n_{nc} \times \left( \frac{1}{|c|} \sum_{m \in \mathcal N(c)} \ell_{mc} \vec n_{mc} \phi_m \right ). \label{eq:curlgrad}
 \end{align}
 This is a linear combination of a finite number of values $\phi_m$. Consider the terms involving $\phi_n$ first (i.e. those for which $m = n$). Observe that $n$ is always contained in $\mathcal N(c)$ if $c \in \mathcal C(n)$, and
 \begin{align}
  -\frac{1}{|c_n|} \sum_{c \in \mathcal C(n)} \ell_{nc} \vec n_{nc} \times \frac{1}{|c|} \ell_{nc} \vec n_{nc}  = 0.
 \end{align}
 Consider now $m \neq n$, but such that it shares an edge with $n$, i.e. $m$ such that $\mathcal N(e) = \{n, m\}$ for some $e$ (see Figure \ref{fig:area}). The summation over $c \in \mathcal C(n)$ is reduced to two cells $c_1, c_2$ that share the edge $e$: 
 \begin{align}
 \frac{1}{|c_1|}  \ell_{nc_1} \vec n_{nc_1} \times \ell_{mc_1} \vec n_{mc_1}  +   \frac{1}{|c_2|}  \ell_{nc_2} \vec n_{nc_2} \times \ell_{mc_2} \vec n_{mc_2}.
 \end{align}
 For one of the cells, the two nodes $n,m$ are ordered counterclockwise, for the other clockwise. By Lemma \ref{lemma:area} both terms are $\pm\frac12$ but differ by a sign, and their sum thus vanishes.
 
 Consider finally, in the case of a quadrangular cell, a node $m$ that does not share an edge with $n$. Then,
 \begin{align}
  \ell_{mc} \vec n_{mc} = - \ell_{nc} \vec n_{nc} \label{eq:oppositenode}
 \end{align}
 
 Thus, the prefactors of all $\phi_m $ in Equation \eqref{eq:curlgrad} vanish.
  \end{enumerate}
% \hfill $\Box$

\end{proof}

As long as $\mathcal E(c) \leq 4$, we thus have identified a discrete node-centered curl operator $\vec C_n$ that identically annihilates the discrete cell-centered gradient $\vec G_c$. As the update of the velocity (Equation \eqref{eq:pstarsolverevoutionofugeneral}) is precisely this gradient applied to $p^*$, the vorticity $\vec C_n \vec v$ is stationary, i.e. a discrete involution.

Note that for cells $c$ with $\mathcal E(c) > 4$ one can always find several nodes $m$ that are not sharing an edge with $n$, but \eqref{eq:oppositenode} is generally no longer true, and neither is \eqref{eq:curlgrad}. This does not mean that there does not exist a discrete curl that would annihilate $\vec G$, but rather that its expression is different from $\vec C$ in the general case. From the analysis at the beginning of this Section one can also conclude that, generically, there will not be enough discrete vorticities. One therefore cannot, in this case, speak of a vorticity preserving method.

Both the statement of Theorem \ref{thm:curl} and its proof are essentially contained in \cite{morton01}, Section 7, where the authors consider a scheme that uses the discrete gradient $\vec G$ in a discretization of $\del_t \vec v + \nabla p = 0$. In the present work, however, this discrete gradient appears as a consequence of the more fundamental assumption of nodal conservation, and thus promotes the latter to a good guiding principle for the construction of structure preserving numerical methods.

\section{The special case of a Cartesian grid} \label{sec:cartesian}

The first aim of this Section is to ensure that the solver suggested in this work can be compared to other vorticity preserving methods, which are mostly derived for Cartesian grids only. The second aim is to pedagogically revisit the concept of nodal conservation in a simple setup. This is why we do not content ourselves with merely stating the Cartesian version of the relevant formulas, but briefly repeat the arguments leading to their derivation. 

We also for the first time derive the Cartesian version of the solver from Section \ref{ssec:solvervelocity}, whose unstructured version has been published in \cite{gallice22}. This is important for the analysis of its structure preservation properties in Section \ref{app:nonvortpres}.

In this Section, the grid will be Cartesian with cell centers $(x_i, y_j)$ and $x_{i+1} - x_i =: \Delta x$ and $y_{j+1} - y_j =: \Delta y$ for all $i,j \in \mathbb Z$. The interfaces are associated with half-integer indices $i+\frac12$ and $j + \frac12$, respectively and nodes are located at $(x_{i+\frac12},y_{j+\frac12})$.

\subsection{Nodal conservation} \label{ssec:nodeconservationcartesian}

Consider the solvers from Sections \ref{ssec:solvervelocity}--\ref{ssec:solverpressure}. The new aspect of these solvers is to introduce discontinuous fluxes across the middle wave. This is to be interpreted in the following way: cell L is updated using the fluxes $\bar p_\text{L}, \bar u_\text{L}$, while cell R is updated using the fluxes $-\bar p_\text{R}, -\bar u_\text{R}$. Every edge is cut in two subedges, with \textit{a priori} different fluxes associated to each subedge, and each of its L/R sides (see Figure \ref{fig:conservationp}). We could denote a flux through the subedge at $x_{i+\frac12}$ which contains the node $(i+\frac12, j+\frac12)$ as $f |_{i+\frac12,j}^{(i+\frac12,j+\frac12)}$. It is clear, however, that $i+\frac12$ is mentioned twice. Removing a few brackets one thus arrives at the unambiguous, but shorter notation $f|_{i+\frac12,j}^{j+\frac12}$, that is going to be used below.

There are four subedges around the node $(i+\frac12, j + \frac12)$. Global conservation will become local by concentrating on the fluxes through these subedges for every node.

\subsubsection{Solver with nodal velocity} \label{ssec:nodalvelocitysolvercartesian}

Consider first the solver \eqref{ssec:solvervelocity} with a free parameter $u^*$. The values of $u$ in the four cells around the node are only updated via $x$-fluxes, as $u$ does not have a $y$-flux. Recall that the $x$-flux of $u$ is $\bar p_\text{L/R}$. The total change of $u$ due to the four fluxes associated to subedges around the node $(i+\frac12, j + \frac12)$ is:

\begin{figure}
\centering
\includegraphics[width=0.45\textwidth]{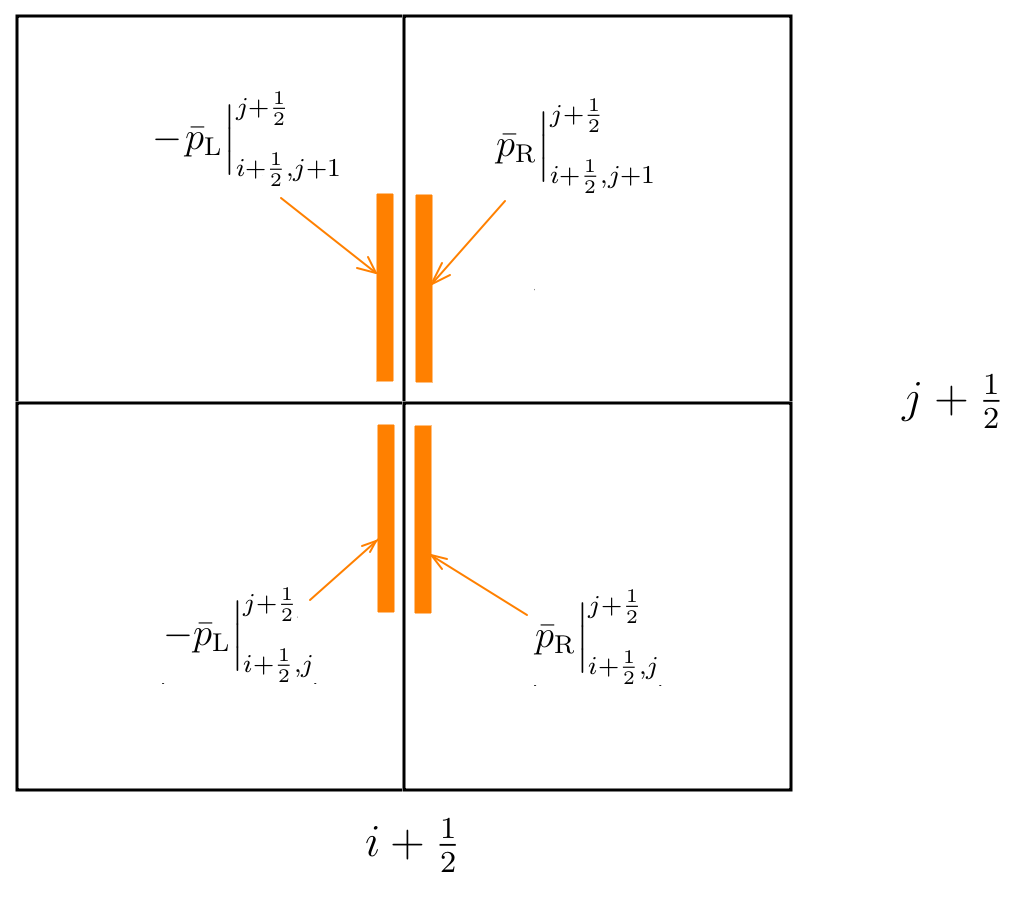}\hfill
\includegraphics[width=0.45\textwidth]{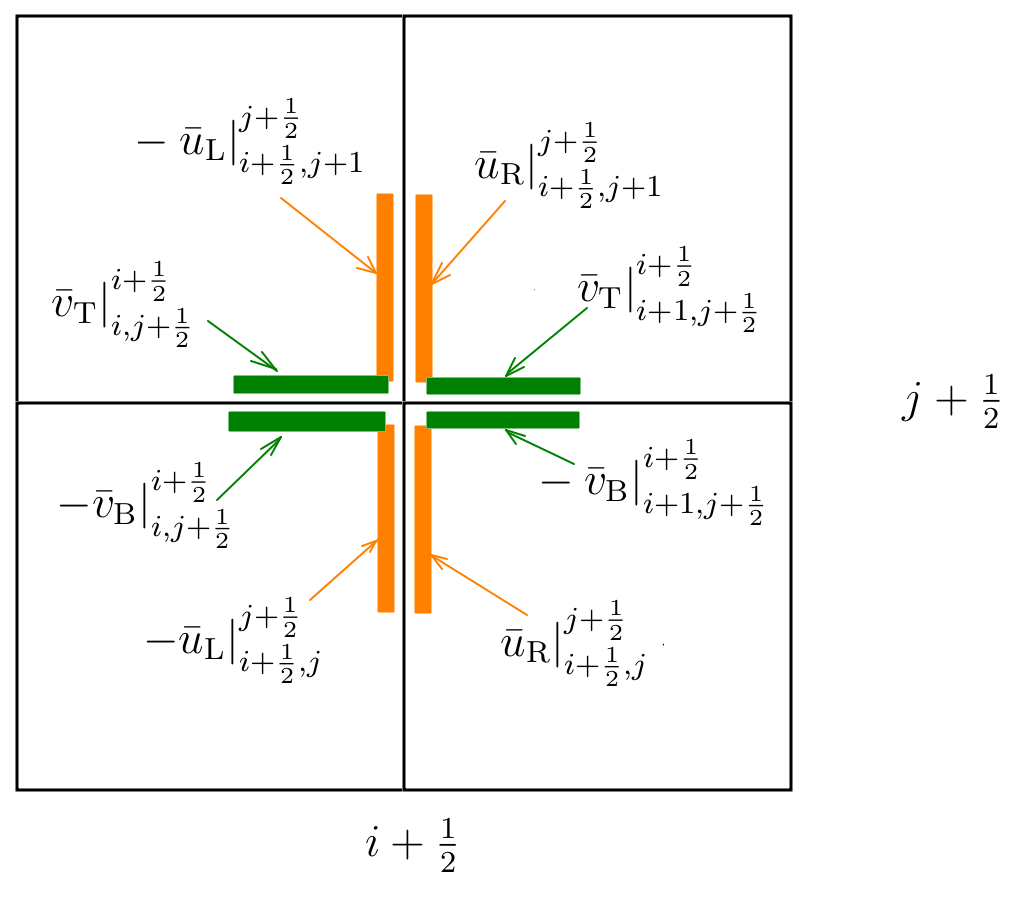}
\caption{\emph{Left}: The fluxes of $u$ involved in the nodal conservation for the solver defined in \ref{ssec:solvervelocity}: The velocity $u^*$ at the node is chosen such that their sum vanishes. \emph{Right}: The fluxes of $p$ involved in the nodal conservation for the solver defined in \ref{ssec:solverpressure}: The pressure $p^*$ at the node is chosen such that their sum vanishes.}
\label{fig:conservationp}
\end{figure}

\begin{align}
 0 &= \underbrace{- \bar p_\text{L} \Big|^{j+\frac12}_{i+\frac12,j}}_{\text{contribution to update of cell} (i,j)} + \underbrace{\bar p_\text{R} \Big|^{j+\frac12}_{i+\frac12,j}}_{\text{update of cell} (i+1,j)}   - \underbrace{\bar p_\text{L} \Big|^{j+\frac12}_{i+\frac12,j+1}}_{\text{update of cell} (i,j+1)} \\
 &+ \underbrace{\bar p_\text{R} \Big|^{j+\frac12}_{i+\frac12,j+1}}_{\text{update of cell} (i+1,j+1)} \nonumber 
 = - \left(p_{ij} - u^*\Big|^{j+\frac12}_{i+\frac12,j} + u_{ij} \right) \\&+ p_{i+1,j} + u^*\Big|^{j+\frac12}_{i+\frac12,j} -u_{i+1,j} -  \left(p_{i,j+1} - u^*\Big|^{j+\frac12}_{i+\frac12,j+1} + u_{i,j+1} \right) \\&\nonumber+ p_{i+1,j+1} + u^*\Big|^{j+\frac12}_{i+\frac12,j+1} -u_{i+1,j+1} \\
 &= 2 u^*\Big|^{j+\frac12}_{i+\frac12,j} + 2 u^*\Big|^{j+\frac12}_{i+\frac12,j+1} + \{ [p]_{i+\frac12} \}_{j+\frac12} - \{ \{ u \}_{i+\frac12} \}_{j+\frac12} 
\end{align}

In the last line, a short-hand notation has been introduced, with $\{ q \}_{i+\frac12} := q_{i+1} + q_i$, $[q]_{i+\frac12} = q_{i+1} - q_i$, and combinations of these. For more details on the notation, see Section 6 in \cite{barsukow21yee}.

Assume now that $u^*\Big|^{j+\frac12}_{i+\frac12,j} = u^*\Big|^{j+\frac12}_{i+\frac12,j+1} =: u^*_{i+\frac12,j+\frac12}$ (thus introducing a velocity associated to a node rather than a subedge), which yields
\begin{align}
 u^*_{i+\frac12,j+\frac12} &= \frac{\{ \{ u \}_{i+\frac12} \}_{j+\frac12}}{4}  - \frac12 \frac{\{ [p]_{i+\frac12} \}_{j+\frac12} }{2} \label{eq:ustardefinition}
\end{align}
An analogous computation gives the perpendicular result
\begin{align}
 v^*_{i+\frac12,j+\frac12} &= \frac{\{ \{ v \}_{i+\frac12} \}_{j+\frac12}}{4}  - \frac12 \frac{[\{ p \}_{i+\frac12}]_{j+\frac12} }{2}
 \label{eq:vstardefinition}
\end{align}

\subsubsection{Solver with nodal pressure}\label{ssec:nodalpressuresolvercartesian}

For the solver defined in \ref{ssec:solverpressure} with a free parameter $p^*$, conservation of $u$ is automatic across the subedge, and so is node-based conservation. Conservation of $p$ involves fluxes $\bar u_\text{L/R}$ across $x$-edges and fluxes $\bar v_\text{T/B}$ across $y$-edges, originating from the $90^\circ$ rotated solver
\begin{align}
 -(u^*_\text{B} - u_\text{B}) &= 0  & & &  u_\text{T} - u^*_\text{T} &= 0\\
 -(v^*_\text{B} - v_\text{B}) &= \bar p_\text{B} - p_\text{B}  &   \bar p_\text{B} &= \bar p_\text{T} &  v_\text{T} - v^*_\text{T} &= p_\text{T} - \bar p_\text{T}\\
 -(p^*_\text{B} - p_\text{B}) &= \bar v_\text{B} - v_\text{B}  & \bar v_\text{B} &\cancel{=} \bar v_\text{T} &  p_\text{T} - p^*_\text{T} &= v_\text{T} - \bar v_\text{T}
\end{align}

One finds
\begin{align}
 \bar u_\text{L} &= -p^* + p_\text{L} + u_\text{L}  & \bar u_\text{R}  &= - p_\text{R} + p^* + u_\text{R} \label{eq:pressurefluxesx} \\
 \bar v_\text{B} &= -p^* + p_\text{B} + v_\text{B}  & \bar v_\text{T} &= -p_\text{T} + p^* + v_\text{T}  \label{eq:pressurefluxesy}
\end{align}
such that node-based conservation around $(i+\frac12,j+\frac12)$ reads
\begin{align}
 \Delta x \left( - \bar u_\text{L} |_{i+\frac12,j}^{j+\frac12} + \bar u_\text{R} |_{i+\frac12,j}^{j+\frac12} - \bar u_\text{L} |_{i+\frac12,j+1}^{j+\frac12} + \bar u_\text{R} |_{i+\frac12,j+1}^{j+\frac12} \right )\\
 + \Delta y \left( - \bar v_\text{B} |_{i,j+\frac12}^{i+\frac12} + \bar v_\text{T} |_{i,j+\frac12}^{i+\frac12} - \bar v_\text{B} |_{i+1,j+\frac12}^{i+\frac12} + \bar v_\text{T} |_{i+1,j+\frac12}^{i+\frac12} \right )= 0
\end{align}

Inserting \eqref{eq:pressurefluxesx}--\eqref{eq:pressurefluxesy} yields
\begin{align}
0&= 2 \frac{p^*|_{i+\frac12,j}^{j+\frac12} + p^*|_{i+\frac12,j+1}^{j+\frac12} }{\Delta x} + 2 \frac{ p^*|_{i,j+\frac12}^{i+\frac12}+ p^*|_{i+1,j+\frac12}^{i+\frac12} }{\Delta y}\\
\nonumber & - \frac{p_{ij} + p_{i+1,j} + p_{i,j+1} + p_{i+1,j+1}}{\Delta x} - \frac{p_{ij} +p_{i,j+1} + p_{i+1,j}+p_{i+1,j+1} }{\Delta y} \\
\nonumber & + \frac{- u_{ij}   + u_{i+1,j}  - u_{i,j+1}  + u_{i+1,j+1} }{\Delta x}  + \frac{- v_{ij} + v_{i,j+1}  - v_{i+1,j}  + v_{i+1,j+1}  }{\Delta y}
\end{align}
Assume now that $p^*$ is the same at all the subedges adjacent to a node. Then
\begin{align}
 p^*_{i+\frac12,j+\frac12} &= \frac14 \{ \{ p \}_{i+\frac12} \}_{j+\frac12} - \frac12 \frac{\left( \frac{\{ [u]_{i+\frac12} \}_{j+\frac12}}{2 \Delta x} + \frac{[\{v \}_{i+\frac12} ]_{j+\frac12}}{2\Delta y} \right)}{\frac{1}{\Delta x} + \frac{1}{\Delta y}} \label{eq:pstardefinition}
\end{align}
or, if $\Delta y = \Delta x$
\begin{align}
 p^*_{i+\frac12,j+\frac12} &= \frac14 \{ \{ p \}_{i+\frac12} \}_{j+\frac12} - \frac14 \left( \frac{\{ [u]_{i+\frac12} \}_{j+\frac12}}{2} + \frac{[\{v \}_{i+\frac12} ]_{j+\frac12}}{2} \right) 
\end{align}
Observe that the term in round brackets is the node-based divergence.

\subsection{Update of the variables}

\subsubsection{Solver with nodal velocity}

Let us collect the fluxes for the update of $u$ for the solver of Section \ref{ssec:solvervelocity} (with $u^*,v^*$ given by \eqref{eq:ustardefinition}--\eqref{eq:vstardefinition}). 

\begin{align}
 2\Delta x\frac{\dd}{\dd t} u_{ij} &= - \bar p_\text{L} \Big|^{j+\frac12}_{i+\frac12,j} - \bar p_\text{L} \Big|^{j-\frac12}_{i+\frac12,j} + \bar p_\text{R} \Big|_{i-\frac12,j}^{j+\frac12} + \bar p_\text{R} \Big|_{i-\frac12,j}^{j-\frac12} \\
 &=  - 4u_{ij} + u^*_{i+\frac12,j+\frac12} + u^*_{i+\frac12,j-\frac12} +  u^*_{i-\frac12,j+\frac12}  + u^*_{i-\frac12,j-\frac12} \\
 &=  - 4u_{ij} + \{ \{ u^*\}_{i\pm\frac12} \}_{j\pm\frac12} \\
 &= - 4 \left( u_{ij} + \frac{\{\{ \{\{ u \}\}_{i\pm\frac12} \}\}_{j\pm\frac12}}{16} \right)  - \frac{\{\{ [p]_{i\pm1} \}\}_{j\pm\frac12} }{4} \label{eq:nodalevelocitysolverveleq}
\end{align}
The leading factor of 2 is due to subedges being involved, instead of full edges. The equation for $v$ is analogous.

\begin{align}
 2\frac{\dd}{\dd t} p_{ij} &= \frac{- \bar u_\text{L} \Big|^{j+\frac12}_{i+\frac12,j} - \bar u_\text{L} \Big|^{j-\frac12}_{i+\frac12,j} + \bar u_\text{R} \Big|_{i-\frac12,j}^{j+\frac12} + \bar u_\text{R} \Big|_{i-\frac12,j}^{j-\frac12} }{\Delta x}+ \frac{\text{perpendicular terms}}{\Delta y}\\
 &=  \frac{- u^*_{i+\frac12,j+\frac12} - u^*_{i+\frac12,j-\frac12} + u^*_{i-\frac12,j+\frac12} + u^*_{i-\frac12,j-\frac12}}{\Delta x} + \frac{\text{terms with $v^*$} }{\Delta y}\\
 &= - \frac{\{[u^*]_{i\pm\frac12}\}_{j\pm\frac12}}{\Delta x} - \frac{[\{v^*\}_{i\pm\frac12} ]_{j\pm\frac12} }{\Delta y}\\
 &= -\left(\frac{\{\{ [ u ]_{i\pm1} \}\}_{j+\frac12}}{4 \Delta x} + \frac{ [ \{\{ v \}\}_{i+\frac12}]_{j\pm1}}{4 \Delta y} \right ) + \frac12 \frac{\{\{ [[p]]_{i\pm\frac12} \}\}_{j\pm\frac12} }{2 \Delta x} \\&\nonumber+  \frac12 \frac{ [[ \{\{p \}\}_{i\pm\frac12}]]_{j\pm\frac12} }{2 \Delta y}
\end{align}

This is the Cartesian version of the solver derived in \cite{gallice22}.

\subsubsection{Solver with nodal pressure} \label{ssec:solvernodalpresurecartesian}

The update of $p$ for the solver of Section \ref{ssec:solverpressure} involves the following fluxes:
\begin{align}
  \frac{\dd}{\dd t} p_{ij} &= \frac{- \bar u_\text{L} |_{i+\frac12,j}^{j+\frac12} + \bar u_\text{R} |_{i-\frac12,j}^{j+\frac12} - \bar u_\text{L} |_{i+\frac12,j}^{j-\frac12} + \bar u_\text{R} |_{i-\frac12,j}^{j-\frac12} }{\Delta x}\\
 &\nonumber +\frac{- \bar v_\text{B} |_{i,j+\frac12}^{i+\frac12} + \bar v_\text{T} |_{i,j-\frac12}^{i+\frac12} - \bar v_\text{B} |_{i,j+\frac12}^{i-\frac12} + \bar v_\text{T} |_{i,j-\frac12}^{i-\frac12}  }{\Delta y}
\end{align}

\begin{figure}
 \centering
 \includegraphics[width=0.65\textwidth]{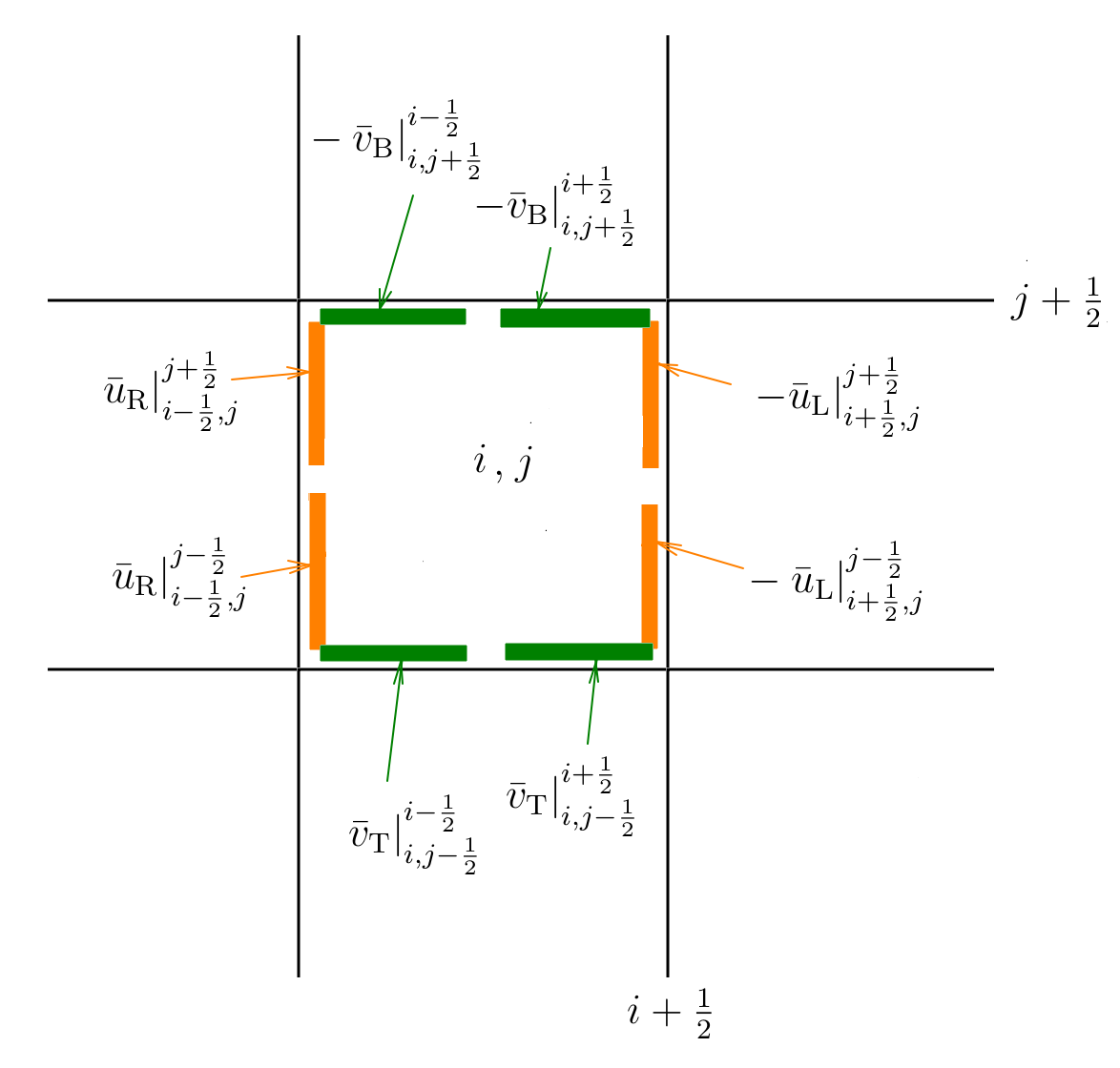}
 \caption{The fluxes involved in the update of cell $(i, j)$.}
 \label{fig:fluxeshalfedges}
\end{figure}

Inserting \eqref{eq:pressurefluxesx}--\eqref{eq:pressurefluxesy} and the expression \eqref{eq:pstardefinition} of $p^*$ yields
\begin{align}
 \frac{\dd}{\dd t} p_{ij} &= \frac12 \left( \frac1{\Delta x} + \frac1{\Delta y}  \right )\left( \{\{ p^* \}_{i\pm\frac12} \}_{j\pm\frac12} - 4 p_{ij} \right) \\
 &= - \frac18 \left( \frac{\{\{ [u]_{i\pm1} \}\}_{j\pm\frac12}}{\Delta x} +\frac{ [\{\{v \}\}_{i\pm\frac12} ]_{j\pm1}}{\Delta y} \right) \\&+\nonumber \frac12 \left( \frac1{\Delta x} + \frac1{\Delta y}  \right ) \left( \frac14 \{\{ \{\{ p \}\}_{i\pm\frac12} \}\}_{j\pm\frac12} - 4 p_{ij} \right )\label{eq:pstarsolverevoutionofp}
\end{align}
If $\Delta y = \Delta x$, one finds
\begin{align}
\Delta x\frac{\dd}{\dd t} p_{ij} &= - \frac18 \left( \{\{ [u]_{i\pm1} \}\}_{j\pm\frac12} + [\{\{v \}\}_{i\pm\frac12} ]_{j\pm1} \right) + \frac14 \{\{ \{\{ p \}\}_{i\pm\frac12} \}\}_{j\pm\frac12} - 4 p_{ij}
 \end{align}
The evolution equation for $p$ involves the divergence obtained as an average of 4 node-based divergences and some (non-standard) diffusion operator on $p$.

The evolution equation for $u$ reads
\begin{align}
 2 \Delta x \frac{\dd}{\dd t} u_{ij} &= - \bar p_\text{L} |_{i+\frac12,j}^{j+\frac12} + \bar p_\text{R} |_{i-\frac12,j}^{j+\frac12} - \bar p_\text{L} |_{i+\frac12,j}^{j-\frac12} + \bar p_\text{R} |_{i-\frac12,j}^{j-\frac12}\\
 &= - p^*_{i+\frac12,j+\frac12} + p^*_{i-\frac12,j+\frac12} - p^*_{i+\frac12,j-\frac12} + p^*_{i-\frac12,j-\frac12}\\
 \Delta x \frac{\dd}{\dd t} u_{ij}  &= - \frac12 \{ [ p^* ]_{i\pm\frac12} \}_{j\pm\frac12} \label{eq:uevolution}\\
 &= - \frac18 \{\{ [ p ]_{i\pm1} \}\}_{j\pm\frac12} + \frac12 \frac{  \frac{\{\{ [[u]]_{i\pm\frac12} \}\}_{j\pm\frac12}}{4\Delta x} + \frac{[[v ]_{i\pm1} ]_{j\pm1}}{4\Delta y} }{\frac1{\Delta x} + \frac{1}{\Delta y}} \label{eq:pstarsolverevoutionofu}
\end{align}
For $\Delta y = \Delta x$ this simplifies to
\begin{align}
 \Delta x \frac{\dd}{\dd t} u_{ij}  &= - \frac18 \{\{ [ p ]_{i\pm1} \}\}_{j\pm\frac12} + \frac14  \left( \frac{\{\{ [[u]]_{i\pm\frac12} \}\}_{j\pm\frac12}}{4} + \frac{[[v ]_{i\pm1} ]_{j\pm1}}{4} \right)
\end{align}

By analogy,
\begin{align}
 \Delta y \frac{\dd}{\dd t} v_{ij}  &= - \frac12 [\{  p^* \}_{i\pm\frac12} ]_{j\pm\frac12} \label{eq:vevolution}
\end{align}

The prefactor in front of the diffusion is only $\frac14$ instead of $\frac12$, therefore the stability condition of the scheme can be expected to be $\text{CFL}< \frac12$.

The diffusion in the velocity equation \eqref{eq:pstarsolverevoutionofu} and the divergence in the pressure equation \eqref{eq:pstarsolverevoutionofp} all are averages and differences of the nodal divergence appearing in $p^*$. This is the essential ingredient for stationarity and vorticity preservation of the method. One easily observes that the discrete vorticity
\begin{align}
 \frac{\dd}{\dd t} \left( \frac{[\{ u \}_{i+\frac12}]_{j+\frac12}}{2\Delta y}  -  \frac{\{ [ v]_{i+\frac12}\}_{j+\frac12}}{2\Delta x}  \right) = 0 \label{eq:discretevorticity}
\end{align}
remains stationary.

\section{Extension to second order of accuracy} \label{sec:secondorder}
\subsection{General case}

Consider in every cell $c$ a polynomial reconstruction $q_{c,\text{r}} = \vecc{\vec v_{c,\text{r}}}{p_{c,\text{r}}}$ in the space $V = \mathrm{span}_{\mathbb R}(1, x, y)$ for each variable. With the ansatz $q_{c,\text{r}}(x,y) = q_c + a_1 x + a_2 y$ the reconstruction is automatically conservative if the cell centroid is chosen as the origin of coordinates $x$ and $y$, i.e. $\frac{1}{|c|}\int_{c}  q_{ij,\text{r}}(\vec x)  \, \dd \vec x = q_{ij}^n$.

It would be desirable to determine the missing parameters $a_1,a_2$ in the reconstruction in such a way that its averages over cells in a stencil $S = \{ c' : c' \text{ neighbour of } c \}$ covering the neighbourhood of $c$ agree with the known cell averages $q_{c'}^n$. This is, in general, not possible because the reconstruction only has two parameters left. Following \cite{barth89,maire09}, we propose to use the best approximation in the sense of least squares, i.e. the reconstruction in $V$ which minimizes
\begin{align}
 \sum_{c \in S} \left( q_c^n  - \frac{1}{|c|}\int_c q_{c,\text{r}}(x, y) \, \dd x\, \dd y  \right)^2 \label{eq:reconminim}
\end{align}
Minimization of \eqref{eq:reconminim} gives rise to a $2 \times 2$ linear system that can be easily solved. For $S$ one can consider either the cells adjacent to all edges of $c$, i.e.
\begin{align}
 S_{\mathcal E} := \{ c' : \mathcal E(c) \cap \mathcal E(c') \neq \emptyset\}
\end{align}
or the cells adjacent to nodes of $c$:
\begin{align}
 S_{\mathcal N} := \{ c' : \mathcal N(c) \cap \mathcal N(c') \neq \emptyset\}.
\end{align}

We then compute the value of $p^*$ according to formula \eqref{eq:nodalpressureunstructuredrawdiv}, but replacing the values of the cell averages by the values of the reconstructions at the location of the node:
\begin{align}
  p^*_{n}  &= \frac{\displaystyle\sum_{s \in \mathcal{SE}(n)} |s|  \frac{p_{\text{R,r}}(\vec x_n - \vec x_{\text{R}}) + p_{\text{L,r}}(\vec x_n - \vec x_{\text{L}})}{2}  -  \frac12 |c_n| D_n \vec v_\text{r}}{\displaystyle\sum_{s \in \mathcal{SE}(n)} |s|}  \\
  D_n \vec v_\text{r} &:= -\frac{1}{|c_n|} \sum_{c \in \mathcal C(n)} \ell_{nc} \vec n_{nc} \cdot \vec v_{c,\text{r}}(\vec x_n - \vec x_c)
\end{align}
where $\vec x_n$ and $\vec x_{\text{R/L}}$ denote the locations of the node and of the cell centroids of the left and right cell (with respect to $\vec n_s$), respectively.

This allows to compute all states of the Riemann solver according to \eqref{eq:pstarstatesu}, replacing again the cell averages with the values of the reconstructions evaluated at the location of the node. (Note that at second order of accuracy the precise location does not matter, and no quadrature along the edge is needed.) The velocity equation \eqref{eq:pstarsolverevoutionofugeneral} simply uses the new nodal pressure $p^*_n$
\begin{align}
 \frac{\dd}{\dd t} \vec v_c &= - \frac{1}{|c|} \sum_{n \in \mathcal N(c)} \sum_{s \in \mathcal{SE}(n, c)} |s| \vec n_{s, c} p^*_n \label{eq:highordervevolution}
\end{align}
while the evolution of $p$ reads (compare to \eqref{eq:pstarsolverevoutionofpgeneral})
\begin{align}
 \frac{\dd}{\dd t} p_c &=  - \frac{1}{|c|} \sum_{n \in \mathcal N(c)} \sum_{s \in \mathcal{SE}(n, c)} |s| (\vec v_{c,\text{r}}(\vec x_n - \vec x_c) \cdot \vec n_{s,c} - p^*_n + p_{c,\text{r}}(\vec x_n - \vec x_c))  
\end{align}

Contrary to \eqref{eq:pstarsolverevoutionofpgeneral} now the contribution of $\vec v_c$ does not disappear, because the values of $\vec v_{c,\text{r}}$ are not the same at all nodes. Write $\vec v_{c,\text{r}}(\vec x) = \vec v_c + \vec w_c^x x + \vec w_c^y y$, where $\vec w_c^x, \vec w_c^y$ are the components of the gradient $\vec w_c$ that are obtained using the least squares procedure above. Then, 
\begin{align}
&\frac{1}{|c|} \sum_{n \in \mathcal N(c)} \sum_{s \in \mathcal{SE}(n, c)} |s| (\vec v_{c,\text{r}}(\vec x_n - \vec x_c) \cdot \vec n_{s,c} )\\
&\nonumber= \frac{1}{|c|} \sum_{n \in \mathcal N(c)} \sum_{s \in \mathcal{SE}(n, c)} |s| \Big(\vec v_c\cdot \vec n_{s,c} + \vec w_c^x \cdot \vec n_{s,c} (x_n-x_c) + \vec w_c^y\cdot \vec n_{s,c} (y_n-y_c) \Big) \\
&\overset{\text{Def. \ref{def:nodenormal}}}{=} \frac{1}{|c|} \sum_{n \in \mathcal N(c)}  \ell_{nc}  \Big(\vec v_c\cdot \vec n_{nc} + \vec w_c^x \cdot \vec n_{nc} (x_n-x_c) + \vec w_c^y\cdot \vec n_{nc} (y_n-y_c) \Big)\\
&\overset{\text{Thm. \ref{thm:normal}(\ref{it:thmnormal1})}}{=}\frac{1}{|c|} \sum_{n \in \mathcal N(c)}  \ell_{nc}  \Big(\vec w_c^x \cdot \vec n_{nc} x_n + \vec w_c^y\cdot \vec n_{nc} y_n \Big)\\
&=\frac{1}{|c|} \sum_{n \in \mathcal N(c)}  \ell_{nc} \mathrm{tr}\Big (( \vec x_n \otimes \vec n_{nc} )\vec w_c \Big)\\
&\overset{\text{Thm. \ref{thm:normal}(\ref{it:thmnormal2})}}{=} \mathrm{tr\,}\vec w_c = \div \vec v_{c,\text{r}}
\end{align}

It is clear that the passage from a first-order to a second-order method does not affect the property of vorticity preservation: Equation \eqref{eq:highordervevolution} is still the discrete gradient $\vec G_n$ applied to some modified nodal pressure. In Section \ref{sec:statioanritypreservation} a discrete curl operator that annihilates this discrete gradient on grids with $\#\mathcal E(c) \leq 4 \, \forall c$ is explicitly constructed. The high-order method therefore is vorticity preserving in this case.

The stationarity preservation property of the high-order method is tightly linked to it being vorticity preserving. This time, however, the discrete divergence kept stationary by the method is not obvious. It would be tempting to choose the reconstruction in such a way that its divergence $\div \vec v_{c,\text{r}}$ would be proportional to the discrete divergence $D_n \vec v_\text{r}$ appearing in $p^*$. Unfortunately, it does not seem feasible to generally link the two as long as the slope is found through a least-squares procedure. This does not mean that there is no stationary discrete divergence, but merely that we are unable to generally identify its functional form. For the special case of Cartesian grids, in Section \ref{ssec:cartesianhighorder} the discrete stationary states are explicitly characterized.

We observe that on grids containing pentagons and hexagons, our second-order method with this particular reconstruction is unstable if the simulation runs over very long times. In certain cases we were able to confirm the finding theoretically by performing a von Neumann stability analysis. On quadrangular-triangular grids we have not observed any instability.

\subsection{Cartesian case} \label{ssec:cartesianhighorder}

On Cartesian grids, the minimization problem reads
\begin{align}
 \sum_{c \in S} \left( q_c^n  - \frac{1}{\Delta x}\int_{-\frac{\Delta x}{2}}^{\frac{\Delta x}{2}} \frac{1}{\Delta y}\int_{-\frac{\Delta y}{2}}^{\frac{\Delta y}{2}} q_{ij,\text{r}}(x, y) \, \dd y\, \dd x  \right)^2 \label{eq:cartesianreconminim}
\end{align}
with
\begin{align}
 q_{ij,\text{r}} \colon \left[-\frac{\Delta x}{2},\frac{\Delta x}{2}\right] \times \left[-\frac{\Delta y}{2},\frac{\Delta y}{2}\right] \to \mathbb R^m
\end{align}
such that $x, y$ take their origin in the centroid of $c$.

The stencil $S_{\mathcal E}$ that involves only those neighbours of $c$ with which $c$ shares an edge is the 5-point stencil
\begin{center}$S_5 :=$\begin{tabular}{c|c|c}
 & $(i,j+1)$ & \\\hline
 $(i-1,j)$ & $(i,j)$ & $(i+1,j)$ \\\hline
 & $(i,j-1)$ &
\end{tabular}\end{center}

$S_{\mathcal N}$ is the 9-point stencil 
\begin{center}$S_9 := $\begin{tabular}{c|c|c}
 $(i-1,j+1)$ & $(i,j+1)$ & $(i+1,j+1)$ \\\hline
 $(i-1,j)$ & $(i,j)$ & $(i+1,j)$ \\\hline
 $(i-1,j-1)$ & $(i,j-1)$ & $(i+1,j-1)$
\end{tabular}\end{center}

The linear system arising in the minimization of \eqref{eq:cartesianreconminim} can be solved explicitly: $q_{ij,\text{r}} = q_{ij}^n + a_1 x + a_2 y$ with
\begin{align}
 a_1^{(S_5)} &= \frac{q_{i+1,j}^n - q_{i-1,j}^n}{2 \Delta x} , & a_1^{(S_9)} &= \frac{  [ q^n ]_{i\pm1,j-1} + [ q^n ]_{i\pm1,j} + [ q^n ]_{i\pm1,j+1}   }{6 \Delta x}, \\
 a_2^{(S_5)} &= \frac{q_{i,j+1}^n - q_{i,j-1}^n}{2 \Delta y} ,& a_2^{(S_9)} &= \frac{  [   q_{i-1}^n + q_{i}^n + q_{i+1}^n ]_{j\pm1}}{6 \Delta y}.
\end{align}

The pressure at a node is
\begin{align}
 p^*_{i+\frac12,j+\frac12} &= \frac14 \left( p_{ij,\text{r}}\left(\frac{\Delta x}{2},\frac{\Delta y}{2}\right) 
 +p_{i+1,j,\text{r}}\left(-\frac{\Delta x}{2},\frac{\Delta y}{2}\right) \right . \\ &+ \nonumber\left .
 p_{i,j+1,\text{r}}\left(\frac{\Delta x}{2},-\frac{\Delta y}{2}\right)
 +p_{i+1,j+1,\text{r}}\left(-\frac{\Delta x}{2},-\frac{\Delta y}{2}\right) \right ) \\\nonumber
 &- \frac18 \left(- u_{ij,\text{r}}\left(\frac{\Delta x}{2},\frac{\Delta y}{2}\right) 
 +u_{i+1,j,\text{r}}\left(-\frac{\Delta x}{2},\frac{\Delta y}{2}\right)  \right . \\ &- \nonumber
 u_{i,j+1,\text{r}}\left(\frac{\Delta x}{2},-\frac{\Delta y}{2}\right) 
 +u_{i+1,j+1,\text{r}}\left(-\frac{\Delta x}{2},-\frac{\Delta y}{2}\right) 
    \\\nonumber
 &  \left. -v_{ij,\text{r}}\left(\frac{\Delta x}{2},\frac{\Delta y}{2}\right) 
 -v_{i+1,j,\text{r}}\left(-\frac{\Delta x}{2},\frac{\Delta y}{2}\right)  \right . \\ &+ \nonumber \left .
 v_{i,j+1,\text{r}}\left(\frac{\Delta x}{2},-\frac{\Delta y}{2}\right) 
 +v_{i+1,j+1,\text{r}}\left(-\frac{\Delta x}{2},-\frac{\Delta y}{2}\right) 
   \right )
\end{align} % nodalNE in mathematica sheet

The update equations read
\begin{align}
 \frac{\dd}{\dd t} u_{ij} &= -\frac12 \frac{\{ [p^*]_{i\pm\frac12}\}_{j\pm\frac12}}{\Delta x} \label{eq:uevolutionhighorder}, \\
 \frac{\dd}{\dd t} v_{ij} &= -\frac12 \frac{[ \{p^*\}_{i\pm\frac12}]_{j\pm\frac12}}{\Delta y} \label{eq:vevolutionhighorder}.
\end{align}
and
\begin{align}
  \frac{\Delta x^2}{\Delta x/2} \frac{\dd}{\dd t} p_{ij} &= \left(p^*_{i+\frac12,j+\frac12} - p_{ij,\text{r}}\left( \frac{\Delta x}{2} , \frac{\Delta y}{2} \right) - u_{ij,\text{r}}\left( \frac{\Delta x}{2} , \frac{\Delta y}{2} \right)\right) \nonumber\\& + \left(- p_{ij,\text{r}}\left( -\frac{\Delta x}{2} , \frac{\Delta y}{2} \right) + p^*_{i-\frac12,j+\frac12} + u_{ij,\text{r}}\left( -\frac{\Delta x}{2} , \frac{\Delta y}{2} \right)\right)  \nonumber \\
  &+ \left(p^*_{i+\frac12,j-\frac12} - p_{ij,\text{r}}\left( \frac{\Delta x}{2} , -\frac{\Delta y}{2} \right) - u_{ij,\text{r}}\left( \frac{\Delta x}{2} , -\frac{\Delta y}{2} \right)\right) \nonumber \\ +\nonumber & \left(- p_{ij,\text{r}}\left( -\frac{\Delta x}{2} , -\frac{\Delta y}{2} \right) + p^*_{i-\frac12,j-\frac12} + u_{ij,\text{r}}\left( -\frac{\Delta x}{2} ,- \frac{\Delta y}{2} \right)\right) \\\nonumber 
 &+ \left(p^*_{i+\frac12,j+\frac12} - p_{ij,\text{r}}\left( \frac{\Delta x}{2} , \frac{\Delta y}{2} \right) - v_{ij,\text{r}}\left( \frac{\Delta x}{2} , \frac{\Delta y}{2} \right)\right) \\\nonumber &+ \left(-p_{ij,\text{r}}\left( \frac{\Delta x}{2} , -\frac{\Delta y}{2} \right) + p^*_{i+\frac12,j-\frac12} + v_{ij,\text{r}}\left( \frac{\Delta x}{2} , -\frac{\Delta y}{2} \right)\right)  \\
  \nonumber &+ \left(p^*_{i-\frac12,j+\frac12} - p_{ij,\text{r}}\left(- \frac{\Delta x}{2} , \frac{\Delta y}{2} \right) - v_{ij,\text{r}}\left( -\frac{\Delta x}{2} , \frac{\Delta y}{2} \right)\right) \\ &+ \left(-p_{ij,\text{r}}\left( -\frac{\Delta x}{2} , -\frac{\Delta y}{2} \right) + p^*_{i-\frac12,j-\frac12} + v_{ij,\text{r}}\left( -\frac{\Delta x}{2} , -\frac{\Delta y}{2} \right)\right)
\end{align}

The high-order version has the very same equations \eqref{eq:uevolutionhighorder}--\eqref{eq:vevolutionhighorder}, even if the definition of $p^*$ is different. This means that it is vorticity preserving, too, with the same discrete vorticity \eqref{eq:discretevorticity}. Upon constructing the evolution matrices $\mathcal E$ one can observe that the right kernels (governing the stationary states) of the first- and second-order schemes are different.

\section{Time integration} \label{sec:timeintegration}

The semi-discrete methods are integrated in time using standard explicit methods (forward Euler in the first-order case, 2nd-order Runge-Kutta in the second-order case). We have performed von Neumann analysis on Cartesian grids in order to study the stability of the first-order methods. The spectrum of the amplification matrix arising in this analysis was investigated by numerically sampling the condition from \cite{schur17,schur18}, for more details see \cite{miller71} and Appendix A of \cite{abgrall22}. In line with what has been remarked in Section \ref{ssec:solvernodalpresurecartesian} already, the new method with nodal pressure has a maximum CFL number $\frac{\Delta t}{\Delta x}$ of $\frac12$. (The method with a nodal velocity is not vorticity preserving, but has a stability up to $\mathrm{CFL} = 1$.) On unstructured grids, we respect the stability condition using twice the radius $r = \frac{2|c|}{|\del c|}$ of the inscribed circle as the relevant length scale.

\section{Numerical results} \label{sec:numerical}

All numerical results have been run with the nodal pressure-based solver using a CFL of 0.3.

\subsection{Oblique sound wave}

\begin{figure}
 \centering
 \includegraphics[width=0.45\textwidth]{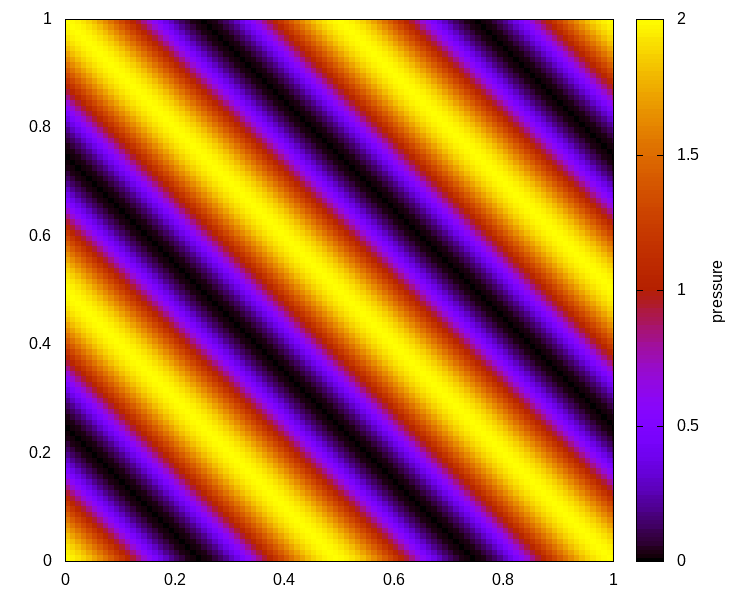}\hfill \includegraphics[width=0.45\textwidth]{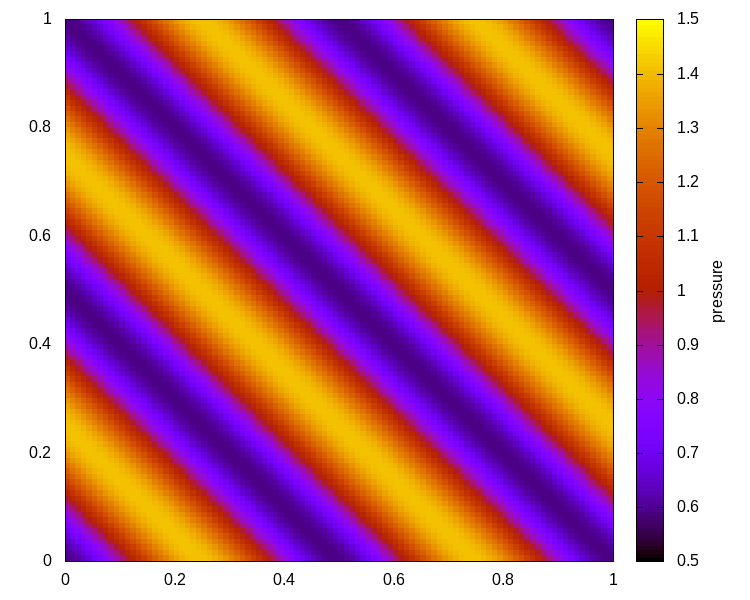}\\
 \includegraphics[width=0.45\textwidth]{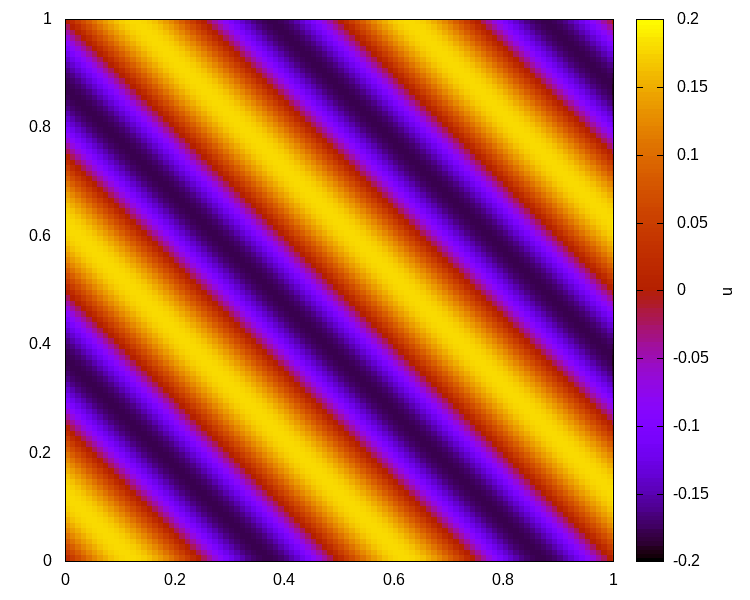}\hfill \includegraphics[width=0.45\textwidth]{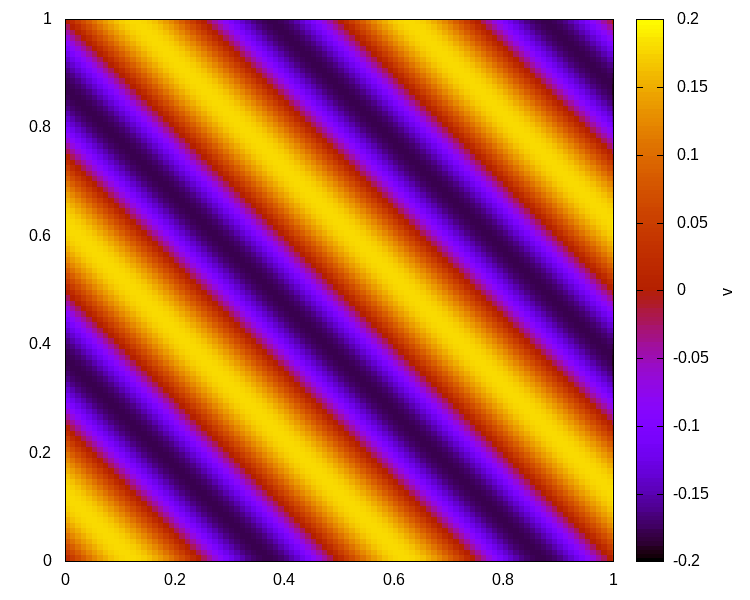}
 \caption{Time evolution of a rotated one-dimensional setup (\emph{Top left}, initial pressure) using the 1st order solver with a nodal pressure variable. \emph{Top right}: Pressure at $t=0.5$. \emph{Bottom left}: Horizontal velocity $u$ at $t=0.5$. \emph{Bottom right}: Vertical velocity $v$ at $t=0.5$.}
 \label{fig:wave1}
\end{figure}

\begin{figure}
 \centering
 \includegraphics[width=0.45\textwidth]{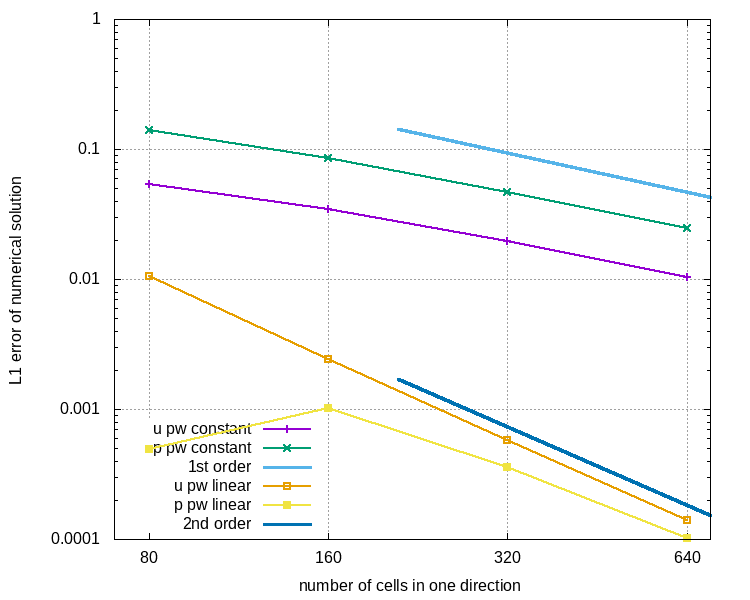}\hfill \includegraphics[width=0.45\textwidth]{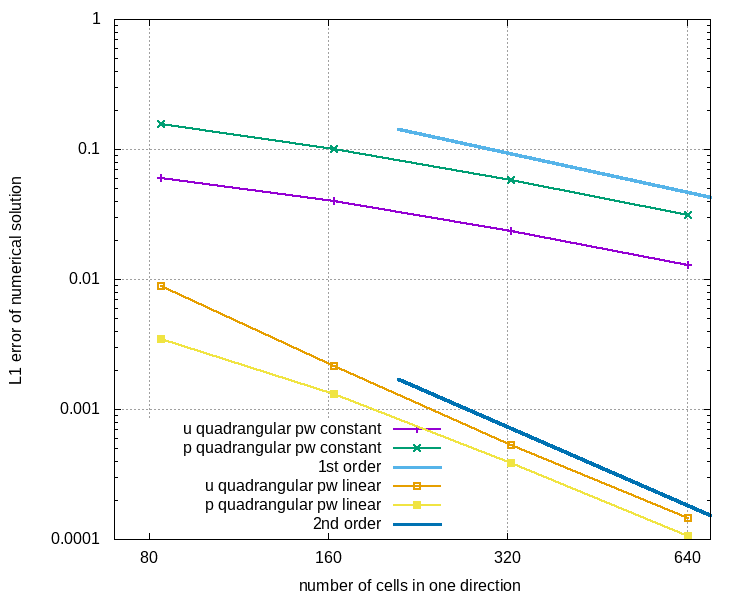}
 \includegraphics[width=0.45\textwidth]{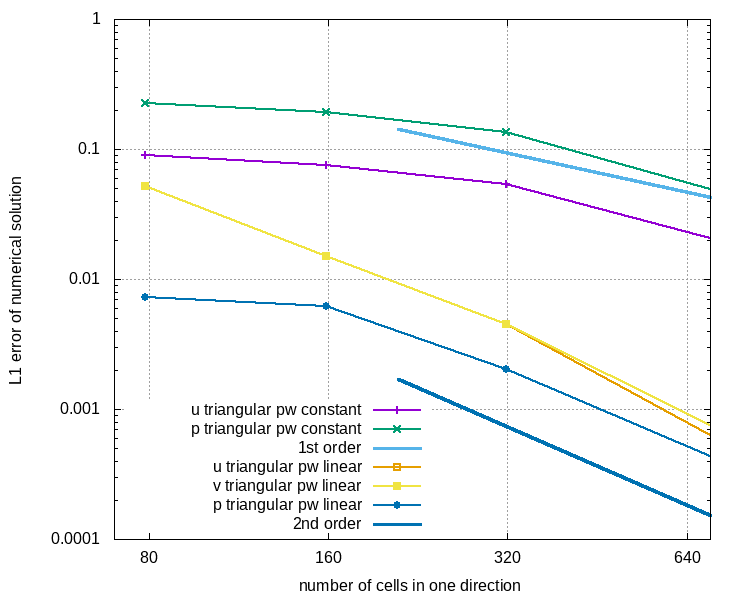}\hfill \includegraphics[width=0.45\textwidth]{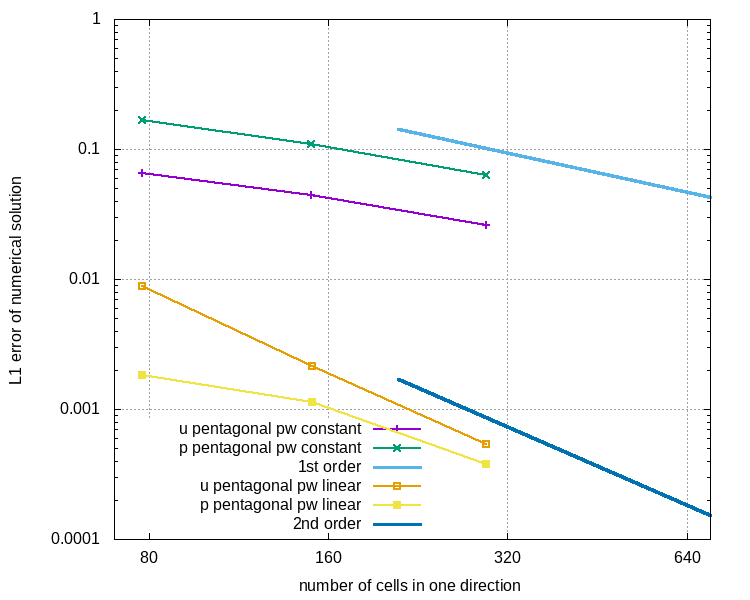}
 \caption{Convergence study at $t=0.5$ for setup \eqref{eq:wave} performed using the solver with a nodal pressure on different grids. The $L^1$ error is shown; as the error curves for the two velocity components are on top of each other only $u$ is shown. \emph{Top left}: Cartesian mesh. The error is shown as function of the number $N$ of grid cells in one direction for $N \times N$ grids. \emph{Top right}: Quadrangular mesh. \emph{Bottom left}: Mixed triangular-quadrangular grid. \emph{Bottom right}: Mixed pentagonal-hexagonal mesh. The unstructured meshes are as in Figure \ref{fig:shockunstructured}. The error is shown as a function of $\sqrt{N}$ for unstructured meshes with a total number $N$ of cells.}
 \label{fig:waveconvergence}
\end{figure}

A way to obtain a family of multi-dimensional setups with exact solutions is to rotate a one-dimensional setup in the $x$-$y$-plane. In one spatial dimension, the initial value problem (IVP) for acoustics can be solved immediately with the method of characteristics. The resulting exact solution allows to perform a convergence analysis. Figure \ref{fig:wave1} shows the initial setup 
\begin{align}
 p_0(x) &=  \cos\left( 2\pi \frac{x}{\lambda \cos(\theta)} \right) & \vec v_0(x) &= 0   \label{eq:wave}
\end{align}
with $\lambda = \frac12$ and rotated by $\theta = \frac{\pi}{4}$. The same Figure also shows its time evolution on a periodic Cartesian grid of $100 \times 100$ cells using the first-order method with nodal pressure. The exact solution at time $t$ is
\begin{align}
    	p(t, x, y) &= \frac12 \left(\cos\left( 2\pi \frac{\xi + t}{\lambda \cos(\theta)} \right) + \cos\left( 2\pi \frac{\xi - t}{\lambda \cos(\theta)} \right) \right) \\
    	u(t, x, y) &= - \frac1{2}\left(\cos\left( 2\pi \frac{\xi + t}{\lambda \cos(\theta)} \right) - \cos\left( 2\pi \frac{\xi - t}{\lambda \cos(\theta)} \right)\right) \cos \theta\\
    	v(t, x, y) &= - \frac1{2}\left(\cos\left( 2\pi \frac{\xi + t}{\lambda \cos(\theta)} \right) - \cos\left( 2\pi \frac{\xi - t}{\lambda \cos(\theta)} \right)\right) \sin \theta
\end{align}
with $ \xi := x \cos \theta + y \sin \theta$. Figure \ref{fig:waveconvergence} shows a convergence study for the method involving a nodal pressure, for both the 1st and 2nd order. The numerical results confirm the theoretical convergence orders. For the 2nd order method, a reconstruction involving $S_{\mathcal N}$, i.e. all neighbours of all nodes of the cell has been used (which amounts to the 9-point stencil on Cartesian grids). The choice of grids is the same as in Figure \ref{fig:shockunstructured}.

\subsection{4-quadrant Riemann problem}

A 4-quadrant Riemann problem for linear acoustics is challenging because certain setups have $L^\infty$-unbounded solutions. One such setup has been studied in detail in \cite{barsukow17} (paralleling the findings in \cite{amadori2015}):
\begin{align}
  p_0(\vec x) &= 0, & u_0(\vec x) &= \begin{cases} 1 & \text{if }x > 0.5 \text{ and }y>0.5, \\ 0 & \text{else,} \end{cases} &   v_0(\vec x) &= 0 .\label{eq:rp}
\end{align}
It has been shown in \cite{barsukow17} that the perpendicular component $v$, initially zero, has a logarithmic singularity at the meeting point $(0.5, 0.5)$ of the four quadrants for all $t> 0$:
\begin{align}
 v(t, x, y) &= \frac{1}{2 \pi} \mathscr L\left( \frac{\sqrt{(x-0.5)^2 + (y-0.5)^2}}{t}  \right )\text{, with} \\
 \mathscr L(s) &:= \ln \frac{1 + \sqrt{1 - s^2}}{s} = - \ln \frac{s}{2} - \frac{s^2}{4} + \mathcal O(s^4)
\end{align}

\begin{figure}
 \centering
 \includegraphics[width=0.45\textwidth]{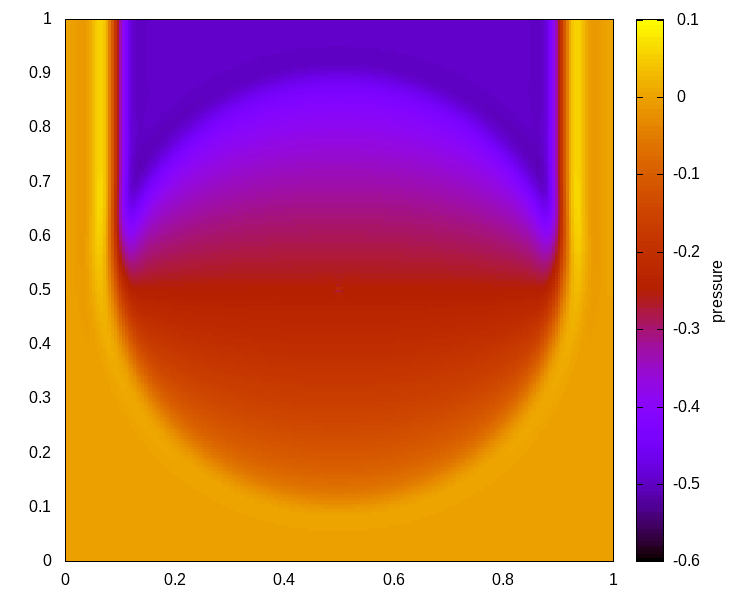}\hfill \includegraphics[width=0.45\textwidth]{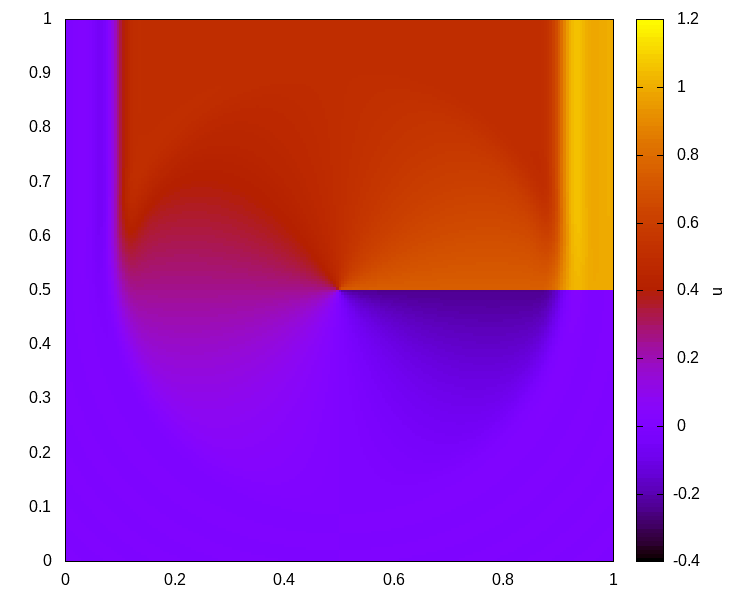}\\
 \includegraphics[width=0.45\textwidth]{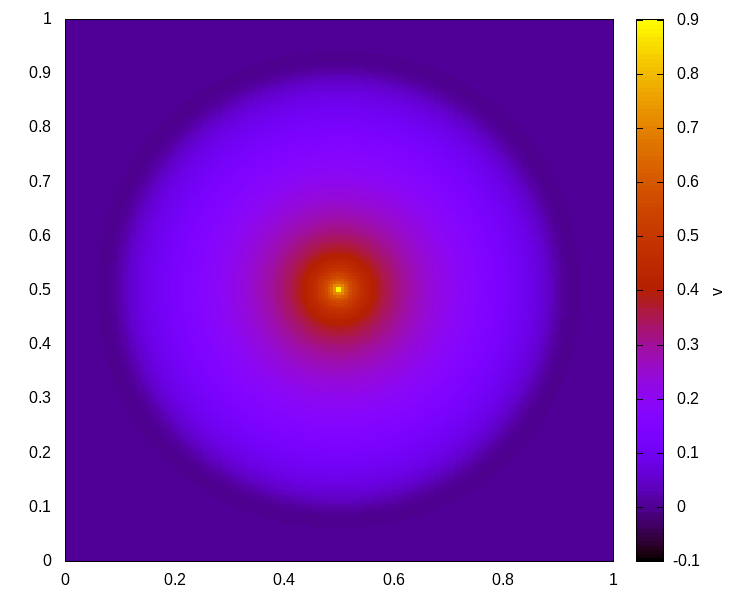}\hfill \includegraphics[width=0.45\textwidth]{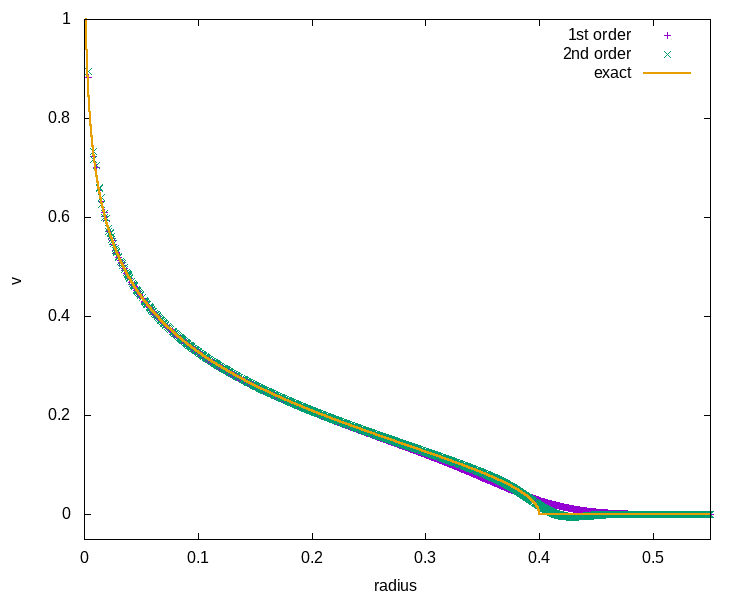}\\
 \includegraphics[width=0.45\textwidth]{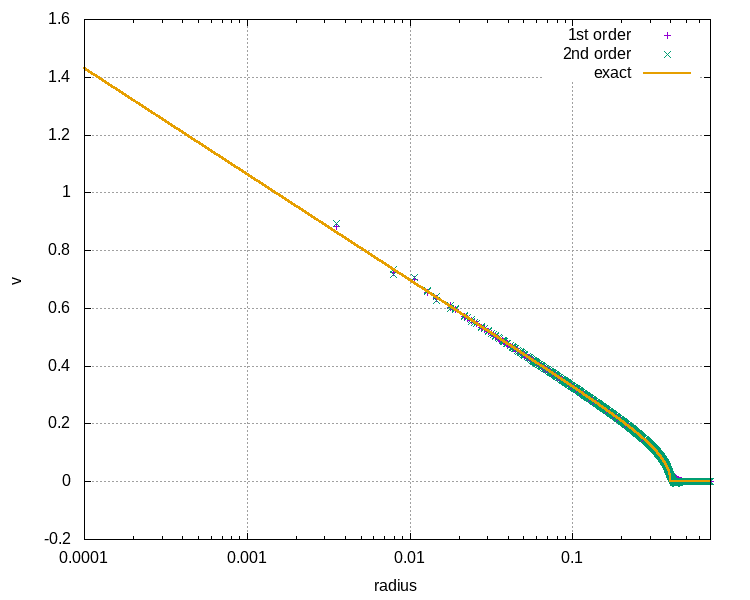}\hfill \includegraphics[width=0.45\textwidth]{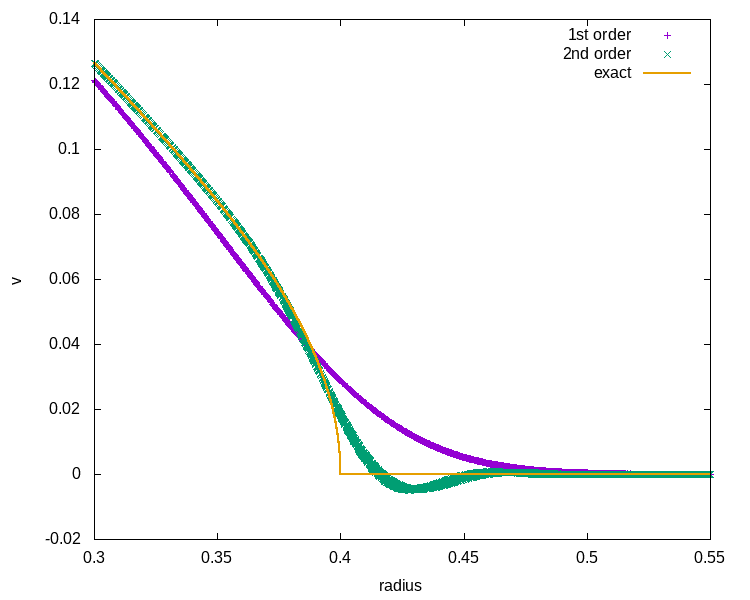}
 \caption{Evolution of the 4-quadrant Riemann problem \eqref{eq:rp}. $p$ (\emph{top left}), $u$ (\emph{top right}) and $v$ (\emph{center left}) at time $t=0.4$, obtained using the second order method with nodal pressure. \emph{Center right}: Radial scatter plot of $v$ for the 1st and 2nd order methods also showing the exact solution. \emph{Bottom}: The same plot, but with the radius shown logarithmically (\emph{left}) and a zoom onto the location $r = t$ (\emph{right}). }
 \label{fig:rp}
\end{figure}

\begin{figure}
 \centering
 \includegraphics[width=0.45\textwidth]{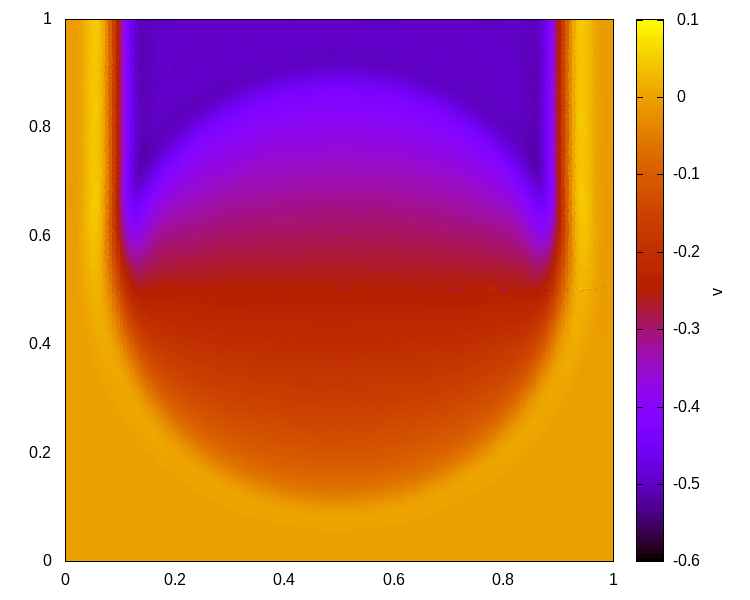}\hfill \includegraphics[width=0.45\textwidth]{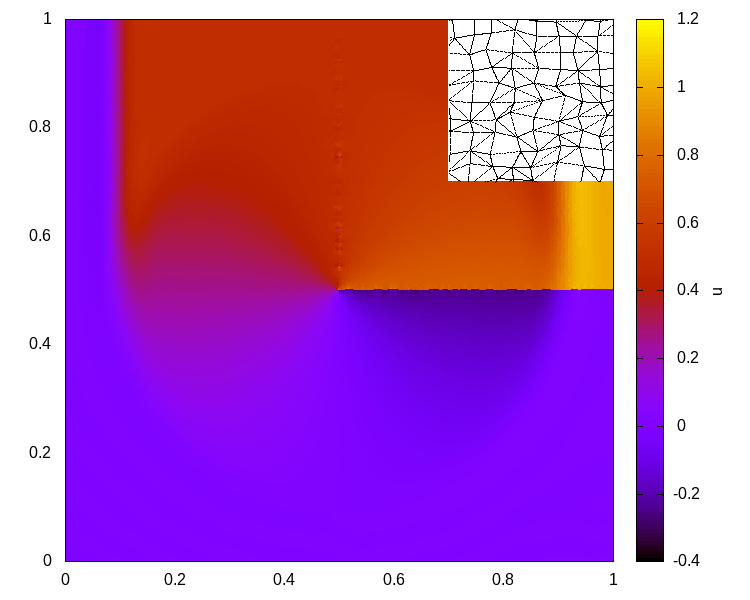}\\
 \includegraphics[width=0.45\textwidth]{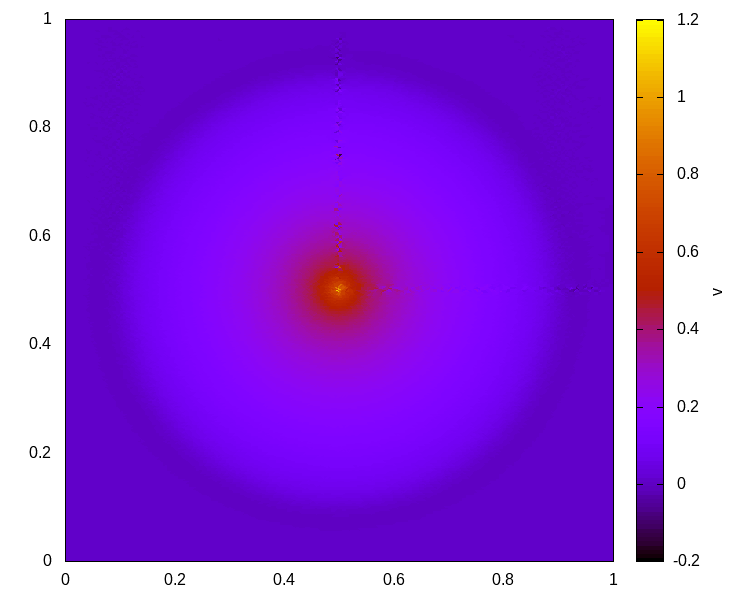}\hfill \includegraphics[width=0.45\textwidth]{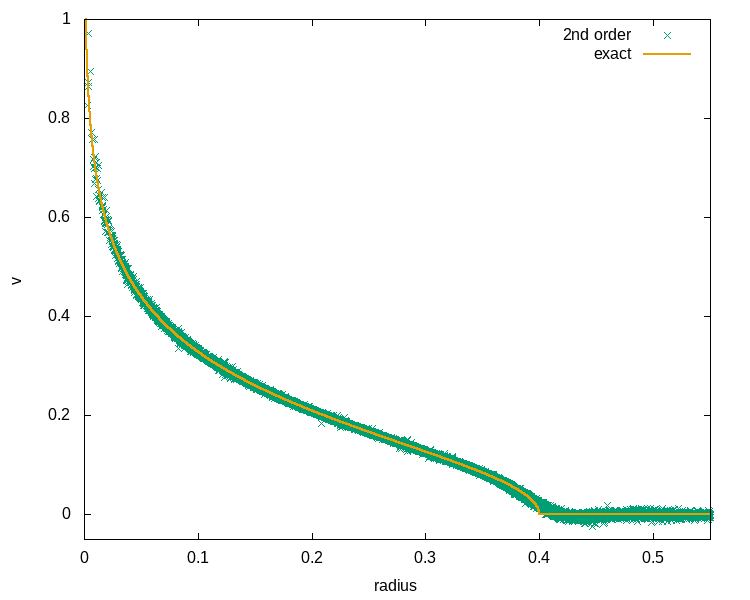}
 \caption{Same as Figure \ref{fig:rp}, but on a mixed triangular-quandrangular grid (a $600\%$ zoom of the region $[0.5, 0.55] \times [0.5,0.55]$ shown as inset) of $\sim 50000$ cells using the second order method with nodal pressure. The radial plot excludes data in $[0.49, 0.51] \times [0.51, 1]$ and $[0.51,1] \times [0.49, 0.51]$, i.e. around the lines of initial discontinuity where some artefacts of the initial discontinuity are visible in all plots.}
 \label{fig:rptri}
\end{figure}

Figure \ref{fig:rp} shows the numerical results at $t=0.4$ obtained on a grid of $200 \times 200$ cells using the solver with a nodal pressure. Zero-gradient boundary conditions were employed. The challenging singular solution is captured within the available resolution, with the 2nd order method being significantly more accurate. Figure \ref{fig:rptri} shows the same setup computed on a grid consisting of triangles and quadrangles using the second order method. The reconstruction involves all neighbours around the cell, also the ones only sharing a node.

\subsection{Spherical Riemann problem}

\begin{figure}
 \centering
 \includegraphics[width=0.45\textwidth]{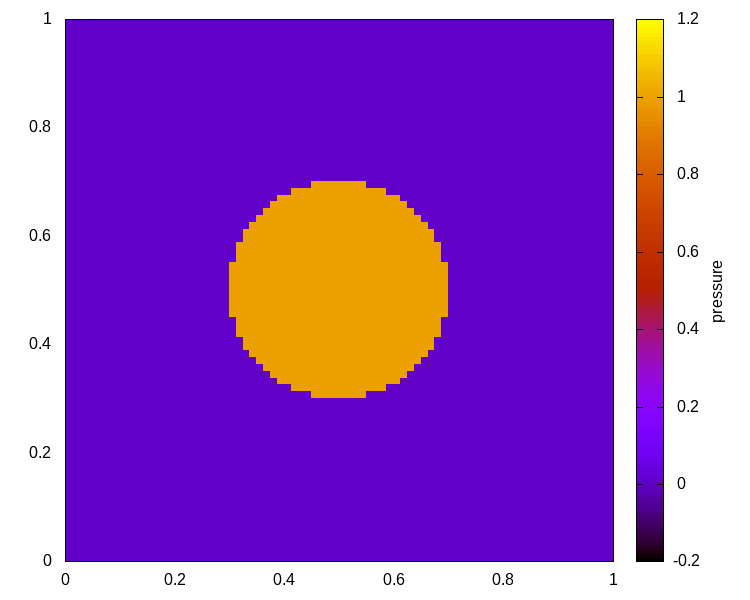}\hfill \includegraphics[width=0.45\textwidth]{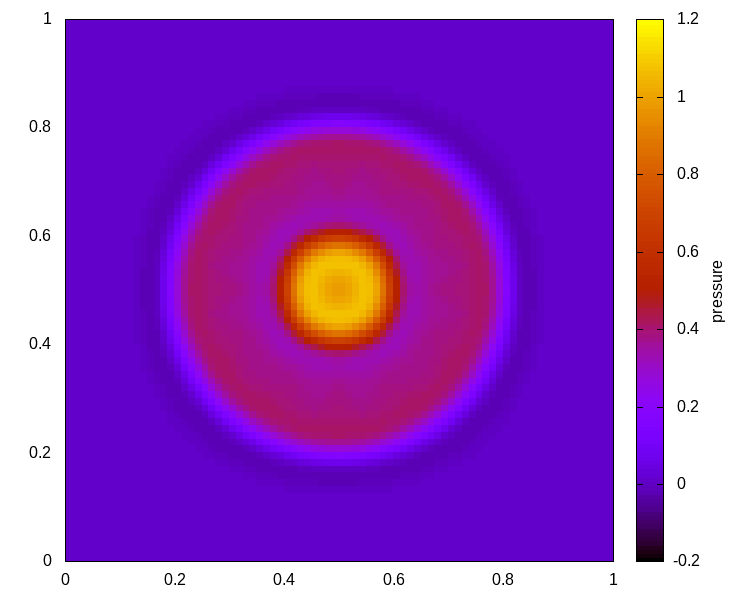}\\
 \includegraphics[width=0.45\textwidth]{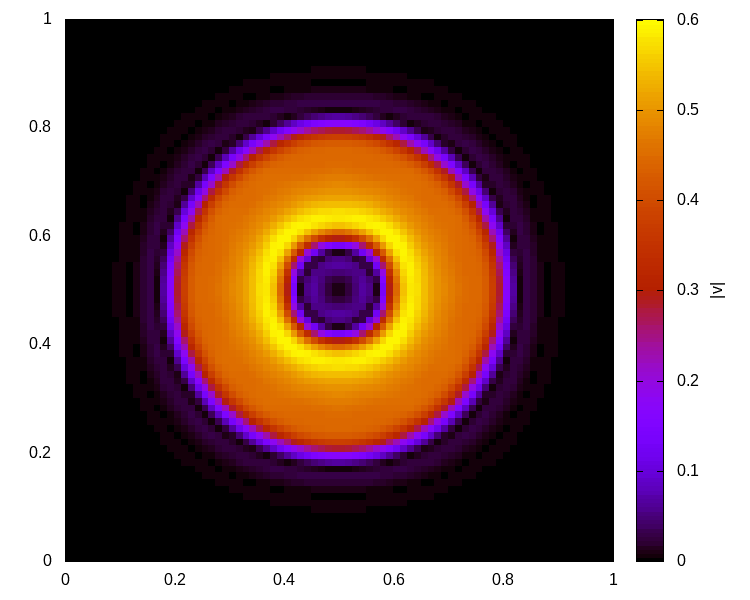}\hfill \includegraphics[width=0.45\textwidth]{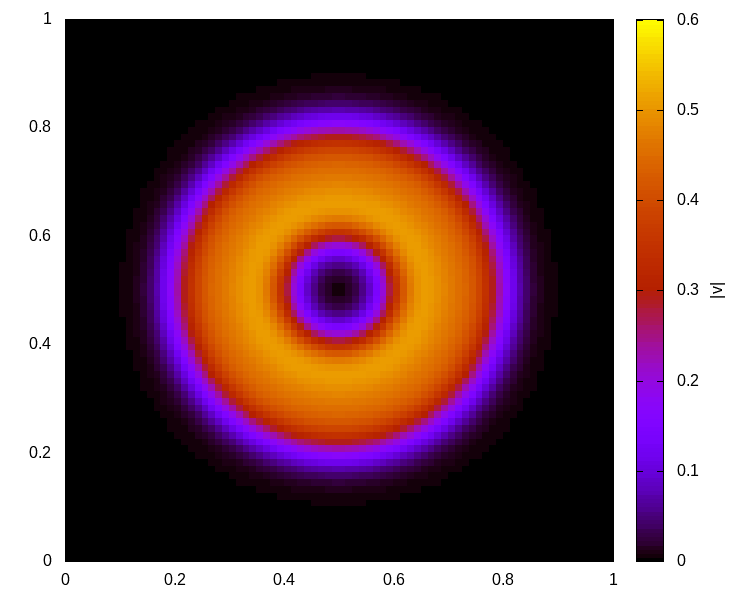}
 \caption{Spherical Riemann problem on a Cartesian grid. \emph{Top}: Pressure $p$ initially (\emph{left}) and at $t = 0.1$ (\emph{right}) for the 2nd order method. \emph{Bottom}: Comparison between the results for $|\vec v| = \sqrt{u^2 + v^2}$ for the 2nd order method (\emph{left}) and the 1st order method (\emph{right}).}
 \label{fig:shock}
\end{figure}

Figure \ref{fig:shock} shows the time evolution of an initial circular discontinuity between 1 and 0 in the pressure $p$ on a Cartesian grid of $80 \times 80$ cells. Zero-gradient boundary conditions were employed. This setup has initially zero vorticity because $u_0 = v_0 =0$. 

\begin{figure}
 \centering
 \includegraphics[width=0.45\textwidth]{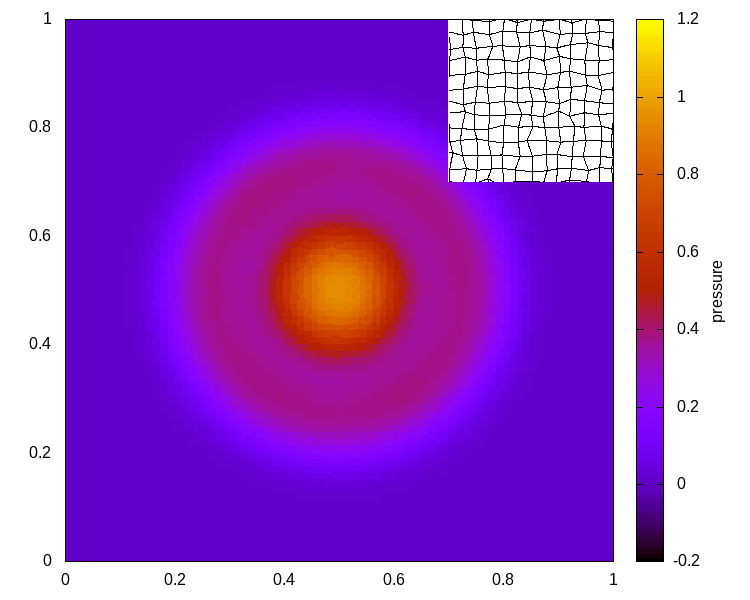}\hfill \includegraphics[width=0.45\textwidth]{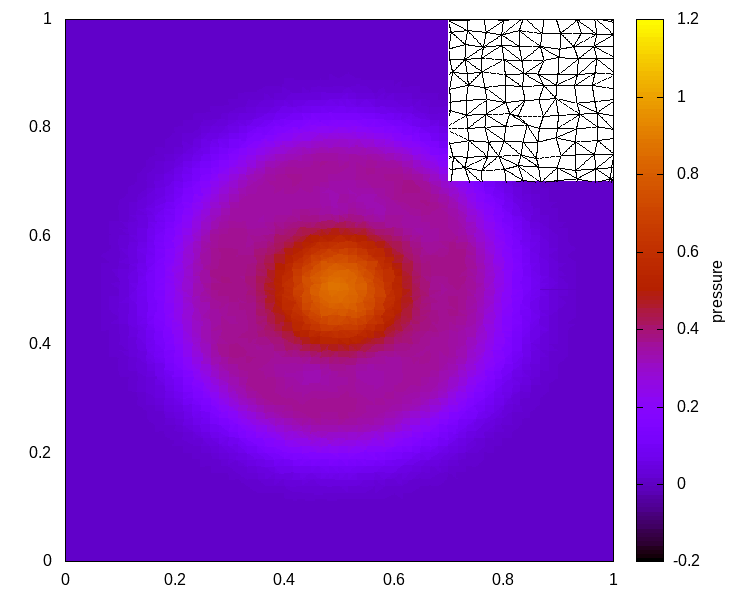}\\
 \includegraphics[width=0.45\textwidth]{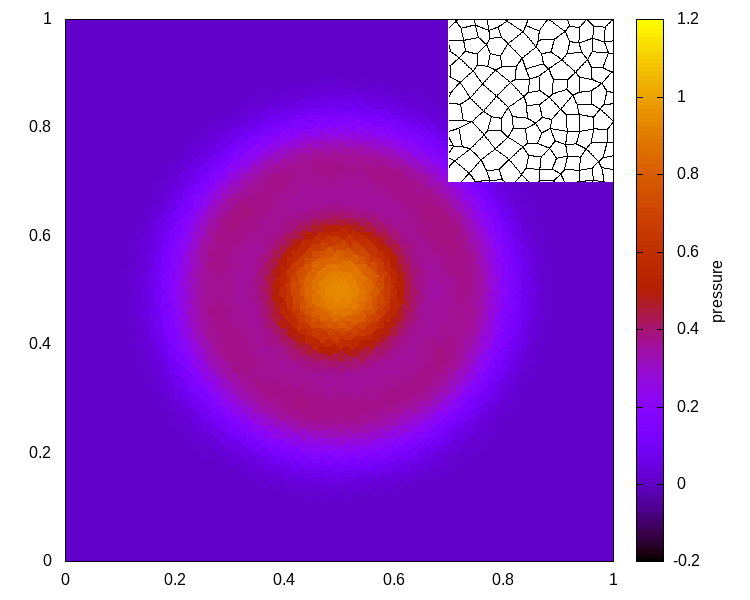} \hfill \includegraphics[width=0.45\textwidth]{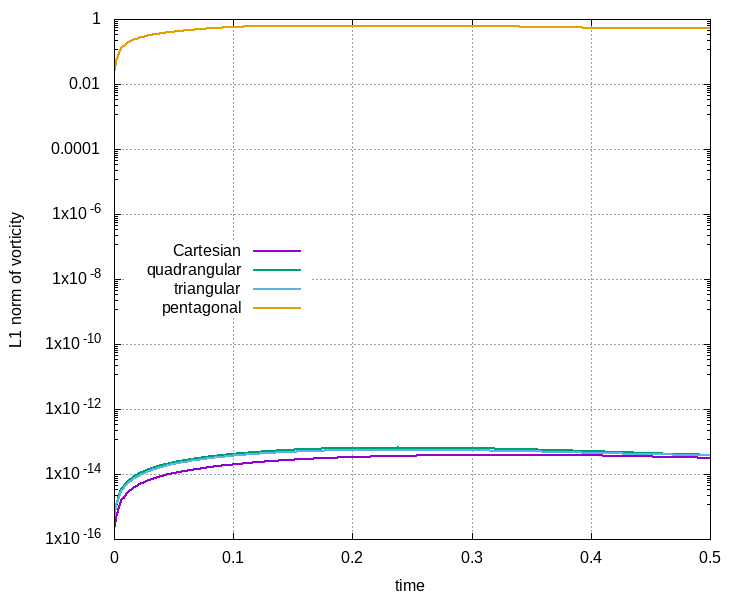}
 \caption{Spherical Riemann problem. Pressure $p$ at $t = 0.1$ for the 1st order method. \emph{Top left}: Perturbed quadrangular grid. \emph{Top right}: Mixed triangular-quadrangular grid. \emph{Bottom left}: Grid containing quadrangles, pentagons and hexagons. The insets show a 200\% zoom of the grid in $[0.5,0.65] \times [0.5,0.65]$. \emph{Bottom right}: $L^1$ norm of the vorticity \eqref{eq:unstructuredcurl} shown as a function of time. One observes that it is not preserved on the grid containing pentagons.}
 \label{fig:shockunstructured}
\end{figure}

\begin{figure}
 \centering
 \includegraphics[width=0.45\textwidth]{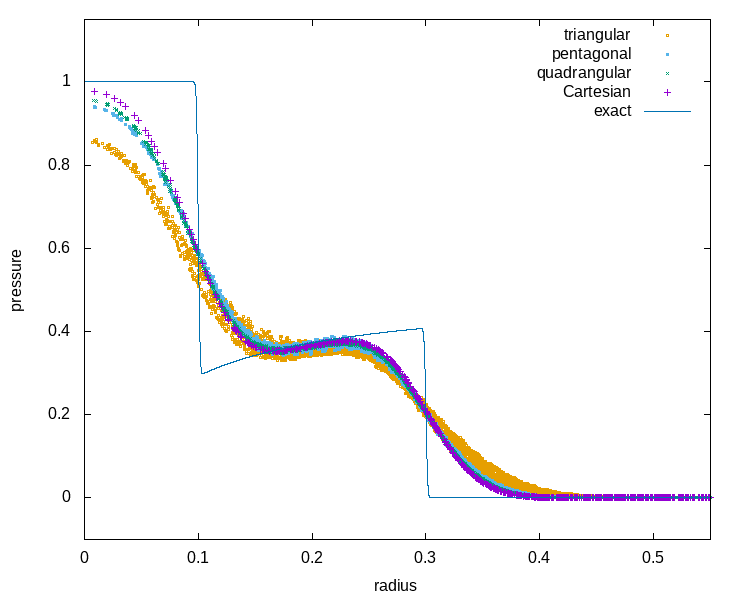}\hfill \includegraphics[width=0.45\textwidth]{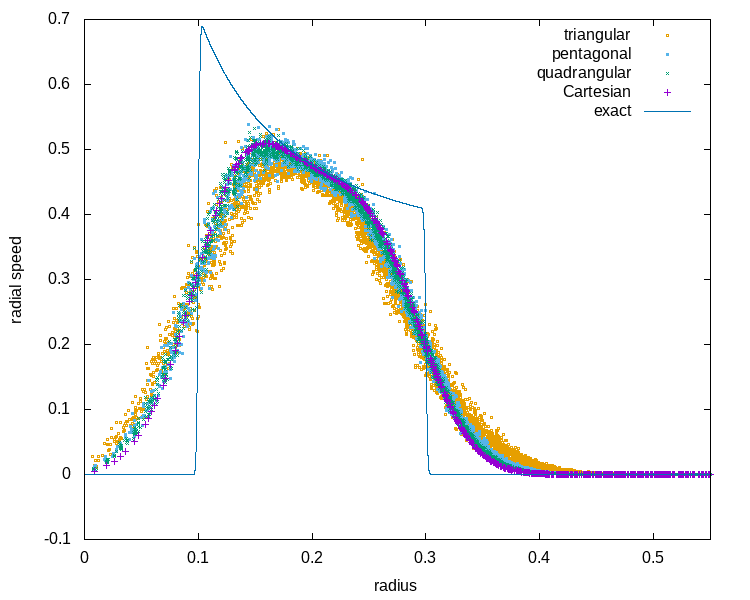}\\
 \includegraphics[width=0.45\textwidth]{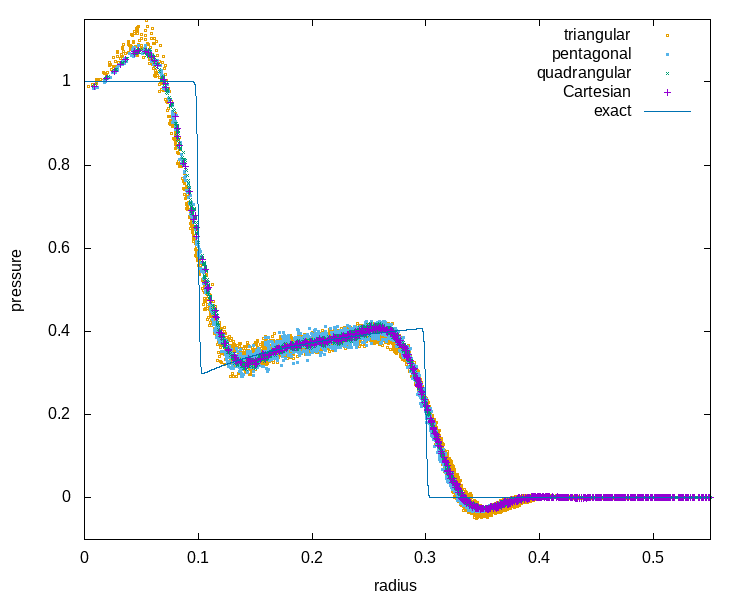}\hfill \includegraphics[width=0.45\textwidth]{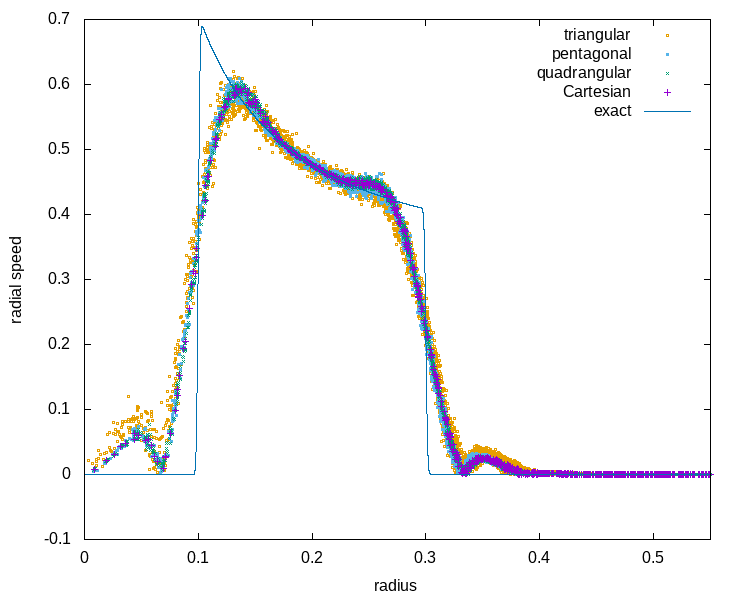}
 \caption{Spherical Riemann problem, radial scatter plots of $p$ (\emph{left}) and $|\vec v|$ (\emph{right}) for the first-order method on different grids (same as in Figure \ref{fig:shockunstructured}). \emph{Top}: First-order method. \emph{Bottom}: Second-order method.}
 \label{fig:shockradial}
\end{figure}

In Figure \ref{fig:shockunstructured} the same setup is shown solved using the 1st order method on a perturbed quadrangular, randomly mixed triangular-quadrangular and a grid with more complex polygons (in particular, pentagons), all of them having comparable resolutions to the Cartesian grid. The numerical solution captured on these grids is very similar to the one on the Cartesian grid. 
Figure \ref{fig:shockradial} shows radial scatter plots of the same setup.

\begin{figure}
 \centering
 \includegraphics[width=0.45\textwidth]{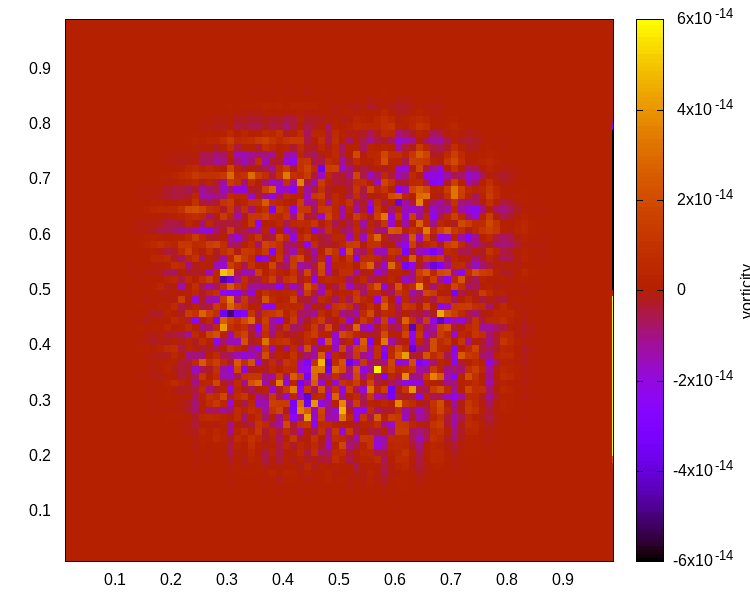}\hfill \includegraphics[width=0.45\textwidth]{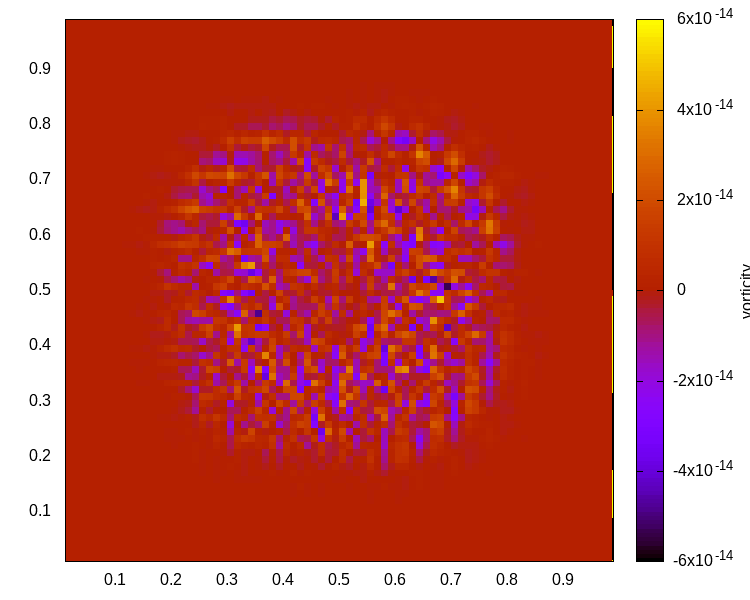}\\
 \includegraphics[width=0.45\textwidth]{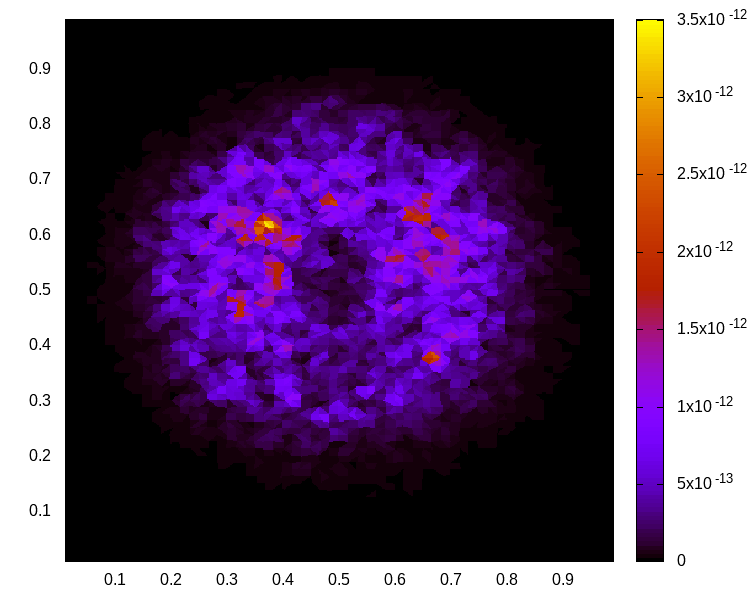}\hfill \includegraphics[width=0.45\textwidth]{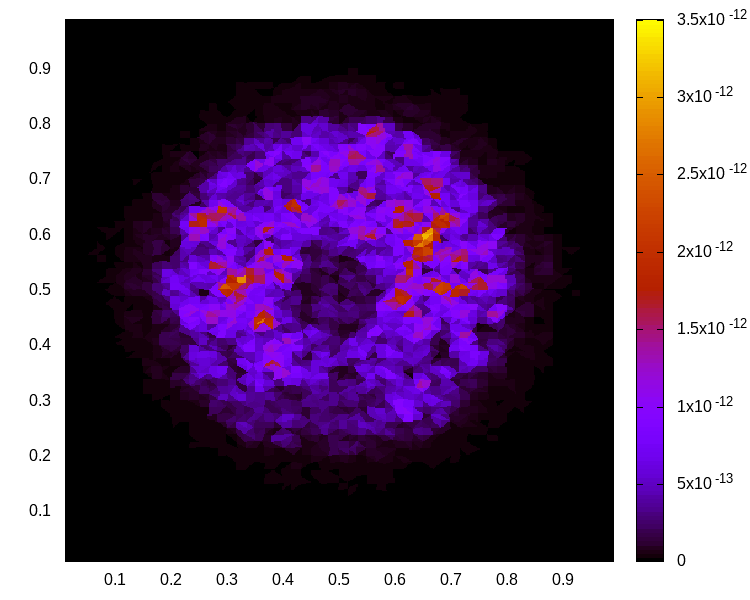}
 \caption{Spherical Riemann problem. Values of the discrete vorticity \eqref{eq:discretevorticity}/\eqref{eq:unstructuredcurl} at time $t=0.1$ for the 1st order method (\emph{left}) and the 2nd order method (\emph{right}) on a Cartesian grid (\emph{top}) and on a grid consisting of triangles and quadrangles (\emph{bottom}). For the latter, $\sum_{n \in \mathcal N(c)} |\vec C_n \vec v|$ in every cell $c$ is shown. They are at the level of the machine error.}
 \label{fig:shockomega}
\end{figure}

Figures \ref{fig:shockunstructured} and \ref{fig:shockomega} also show the discrete vorticity $\vec C_n \vec v$ (see Equation \eqref{eq:unstructuredcurl}) known to remain stationary for both the 1st and 2nd order methods for grids involving at most quadrangles. It is indeed of the order of the machine error (recall that the velocity $\vec v$ and thus also the vorticity are zero initially), which is not the case for the grid containing pentagons. 

\subsection{Stationary vortex}

\begin{figure}
 \centering
 \includegraphics[width=0.45\textwidth]{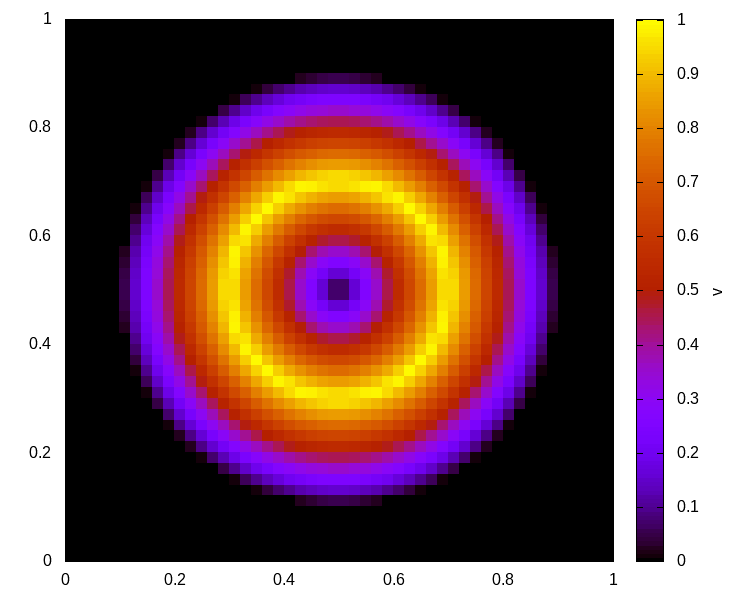}\hfill \includegraphics[width=0.45\textwidth]{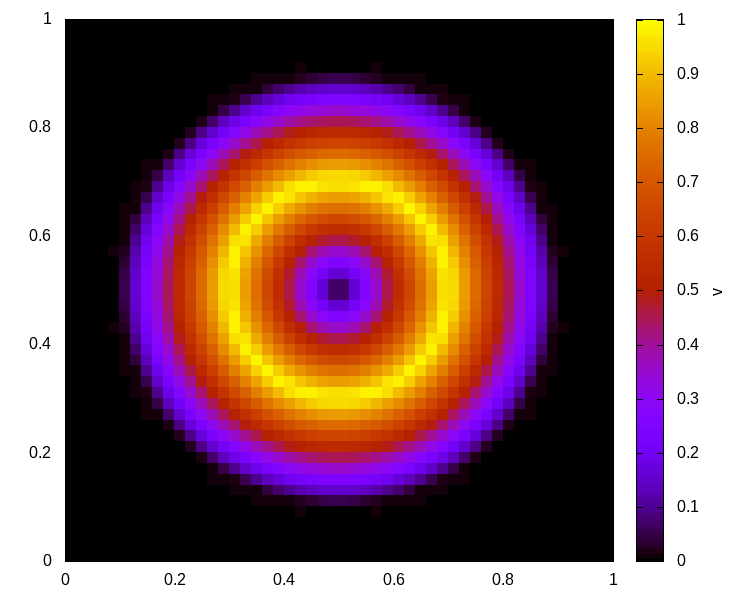}
 \includegraphics[width=0.45\textwidth]{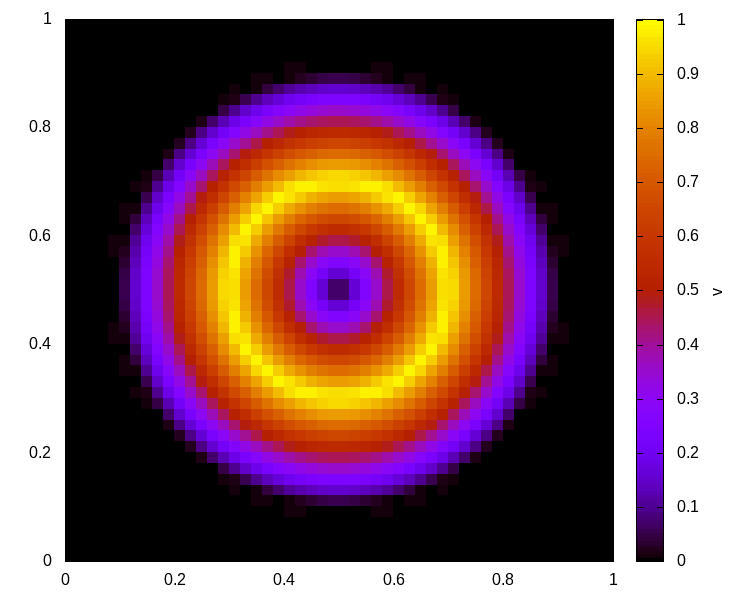} \hfill \includegraphics[width=0.45\textwidth]{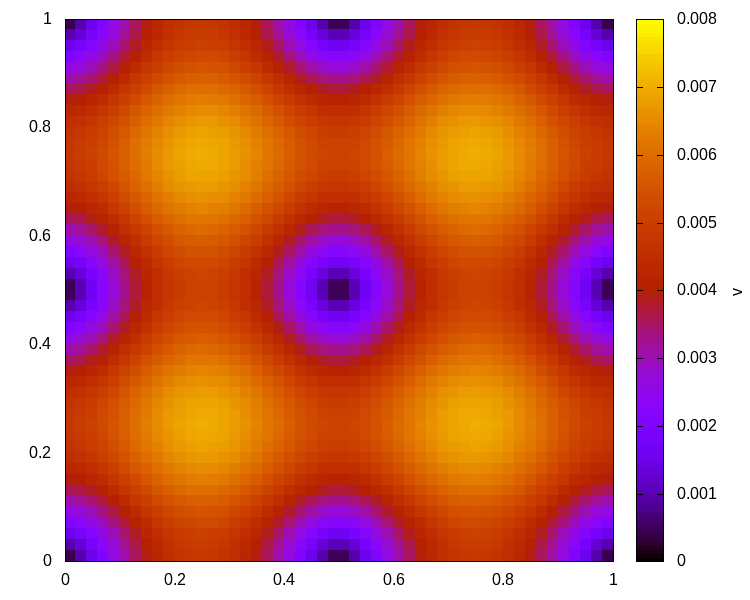}
 \caption{Results for the stationary vortex setup \eqref{eq:vortex} on a Cartesian grid. \emph{Top left}: Initial setup, $|\vec v| = \sqrt{u^2 + v^2}$ is shown. \emph{Top right}: Numerical solution at $t=100$ for the 1st order method. \emph{Bottom left}: The same at $t=100$ for the 2nd order method. \emph{Bottom right}: For comparison, the same setup at $t=10$ for the solver \eqref{eq:velocitysolver1}--\eqref{eq:velocitysolver2} employing a nodal velocity, i.e. a method that does not preserve any discrete vorticity. Observe the different scale of the colour bar. At $t=100$ the vortex has been diffused away to the level of machine error (not shown).}
 \label{fig:vortex}
\end{figure}

\begin{figure}
 \centering
 \includegraphics[width=0.45\textwidth]{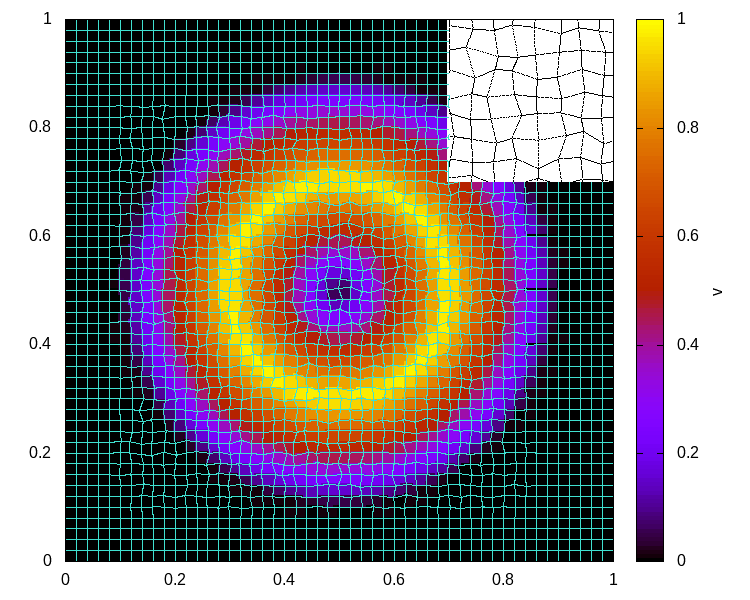}\hfill \includegraphics[width=0.45\textwidth]{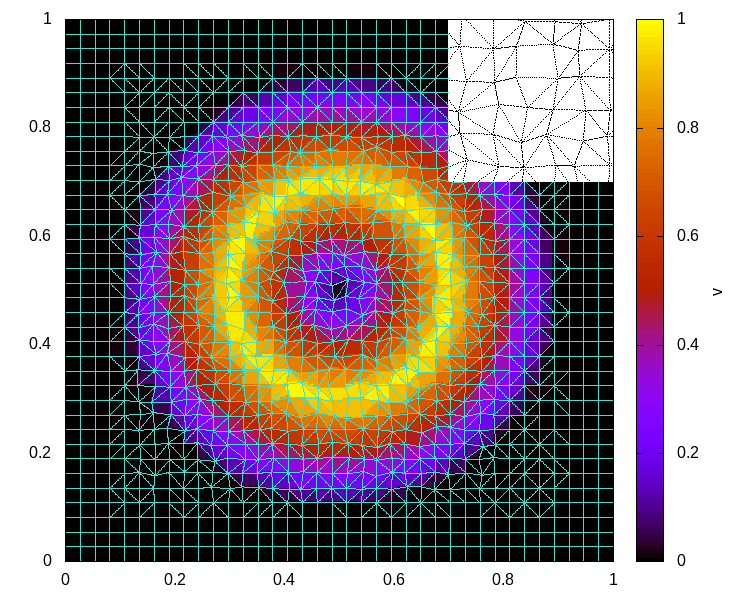}
 \includegraphics[width=0.45\textwidth]{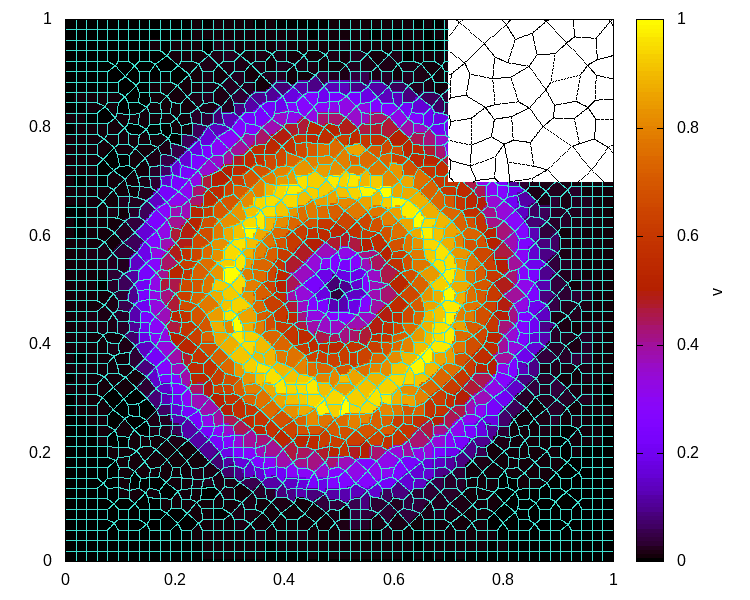} \hfill \includegraphics[width=0.45\textwidth]{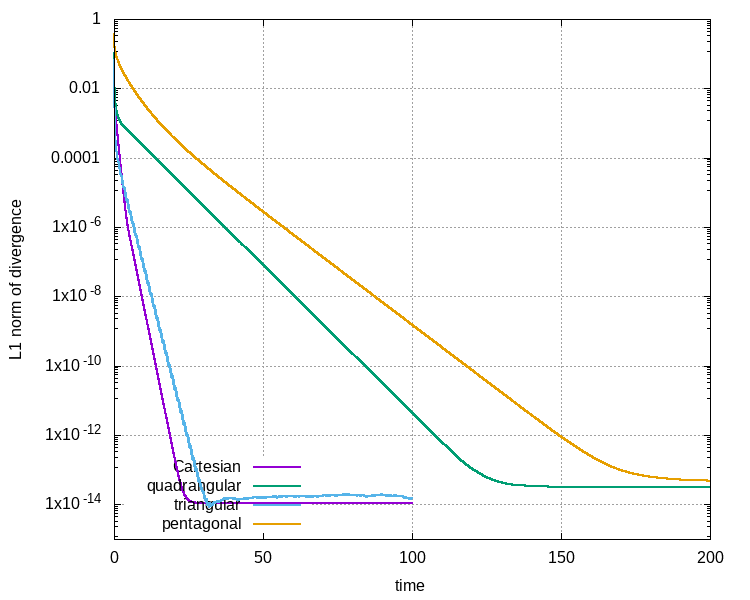}
 \caption{Results for the vortex setup \eqref{eq:vortex}. $|\vec v| = \sqrt{u^2 + v^2}$ is shown at $t=100$ for the 1st order method. \emph{Top left}: Quadrangular grid. \emph{Top right}: Mixed triangular-quadrangular grid. \emph{Bottom left}: Grid containing quadrangles, pentagons and hexagons. The insets show a 200\% zoom of the grid in $[0.5,0.65] \times [0.5,0.65]$. \emph{Bottom right}: $L^1$ norm of the divergence \eqref{eq:discretedivunstructured} shown as a function of time. }
 \label{fig:vortexunstructured}
\end{figure}

\begin{figure}
 \centering
 \includegraphics[width=0.45\textwidth]{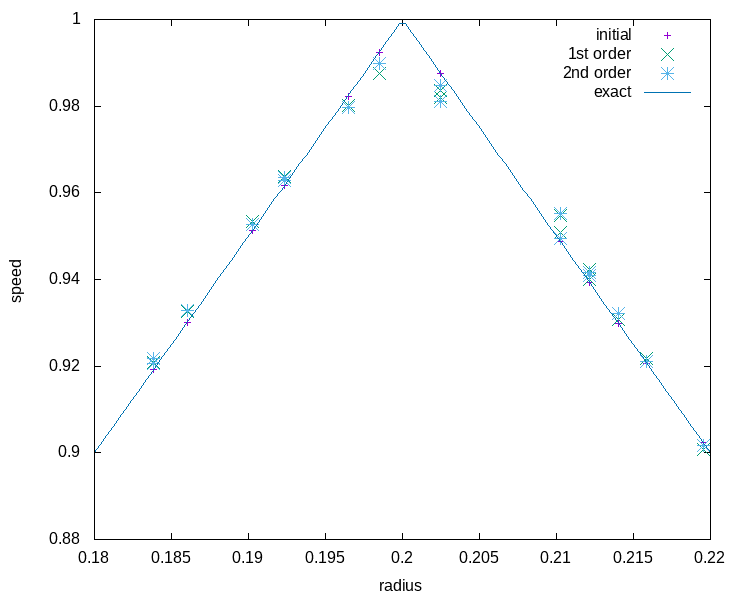} \includegraphics[width=0.45\textwidth]{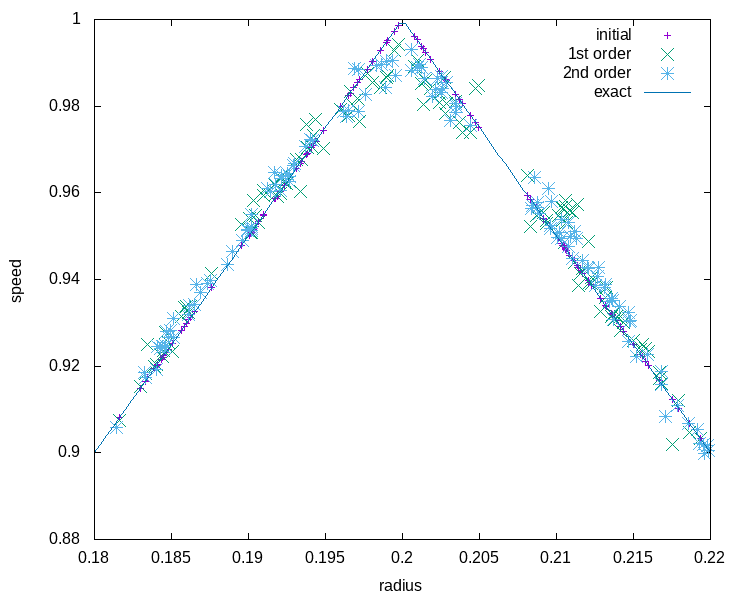}\\
 \includegraphics[width=0.45\textwidth]{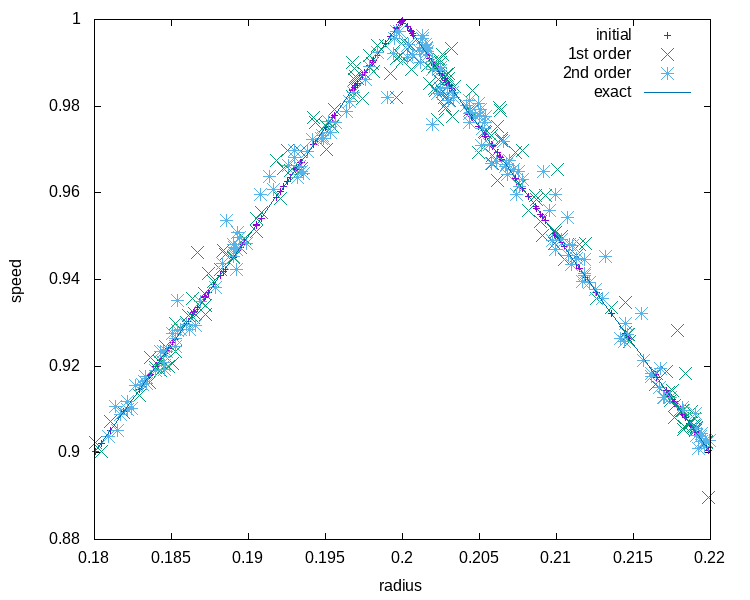} \includegraphics[width=0.45\textwidth]{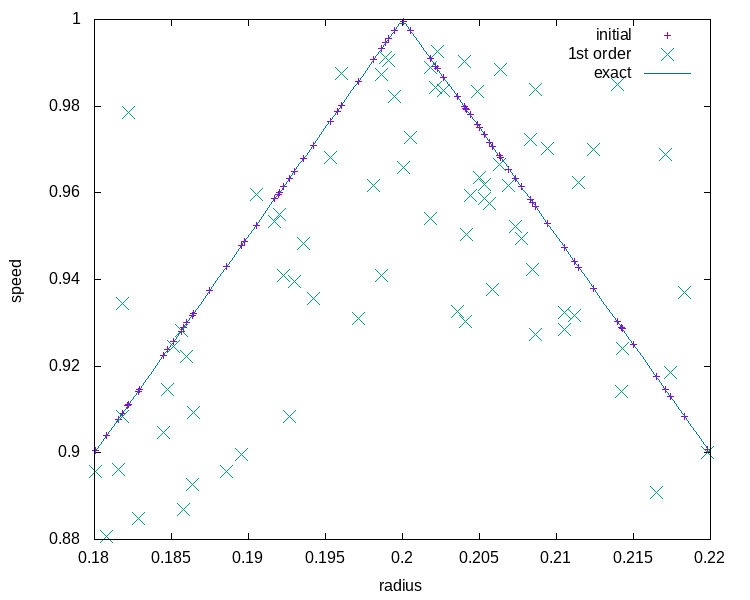}
 \caption{Radial scatter plots of $|\vec v|$ on different grids (observe the modified axes). \emph{Top left}: Cartesian grid. \emph{Top right}: Quadrangular grid. \emph{Bottom left}: Mixed triangular-quadrangular grid. \emph{Bottom right}: Grid containing quadrangles, pentagons and hexagons. The setup becomes numerically stationary, the solution no longer changes.}
 \label{fig:vortexradial}
\end{figure}

Finally, a stationary vortex setup will be used to assess the stationarity preservation property of the methods. The exact solution is given for all times by $p_0(\vec x) = 0$ and
\begin{align}
 \vec v_0(\vec x) &= v(|\vec x|) \vecc{-y/|\vec x|}{x/|\vec x|} ,& v(r) &= \begin{cases} r/w & r < w, \\ 2 - \frac{r}{w} & w \leq r < 2w, \\ 0 &\text{else.} \end{cases} \label{eq:vortex}
\end{align}
Figure \ref{fig:vortex} shows the initial setup with $w=0.2$ and the solutions at $t=100$ for the 1st and 2nd order methods on a Cartesian grid. No visible difference can be detected, as usual for stationarity preserving methods, compare e.g. to results in \cite{barsukow17a,barsukow18activeflux}. Figure \ref{fig:vortexunstructured} shows the same for different unstructured grids with comparable resolution. It also shows the evolution of the discrete divergence \eqref{eq:discretedivunstructured}. Even though these initial data originate from a stationary solution of the PDE, the discrete divergence is not zero initially, because we do not prepare the discrete initial data in any specific way. However, by stability arguments, all those modes that are evolving will eventually decay (exponentially). One can see in Figure \ref{fig:vortexradial} that after some time the numerical solution is no longer changing.

As time goes by, one also expects the discrete divergence to become machine-level, as is the case. Other discretizations of the divergence, in general, will become stationary (as the entire setup), but not reach machine-zero. A similar discussion with comparable numerical results can be found in \cite{barsukow20hypproceeding} for a different stationarity preserving method. The characterization of stationary states, of course, does not tell anything about how quickly they would be reached in the limit of long time. The different decay rates observable in Figure \ref{fig:vortexunstructured} therefore remain without explanation for the moment.

\section{Conclusions and outlook}

Systems of hyperbolic PDEs in multiple spatial dimensions have many more phenomena than their one-dimensional counterparts, and it has been recognized that numerical methods should be able to capture them. Prominent examples are vortices and non-trivial stationary states governed by a zero divergence. Both concepts are non-existent or trivial in 1-d. While Riemann solvers and Godunov methods are very successful in the one-dimensional setting, they are unable to capture many of those multi-dimensional properties. It has been shown in \cite{barsukow17} that this is not due to the approximate nature of many Riemann solvers used in practice, but occurs even if a Godunov method is used, and even if the evolution step in \emph{reconstruction-evolution-averaging} is computed exactly. Structure preserving finite volume methods at the moment typically rely on \textit{ad-hoc} modifications, such as the low Mach fixes for the Euler equations, and are not derived from any kind of first principle.

This paper presents the derivation of a method from the truly multi\-di\-men\-sio\-nal first principle of node-based conservation, appearing e.g. already in \cite{despres05,maire07}. Global conservation is usually made local by ensuring that the fluxes add up to zero at every edge. Nodal conservation reorganizes the global summation by halving every edge and making sure that the fluxes around a node sum up to zero. Such an approach does not have a one-dimensional counterpart, i.e. it is truly multi-dimensional. This requires the usage of non-standard Riemann solvers with free parameters that allow to enforce nodal conservation. For linear acoustics, this paper presents a new method for unstructured grids that is vorticity preserving. Its derivation highlights the fact that the choice of the free parameter in the Riemann solver is important, and only certain choices lead to structure preservation. We demonstrate the favourable properties of the new method theoretically and numerically.

Future work will be devoted to an extension of the method to three-dimensional elements and to finding the rules for the derivation of the non-standard Riemann solver in order to ensure structure preservation for more general systems of hyperbolic PDEs.

\newpage

\newcommand{\etalchar}[1]{$^{#1}$}

\appendix

\section{Summary of the methods}

\subsection{Unstructured grid}

The new method based on a nodal pressure $p^*$ (Section \ref{ssec:solverpressure} and \ref{ssec:nodalconservation}) reads

\begin{align}
  p^*_{n}  &= \frac{\displaystyle\sum_{s \in \mathcal{SE}(n)} |s| \left( \frac{p_\text{R} + p_\text{L}}{2} - \frac{(\vec v_\text{R}- \vec v_\text{L}) \cdot \vec n_{s}}{2} \right)}{\displaystyle\sum_{s \in \mathcal{SE}(n)} |s|} 
  = \frac{\displaystyle\sum_{s \in \mathcal{SE}(n)} |s|  \frac{p_\text{R} + p_\text{L}}{2}  -  \frac12 |c_n| D_n \vec v}{\displaystyle\sum_{s \in \mathcal{SE}(n)} |s|}
\end{align}

\begin{align}
 \frac{\dd}{\dd t} \vec v_c &= - \frac{1}{|c|} \sum_{n \in \mathcal N(c)} \ell_{nc} \vec n_{nc} p^*_n = - \vec G_c p^*\\
 \frac{\dd}{\dd t} p_c &=  - \frac{1}{|c|} \sum_{n \in \mathcal N(c)} \sum_{s \in \mathcal{SE}(n, c)} |s| (- p^*_n + p_c)
\end{align}
It is stationarity and vorticity preserving.

\subsection{Cartesian grid}

For the notation employed, see Section 6 in \cite{barsukow21yee}.

The method based on a nodal velocity $u^*, v^*$ on Cartesian grids (Sections \ref{ssec:solvervelocity} and \ref{ssec:nodalvelocitysolvercartesian}) reads

\begin{align}
 u^*_{i+\frac12,j+\frac12} &= \frac{\{ \{ u \}_{i+\frac12} \}_{j+\frac12}}{4}  - \frac12 \frac{\{ [p]_{i+\frac12} \}_{j+\frac12} }{2} \\
 v^*_{i+\frac12,j+\frac12} &= \frac{\{ \{ v \}_{i+\frac12} \}_{j+\frac12}}{4}  - \frac12 \frac{[\{ p \}_{i+\frac12}]_{j+\frac12} }{2}
\end{align}

\begin{align}
 \frac{\dd}{\dd t} u_{ij} &= \frac{- 4u_{ij} + \{ \{ u^*\}_{i\pm\frac12} \}_{j\pm\frac12} }{2\Delta x}
 = - \frac{\{\{ [p]_{i\pm1} \}\}_{j\pm\frac12} }{8 \Delta x} - \frac{ 2 \left( u_{ij} + \frac{\{\{ \{\{ u \}\}_{i\pm\frac12} \}\}_{j\pm\frac12}}{16} \right) }{\Delta x} \label{eq:velocitysolver1} \\
 \frac{\dd}{\dd t} v_{ij} &= \frac{- 4v_{ij} + \{ \{ v^*\}_{i\pm\frac12} \}_{j\pm\frac12} }{2\Delta x}\\
 \frac{\dd}{\dd t} p_{ij} &= - \frac{\{[u^*]_{i\pm\frac12}\}_{j\pm\frac12}}{2\Delta x} - \frac{[\{v^*\}_{i\pm\frac12} ]_{j\pm\frac12} }{2\Delta y}\\
 &= -\left(\frac{\{\{ [ u ]_{i\pm1} \}\}_{j+\frac12}}{8 \Delta x} + \frac{ [ \{\{ v \}\}_{i+\frac12}]_{j\pm1}}{8 \Delta y} \right ) + \frac12 \frac{\{\{ [[p]]_{i\pm\frac12} \}\}_{j\pm\frac12} }{4 \Delta x} +  \frac12 \frac{ [[ \{\{p \}\}_{i\pm\frac12}]]_{j\pm\frac12} }{4 \Delta y} \label{eq:velocitysolver2}
\end{align}
This solver is not vorticity preserving.

The new method based on a nodal pressure $p^*$ on Cartesian grids (Sections \ref{ssec:solverpressure} and \ref{ssec:nodalvelocitysolvercartesian}) reads

\begin{align}
 p^*_{i+\frac12,j+\frac12} &= \frac14 \{ \{ p \}_{i+\frac12} \}_{j+\frac12} - \frac12 \frac{\left( \frac{\{ [u]_{i+\frac12} \}_{j+\frac12}}{2 \Delta x} + \frac{[\{v \}_{i+\frac12} ]_{j+\frac12}}{2\Delta y} \right)}{\frac{1}{\Delta x} + \frac{1}{\Delta y}} 
\end{align}

\begin{align}
 \frac{\dd}{\dd t} u_{ij}  &= - \frac{ \{ [ p^* ]_{i\pm\frac12} \}_{j\pm\frac12}  }{2 \Delta x} 
 = - \frac{ \{\{ [ p ]_{i\pm1} \}\}_{j\pm\frac12} }{8 \Delta x} + \frac1{2 \Delta x} \frac{  \frac{\{\{ [[u]]_{i\pm\frac12} \}\}_{j\pm\frac12}}{4\Delta x} + \frac{[[v ]_{i\pm1} ]_{j\pm1}}{4\Delta y} }{\frac1{\Delta x} + \frac{1}{\Delta y}} \\
 \frac{\dd}{\dd t} v_{ij}  &= - \frac{ [\{  p^* \}_{i\pm\frac12} ]_{j\pm\frac12}  }{2 \Delta y}\\
 \frac{\dd}{\dd t} p_{ij}  &= \frac12 \left( \frac1{\Delta x} + \frac1{\Delta y}  \right )\left( \{\{ p^* \}_{i\pm\frac12} \}_{j\pm\frac12} - 4 p_{ij} \right) \\
 &= - \frac18 \left( \frac{\{\{ [u]_{i\pm1} \}\}_{j\pm\frac12}}{\Delta x} +\frac{ [\{\{v \}\}_{i\pm\frac12} ]_{j\pm1}}{\Delta y} \right) \\&\phantom{mmmmmmm}+\nonumber \frac12 \left( \frac1{\Delta x} + \frac1{\Delta y}  \right ) \left( \frac14 \{\{ \{\{ p \}\}_{i\pm\frac12} \}\}_{j\pm\frac12} - 4 p_{ij} \right )
\end{align}
This solver is vorticity preserving.

\section{The non-preservation of vorticity for the solver with velocities at nodes} \label{app:nonvortpres}

In Section \ref{sec:statioanritypreservation} it has been shown that the new method with pressures at nodes (derived in Sections \ref{ssec:solverpressure}/\ref{ssec:nodalconservation}) is stationarity and vorticity preserving, and an explicit discrete vorticity was constructed for grids with 3- and 4-sided cells, of which Cartesian grids are a special case. The proof of stationarity preservation was relatively easy, because it was obvious which discrete divergence would be involved. Contrarily, it was already much less obvious which discrete vorticity is preserved.

What remains to be shown is that no comparable result holds for the solver of Section \ref{ssec:solvervelocity}/\ref{ssec:methodgallice} with velocities at nodes. Here, one faces the conceptual difficulty that not having found a discrete divergence that is kept stationary might simply mean that one has to continue looking for it. Luckily, on Cartesian grids the much more powerful tool of the discrete Fourier transform is available, and its use in the context of structure preserving numerical methods was first advertised in \cite{barsukow17a}. With its help it is indeed possible to prove the non-existence of a discrete stationary divergence or, equivalently, the non-existence of a discrete vorticity preserved by the method.

One observation can be made immediately though: In \eqref{eq:nodalevelocitysolverveleq} the diffusion in the $u$-equation only involves $u$. Since the component $v$ is not involved, the modified equation cannot be \eqref{eq:multidmodifiedu}. This is not a proof yet, because the method might be achieving vorticity/stationarity preservation in some other way.

\subsection{A brief review of the continuous Fourier Transform}

A linear system $\del_t q + (\vec J \cdot \nabla) q = 0$ of PDEs is called hyperbolic if $\vec J \cdot \vec k$ is real diagonalizable $\forall \vec k$. Linear hyperbolic systems with $\det \vec J \cdot \vec k = 0$ for all $\vec k$ are of particular interest because they possess a rich set of stationary states, and also linear involutions. The stationary states are governed by
\begin{align}
 \vec J \cdot \nabla q_\text{stat} = 0.
\end{align}
$q_\text{stat} = \mathrm{const}$ is an obvious solution, but many more exist if $\det \vec J \cdot \vec k = 0$ for all $\vec k$, as is discussed later. Involutions are associated to a row-vector $r$ of linear operators such that
\begin{align}
 \del_t (r q) = - r (\vec J \cdot \nabla) q = 0 \qquad \forall q \label{eq:involutionrealspace}
\end{align}

It is easier to understand these properties in Fourier space, i.e. by assuming that the solution is of the form $\hat q \exp(\ii \vec k \cdot \vec x  - \ii \omega t)$ with $\ii$ the imaginary unit, $\hat q(\vec k) \in \mathbb R^m, \vec k \in \mathbb R^d, \omega \in \mathbb R$. It is a solution iff $\omega$ is an eigenvalue of $\vec J \cdot \vec k$. The constant states correspond to $\vec k = 0$, and will not be considered further in the analysis. By linearity, any two solutions (possibly with different $\vec k$) can be summed to give rise to a new solution.

The statement \eqref{eq:involutionrealspace} translates into
\begin{align}
 0 &= - \ii \hat r (\vec J \cdot \vec k) q \qquad \forall \vec k \in \mathbb R^d, \hat q \in \mathbb R^m,
\end{align}
for some row vector $\hat r$; the involution then reads $\del_t (\hat  r \hat q) = 0$. The involutions are thus associated to the left kernel of $\vec J \cdot \vec k$. If we were not to insist on being able to choose $\vec k$ freely, then the involution would not be valid for all solutions, but only for those with a specific spatial dependence.

To be a stationary state, $\hat q_\text{stat}$ has to be in the right kernel of $\vec J \cdot \vec k$:
\begin{align}
 \del_t \hat q_\text{stat} &= - \ii (\vec J \cdot \vec k) \hat q_\text{stat} = 0
\end{align}
To begin with, one does not need to assume that $\vec J \cdot \vec k$ has a non-trivial (right) kernel for all $\vec k$. This is why \cite{barsukow17a} distinguishes trivial and non-trivial stationary states. A linear system of PDEs is said to possess non-trivial stationary states if for \emph{every} $\vec k$ one can find a $q_\text{stat}\neq 0$ in the right kernel of $\vec J \cdot \vec k$. This does not mean that there are no conditions on the spatial dependence of stationary states; indeed such a statement would imply that any initial data are stationary. $q_\text{stat}$ in general depends on $\vec k$. Broadly speaking, non-trivial stationary states are such that at least one component of $q_\text{stat}$ can be chosen freely as a function of $\vec x$. This is not the case for equations which have only trivial stationary states. An example of the latter is linear advection $\del_ t q + U \del_x q + V \del_y q = 0$. Its $1\times 1$ matrix $\vec J \cdot \vec k$ reads $U k_x + V k_y$ and vanishes if either $\vec k = 0$, or if $\vec k \perp (U,V)^\text{T}$. This means that the stationary state can only have spatial variation perpendicular to the velocity vector.

\subsection{A brief review of the discrete Fourier transform for finite differences}

The discrete Fourier transform allows to rephrase a numerical method in the same language. Because the shift operator $T$ which maps $q_i \mapsto q_{i+1}$ is diagonal on functions $\exp(\ii k_x \Delta x i)$ with eigenvalue $\exp(\ii k_x \Delta x)$, the discrete Fourier transform allows to replace finite differences by algebraic factors. Define therefore
\begin{align}
 t_x &:= \exp(\ii k_x \Delta x) , & t_y &:= \exp(\ii k_y \Delta y),
\end{align}
and by writing $q_{ij} = \hat q(t) t_x t_y$ observe that
\begin{align}
 q_{i + a,j+b} = q_{ij} t_x^a t_y^b.
\end{align}
Any linear method can therefore be written as 
\begin{align}
 \del_t \hat q(t) + \mathcal E \hat q(t) &= 0,
\end{align}
with the \emph{evolution matrix} $\mathcal E$ a complex-valued $m \times m$ matrix which depends on $t_x, t_y$, and thus on $\vec k$. It is the discrete counterpart of $\vec J \cdot \vec k$. The numerical method is called involution/stationarity preserving if $\dim \ker \mathcal E = \dim \ker \vec J \cdot \vec k$ for all $\vec k$, with discrete involutions being associated to the left kernel of $\mathcal E$ and discrete non-trivial stationary states to its right kernel. Thus, involution preserving numerical methods are stationarity preserving; the two properties are equivalent, although in practice sometimes one is easier to show than the other. For more details, see \cite{barsukow17a}.

The evolution matrix for the solver of Section \ref{ssec:solvervelocity} is
\begin{align}
 \mathcal E = \left ( \begin{array}{ccc} \frac{-4+\frac{(t_x+1)^2 (t_y+1)^2}{4 t_x t_y}}{2 \Delta x} & 0 & -\frac{(t_x-1) (t_x+1)
(t_y+1)^2}{8 \Delta x t_x t_y} \\ 0 & \frac{-4+\frac{(t_x+1)^2 (t_y+1)^2}{4 t_x t_y}}{2 \Delta y} & -\frac{(t_x+1)^2 (t_y-1) (t_y+1)}{8 \Delta y t_x t_y} \\ -\frac{(t_x-1) (t_x+1) (t_y+1)^2}{8 \Delta x t_x
t_y} & -\frac{(t_x+1)^2 (t_y-1) (t_y+1)}{8 \Delta y t_x t_y} & \frac{(t_x+1)^2 (t_y-1)^2}{8 \Delta y
t_x t_y}+\frac{(t_x-1)^2 (t_y+1)^2}{8 \Delta x t_x t_y} \end{array} \right )
\end{align}
and its determinant does not vanish:
\begin{align}
 \det \mathcal E &= -\frac{\left(\Delta x (t_x+1)^2 (t_y-1)^2+\Delta y (t_x-1)^2 (t_y+1)^2\right) }{32 \Delta x^2 \Delta y^2 t_x^2 t_y^2}  \times \\& \nonumber \phantom{mmmmmmmmmm}\Big((t_y+1)^2+t_x^2 (t_y+1)^2+2 t_x (1+(-6+t_y) t_y)\Big)
\end{align}
The method thus only has a non-trivial kernel for certain $\vec k$ -- it only possesses trivial stationary states and no discrete vorticity.

Equations \eqref{eq:uevolution}--\eqref{eq:vevolution} (solver of Section \ref{ssec:solverpressure} with $p^*$), upon the discrete Fourier transform, read
\begin{align}
\Delta x \del_t \hat u  &= - \frac12 (t_y+1)(t_x-1) \hat p^*,\\
\Delta y \del_t \hat v  &= - \frac12 (t_y-1)(t_x+1) \hat p^*. 
\end{align}
Even without taking into account the equation for $p$ and a detailed study of the evolution matrix, one immediately observes that
\begin{align}
 \del_t \left(\frac{t_y-1}{\Delta y} \frac{t_x+1}{2}  \hat u - \frac{t_x-1}{\Delta x} \frac{t_y+1}{2} \hat v \right)&= 0 .
\end{align}
This is the Fourier transform of the following discrete involution:
\begin{align}
 \del_t \left( \frac{[\{ u \}_{i+\frac12}]_{j+\frac12}}{2\Delta y}  -  \frac{\{ [ v]_{i+\frac12}\}_{j+\frac12}}{2\Delta x}  \right) = 0 
\end{align}

Again, even without constructing $\mathcal E$ explicitly, one can easily see that the stationary states of the first-order method are given by a vanishing discrete (node-centered) divergence $\mathscr D_{i+\frac12,j+\frac12} = \frac{\{ [ u ]_{i+\frac12}\}_{j+\frac12}}{\Delta x}  +  \frac{[ \{ v\}_{i+\frac12}]_{j+\frac12}}{\Delta y}$ and constant pressure. This is because the divergence appearing in the evolution equation \eqref{eq:pstarsolverevoutionofp} of $p$ is the average $\{\{ \mathscr D \}_{i\pm\frac12} \}_{j\pm\frac12}$, while the diffusion appearing in \eqref{eq:pstarsolverevoutionofu} is the discrete derivative $\{[ \mathscr D ]_{i\pm\frac12} \}_{j\pm\frac12}$. In fact, the method is very close to the multi-dimensional vorticity/stationarity preserving method from \cite{barsukow17a}, where this discrete divergence has been shown to be the only one allowing for such a construction if one restricts oneself to 9-point stencils and (anti)symmetric operators.

\end{document}